\documentclass[11pt]{amsart}

\usepackage[
paper=a4paper,
text={147mm,230mm},centering
]{geometry}

\iftrue
\makeatletter
\def\@settitle{%
  \vspace*{-20pt}
  \begin{flushleft}%
    \baselineskip14\p@\relax
    \normalfont\bfseries\LARGE
    \@title
  \end{flushleft}%
}
\def\@setauthors{%
  \begingroup
  \def\thanks{\protect\thanks@warning}%
  \trivlist
  \large \@topsep30\p@\relax
  \advance\@topsep by -\baselineskip
  \item\relax
  \author@andify\authors
  \def\\{\protect\linebreak}%
  \authors
  \ifx\@empty\contribs
  \else
    ,\penalty-3 \space \@setcontribs
    \@closetoccontribs
  \fi
  \normalfont
  \@setaddresses
  \endtrivlist
  \endgroup
}
\def\@setaddresses{\par
  \nobreak \begingroup\raggedright
  \small
  \def\author##1{\nobreak\addvspace\smallskipamount}%
  \def\\{\unskip, \ignorespaces}%
  \interlinepenalty\@M
  \def\address##1##2{\begingroup
    \par\addvspace\bigskipamount\noindent
    \@ifnotempty{##1}{(\ignorespaces##1\unskip) }%
    {\ignorespaces##2}\par\endgroup}%
  \def\curraddr##1##2{\begingroup
    \@ifnotempty{##2}{\nobreak\noindent\curraddrname
      \@ifnotempty{##1}{, \ignorespaces##1\unskip}\/:\space
      ##2\par}\endgroup}%
  \def\email##1##2{\begingroup
    \@ifnotempty{##2}{\smallskip\nobreak\noindent E-mail address%
      \@ifnotempty{##1}{, \ignorespaces##1\unskip}\/:\space
      \ttfamily##2\par}\endgroup}%
  \def\urladdr##1##2{\begingroup
    \def~{\char`\~}%
    \@ifnotempty{##2}{\nobreak\noindent\urladdrname
      \@ifnotempty{##1}{, \ignorespaces##1\unskip}\/:\space
      \ttfamily##2\par}\endgroup}%
  \addresses
  \endgroup
  \global\let\addresses=\@empty
}
\def\@setabstracta{%
    \ifvoid\abstractbox
  \else
    \skip@25\p@ \advance\skip@-\lastskip
    \advance\skip@-\baselineskip \vskip\skip@
    \box\abstractbox
    \prevdepth\z@ 
    \vskip-10pt
  \fi
}
\renewenvironment{abstract}{%
  \ifx\maketitle\relax
    \ClassWarning{\@classname}{Abstract should precede
      \protect\maketitle\space in AMS document classes; reported}%
  \fi
  \global\setbox\abstractbox=\vtop \bgroup
    \normalfont\small
    \list{}{\labelwidth\z@
      \leftmargin0pc \rightmargin\leftmargin
      \listparindent\normalparindent \itemindent\z@
      \parsep\z@ \@plus\p@
      
    }%
    \item[\hskip\labelsep\bfseries\abstractname.]%
}{%
  \endlist\egroup
  \ifx\@setabstract\relax \@setabstracta \fi
}
\def\section{\@startsection{section}{1}%
  \z@{-1.2\linespacing\@plus-.5\linespacing}{.8\linespacing}%
  {\normalfont\bfseries\large}}
\def\subsection{\@startsection{subsection}{2}%
  \z@{-.8\linespacing\@plus-.3\linespacing}{.3\linespacing\@plus.2\linespacing}%
  {\normalfont\bfseries}}
\def\subsubsection{\@startsection{subsubsection}{3}%
  \z@{.7\linespacing\@plus.1\linespacing}{-1.5ex}%
  {\normalfont\itshape}}
\def\@secnumfont{\bfseries}
\makeatother
\fi 

\usepackage{amssymb}
\usepackage{mathrsfs}
\usepackage[all]{xy}
\usepackage{graphicx,color,float}
\usepackage[bookmarks]{hyperref}
\usepackage{pinlabel}
\usepackage[textwidth=1.3in,color=yellow]{todonotes}
\usepackage{amsmath}
\usepackage{pb-diagram,pb-xy}
\usetikzlibrary{calc}
\textwidth=\paperwidth \advance\textwidth by-2\oddsidemargin
\advance\oddsidemargin by-1in
\evensidemargin=\oddsidemargin
\calclayout

\theoremstyle{plain}
\newtheorem{theorem}{Theorem}[section]
 
\newtheorem{thmx}{Theorem}

\newtheorem*{lem}{Lemma}

\newtheorem{proposition}[theorem]{Proposition}
\newtheorem{lemma}[theorem]{Lemma}
\newtheorem{corollary}[theorem]{Corollary}

\theoremstyle{definition}

\newtheorem{definition}[theorem]{Definition}
\newtheorem{example}[theorem]{Example}

\theoremstyle{remark}
\newtheorem{remark}[theorem]{Remark}

\renewcommand{\bar}{\overline}
\renewcommand{\tilde}{\widetilde}

\renewcommand{\hat}{\widehat}

\newcommand{\C}{\mathbb{C}}

\newcommand{\Q}{\mathbb{Q}}
\newcommand{\bL}{\mathbb{L}}
\newcommand{\BL}{\mathbb{L}}
\newcommand{\cS}{\mathcal{S}}
\newcommand{\bS}{\mathbb{S}}
\newcommand{\Gr}{\mathrm{Gr}}

\newcommand{\pt}{\mathrm{pt}}

\newcommand{\R}{\mathbb{R}}

\newcommand{\Z}{\mathbb{Z}}

\newcommand{\bP}{\mathbb{P}}

\newcommand{\CP}{\mathbb{C}P}

\newcommand{\one}{\mathbf{1}}
\newcommand{\Fuk}{\mathrm{Fuk}}
\newcommand{\OG}{\mathrm{OG}}
\newcommand{\MF}{\mathrm{MF}}

\def\pa{\partial}
\def\mcal{\mathcal}
\def\frak{\mathfrak}
\def\scr{\mathscr}

\numberwithin{equation}{section} \numberwithin{table}{section}

\def\to{\mathchoice{\longrightarrow}{\rightarrow}{\rightarrow}{\rightarrow}}
\makeatletter
\newcommand{\shortxra}[2][]{\ext@arrow 0359\rightarrowfill@{#1}{#2}}
\def\longrightarrowfill@{\arrowfill@\relbar\relbar\longrightarrow}
\newcommand{\longxra}[2][]{\ext@arrow 0359\longrightarrowfill@{#1}{#2}}

\makeatother
\numberwithin{equation}{section}

\begin{document}                                                                          
\title[Immersed two-spheres and SYZ with Application to Grassmannians]{Immersed two-spheres and SYZ \\with Application to Grassmannians}

\author{Hansol Hong}
\address{Department of Mathematics, Yonsei University}
\email{hansolhong@yonsei.ac.kr}

\author{Yoosik Kim}
\address{ Department of Mathematics, Brandeis University \hfill \break
Center of Mathematical Sciences and Applications, Harvard University} 
\email{yoosik@brandeis.edu, yoosik@cmsa.fas.harvard.edu}

\author{Siu-Cheong Lau}
\address{Department of Mathematics and Statistics, Boston University}
\email{lau@math.bu.edu}


\begin{abstract}
	We develop a Floer theoretical gluing technique and apply it to deal with the most generic singular fiber in the SYZ program, namely the product of a torus with the immersed two-sphere with a single nodal self-intersection.  As an application, we construct immersed Lagrangians in $\mathrm{Gr}(2,\mathbb{C}^n)$ and $\mathrm{OG}(1,\mathbb{C}^5)$ and derive their SYZ mirrors.  It recovers the Lie theoretical mirrors constructed by Rietsch.  It also gives an effective way to compute stable disks (with non-trivial obstructions) bounded by immersed Lagrangians.
\end{abstract}


\maketitle
\setcounter{tocdepth}{1} 
\tableofcontents

\section{Introduction}

Strominger-Yau-Zaslow (SYZ) \cite{SYZ} proposed in a general context that \emph{mirror symmetry is T-duality}.  It conjectures a geometric way to construct mirrors and derive homological mirror symmetry via duality of special Lagrangian torus fibrations.  

Singular fibers in a Lagrangian fibration are the main difficulty in realizing the SYZ program. They lead to quantum corrections and the wall-crossing phenomena.  In the previous literatures concerning the SYZ construction via symplectic geometry, smooth torus fibers were mainly dealt with.  As a result, the precise Floer theoretical mirrors are incomplete since they have missing points in codimension two corresponding to the singular fibers.  In particular, they are not sufficient for homological mirror symmetry, since the missing strata include non-trivial objects in the derived category in most of situations.

The purpose of this paper is twofold.  First, we develop a Floer-theoretical gluing technique and fill in the missing points coming from the most generic singular SYZ fibers, which are the product of a torus and the immersed sphere with one nodal self-intersection.  They are closely related to the wall-crossing formula.

Second, we apply the technique to construct the complete SYZ mirrors for the type-A and type-B Grassmannians $\Gr(2,\C^n)$ and $\OG(1,\C^5)$.  It agrees with the Lie theoretical mirror of Rietsch, and in particular provides a method to extract open Gromov-Witten potential of Lagrangian branes from the Rietsch mirror, see Remark \ref{rmk:Rie} and Example \ref{ex-intro}.  This is important for studying symplectic geometry, particularly the Floer cohomology and non-displaceability of Lagrangian branes.

The family Floer theory proposed by Fukaya \cite{Fukaya-famFl} and developed by Tu \cite{Tu-reconstruction} and Abouzaid \cite{Ab-famFl} provides a powerful way to construct the SYZ mirrors.  They utilize the deformation theory of the SYZ torus fibers to construct the mirror charts and glue them up to a Floer-theoretical mirror.

In this paper, we consider the monotone torus fiber of the Gelfand-Cetlin Lagrangian fibration \cite{GS83}, together with finitely many Lagrangian immersions, which are immersed fibers of certain Lagrangian fibrations which interpolate different toric degenerations.  Note that the Lagrangian immersions that we use are \emph{not} fibers of the Gelfand-Cetlin system.  Our method has the following advantages.

First, the singular fibers of the Gelfand-Cetlin systems only occur at the boundaries of the base polytopes.  For instance, in $\Gr(2,\C^4)$, there is a Lagrangian $\mathrm{U}(2)$ sitting over a boundary edge of the Gelfand-Cetlin base polytope, which is a non-zero object in the Fukaya category over the field of characteristic zero by the work of Nohara-Ueda \cite{NU16}.  It is a generator for a summand of the quantum cohomology ring by the work of Evans-Lekili \cite{EL19}.  If we use the family Floer theory of the fibers over the interior to construct an LG mirror, then there are missing critical points corresponding to these non-trivial objects over the boundary.  Unfortunately, these singular fibers have smaller deformation spaces which do not glue well with this LG mirror.  Thus we shall not merely work with fibers of a Gelfand-Cetlin system.

Second, for the family Floer theory of Lagrangian torus fibers of one fibration, infinitely many mirror charts coming from the torus fibers are required due to the scattering phenomenon \cite{KS-affine, GS07}.  On the other hand, our method only needs to employ finitely many mirror charts coming from the immersed Lagrangians.  Most of the mirror charts coming from the torus fibers are indeed redundant and are already covered by the charts of immersed Lagrangians.  This method is particularly efficient and practical for constructing $\C$-valued mirrors of Fano and general-type manifolds.

Below we explain our method and the main theorems.

\subsection{Immersed SYZ fibers and wall-crossing}
As mentioned above, immersed SYZ fibers are crucial objects to complete the SYZ construction.  A local (and noncommutative) mirror construction using immersed Lagrangians was developed by the joint works \cite{CHLabc,CHLnc} of the first named and third named authors with C.-H. Cho, based on the foundational works of Floer theory by Fukaya-Oh-Ohta-Ono \cite{FOOO} and Akaho-Joyce \cite{AJ}.  In this paper we reveal the relation between the immersed two-sphere and the wall-crossing phenomenon.

The main difficulty in working with Lagrangian immersions is that constant disk bubblings occur at immersed points.  They lead to highly non-trivial obstructions in the moduli spaces.  The choices of perturbations to handle these obstructions are not explicit.  In particular, the Floer-theoretical gluing between the Lagrangian immersion and its smoothing depends on choices of the Kuranishi perturbations, and it is a very challenging task to derive the corresponding gluing formulas.

There are two main ingredients in our method to overcome this difficulty.  First, we use perfect Morse functions and pearl trajectories \cite{BC,FOOO-can,Sheridan-CY,Sc16} in the formulation of Floer theory.  A perfect Morse function provides a concrete minimal model for the Fukaya algebra.  It also makes the choices of perturbations more explicit in our situation.  Namely there is a delicate relation between the choice of a Morse function and the gluing formula.

The second ingredient is the functoriality of our gluing construction.  The gluing employs quasi-isomorphisms for objects in the Fukaya category.  Since they are functorial, it enables us to derive an explicit gluing formula via the wall-crossing between different smoothings of the immersed Lagrangians (which are the Clifford and Chekanov tori in our situation).  It gives a way to compute the open Gromov-Witten invariants of the immersed Lagrangian in the presence of highly non-trivial obstructions in the moduli.

The procedure of filling in the missing strata is the following.  First we need to construct the deformation space of the immersed two-sphere.  The key step is to achieve weakly unobstructedness, namely $m_0^b = W(b) \cdot \one_{\bL}$.  It ensures that the Floer cohomology is well-defined, and hence $(\bL,b)$ gives a well-defined object in the Fukaya category.

\begin{lem}
	Let $\bL$ be a product of an immersed two-sphere (with exactly one immersed point) and a torus which only bounds non-constant holomorphic disks with positive Maslov index.  Then $(\bL,uU+vV,\nabla)$ is weakly unobstructed, where $U,V$ are the immersed generators corresponding to the immersed point, and $\nabla$ is a flat $\C^\times$ connection on the torus component.
\end{lem}

Then we construct a quasi-isomorphism between the immersed sphere and the neighboring Chekanov and Clifford tori.  Our gluing construction produces the following gluing formula.

\begin{thmx}[Theorem~\ref{thm:q-iso}]\label{thm:gluing}
	Let $\bL_0$ be an immersed Lagrangian two-sphere, which has two degree-1 immersed generators $u,v$ coming from the immersed point.  Let $\bL_1$ be a Lagrangian torus obtained from smoothing of $\bL_0$ at the generator $u$.  Let $x$ be the holonomy variable associated to the vanishing circle of $\bL_1$, and $y$ be the holonomy variable of another circle such that they form a basis of $H_1(\bL_1,\Z)$.  
	Then there exists a perfect Morse function of the two-sphere such that 
	$(\bL_0,uU+vV) \cong (\bL_1,\nabla^{(x,y)})$  if and only if
	$$ u=y,\,\, uv = 1+x $$
	where $(u,v)\in\Lambda_0 \times \Lambda_+$,$(x,y)\in\Lambda_{\rm U}^2$.
\end{thmx}

\begin{remark}
	Fukaya has explained the application of immersed Lagrangian spheres and their smoothings to wall-crossing and mirror symmetry in his talks.  The above theorem realizes his idea by using quasi-isomorphisms between Lagrangian branes.  This method was announced and briefly explained in \cite[Section 5]{HL18}.  We provide the details and proofs in this paper.
\end{remark}

\begin{remark} 
	Very recently, Dimitroglou Rizell-Ekholm-Tonkonog \cite{DET} took the Chekanov-Eliashberg algebra approach to immersed Lagrangian surfaces.  It uses the interesting relationship between Legendrians and Lagrangians which was investigated in \cite{EkLe,EkS14,STW}. 
	
	Comparing to the Lagrangian Floer theory in \cite{FOOO} and \cite{AJ}, the Chekanov-Eliashberg approach has an additional invertible generator $t$ and the built-in relation $uv=1-t$.  It is particularly adapted to the surface case, and has the nice feature that $u,v$ can be taken to be $\C$-valued.  They applied this to the study of Lagrangian surgeries and mutations.
	
	On the other hand it is essential for us to work over the Novikov field  $\Lambda$, and homological mirror symmetry relies delicately on the associated $T$-adic topology. We first take $(u,v)$ in $(\Lambda_+ \times \Lambda_0) \cup (\Lambda_0 \times \Lambda_+)$ to ensure convergence.  After we deduce the gluing formula in Theorem \ref{thm:gluing}, it follows that $(\bL_0,uU+vV)$ for $u,v\in\C^2$ with $uv\not=1$ are well-defined objects in the subcategory generated by $\{(\bL_0,uU+vV),(\bL_1,\nabla_{(x,y)})\}$.  (See Remark \ref{rmk:t-power} for the relation of the $t$ variable in \cite{DET} and our formulation.)
\end{remark}

\subsection{Floer theoretical mirror of Grassmannians}
As an application, we construct mirror symmetry for flag manifolds.  Partial flag manifolds serve as an important class of Fano varieties and have attracted a lot of attention in the study of mirror symmetry.  In the pioneering work of Hori-Vafa \cite{HV}, Landau-Ginzburg (LG for short) mirrors for flag manifolds were proposed using \emph{T-duality} and physical derivations.  Closed string mirror symmetry for flags and their complete intersections was derived by Givental-Kim \cite{GK95}, Lian-Liu-Yau \cite{LLY1,LLY2}, Kim \cite{Kim99} and Joe-Kim \cite{JK03}.  Motivated by the Peterson variety representation of the quantum cohomology, Rietsch \cite{Rie} gave a Lie-theoretical construction of LG mirrors which can be understood as partial compactifications of those in Eguchi-Hori-Xiong and Hori-Vafa \cite{EHX97, HV}.

Nishinou-Nohara-Ueda \cite{NNU} computed the disk potential of a regular Gelfand-Cetlin fiber of partial flag manifolds including $\Gr(2,\C^n)$ which agrees with the prediction of Hori-Vafa.
However, \emph{the number of critical points of the disk potential of a regular fiber is in general smaller than the dimension of the quantum cohomology ring.}  There are not enough regular fibers to generate the Fukaya category and hence they are insufficient for the study of mirror symmetry.  

The reason is that there are Lagrangian spheres over the boundary of the base polytope, which do not intersect torus fibers over the interior and hence cannot be probed.  However Lagrangian spheres are rigid and cannot be directly used for the SYZ construction (which requires Lagrangian deformations).  In order to construct the complete SYZ mirror, one needs to perturb the Gelfand-Cetlin system to construct another Lagrangian fibration which `pushes in' the Lagrangian spheres sitting in the boundary. 

For this purpose, Nohara-Ueda \cite{NU14} constructed generalized Gelfand-Cetlin systems for $\Gr(2,\C^n)$. Interpolations between different generalized Gelfand-Cetlin systems give Lagrangian fibrations with interior discriminant loci. In \cite{NUclu}, they glued deformation spaces of toric fibers to cover a part of the Rietsch's mirror.  

We need to glue in the singular fibers in order to obtain the complete mirror.  The singular fibers over discriminant loci are products of immersed spheres with tori in this case.  They can be used in place of Lagrangian spheres in Gelfand-Cetlin systems to generate the Fukaya category.

We apply our gluing technique to $\Gr(2,\C^n)$ and obtain the Rietsch's mirror. 
The classification of Lagrangian Gelfand-Cetlin fibers in \cite{CKO, CK-mono} alludes to the locations of Lagrangians that are expected to be essential in the Fukaya category. 
 
We construct the Maurer-Cartan deformation spaces of Lagrangians in each chart, namely the monotone Chekanov torus $\scr{L}_1$, the monotone Clifford torus $\scr{L}_2$, and an immersed Lagrangian $\scr{L}_0$ which is topologically a product of two-dimensional immersed spheres with a torus.  
Then we use Floer theoretical isomorphisms to glue their deformation spaces and obtain the following.

\begin{thmx}[Theorem~\ref{theorem_RietschSYZ} and Theorem~\ref{theorem_RietschSYZgr2n}]\label{thm:Gr24}	
The mirror glued from the Maurer-Cartan deformation spaces of Lagrangians in local charts of $\Gr(2,\C^n)$ is equals to the Rietsch's LG mirror $(\check{X}, W_\textup{Rie})$ in \cite{MR}, where
$$\check{X}=\mathrm{Gr}(2, \C^n) \backslash \scr{D},$$ ${\scr{D}} := \{p_{1,2} \cdot p_{2,3} \cdots p_{n-1,n} \cdot p_{1,n} = 0\}$, $[p_{i,j}]$ are the Pl{\"u}cker coordinates of the dual Grassmannian $\mathrm{Gr}(2, \C^n)$, and
$$
W_{\textup{Rie}}\left( \left[ p_{i,j} \right] \right) := q \frac{p_{2,n}}{p_{1,2}} + 
\sum_{j=2}^{n-1} \frac{p_{j-1,j+1}}{p_{j,j+1}} + \frac{p_{1,n-1}}{p_{1,n}}.
$$

Moreover, the monotone immersed Lagrangians used in the mirror construction are non-displaceable, as they support critical points of the potential.
\end{thmx}


\begin{remark} \label{rmk:Rie}
	Our result gives an enumerative meaning of the Rietsch's LG mirror.  Namely, the open Gromov-Witten invariants of our immersed Lagrangians can be extracted combinatorially from the Rietsch superpotential. For instance, each of terms in \eqref{eqn:Wforuvzw24} corresponds to a count of Maslov-2 holomorphic disk whose boundary-classes can be read off from $z_0$ and $w_0$, and the disk should have corners precisely in accordance with powers of $u,v$ appearing in the term. See Corollary~\ref{thm:W^L} for more general cases.
	
	It is informative to compare with the works \cite{CLLT,CCLT13}, which extract open Gromov-Witten invariants from the LG mirrors of toric Calabi-Yau and semi-Fano manifolds.  Even though the mirror map is trivial in this situation, there are still non-trivial coordinate changes for the open parameters, and there are infinitely many non-zero invariants.  
	
	Technically the new ingredient in this case is \emph{the obstruction coming from constant polygons}.  Our method provides an effective way to compute these invariants.\end{remark} 

To illustrate, let us consider the simplest non-trivial example $\mathrm{Gr}(2,\C^4)$.  Only one immersed Lagrangian is used in the construction.  The disc potential can be identified with $W_{\mathrm{Rie}}$ as follows.  

\begin{example} \label{ex-intro}
	For $\mathrm{Gr}(2,\C^4)$, 
	the disk potential of the monotone immersed Lagrangian $\scr{L}_0 \cong \mathcal{S}_2\times T^2$, where $\mathcal{S}_2$ denotes the immersed sphere with exactly one-self nodal point, equals to
	\begin{equation}\label{eqn:Wforuvzw24}
	W_{\scr{L}_0} (u, v, z_0, w_0) = T \cdot \left(\frac{v}{(uv-1)z_0}  + u + {v w_0}+ \frac{uz_0}{w_0} 
	\right)= T \cdot \left(-\frac{v}{z_0} \sum_{i=0}^\infty (uv)^i  + u + {v w_0} + \frac{uz_0}{w_0} \right)
	\end{equation}
	where $T$ is the Novikov parameter, $u,v$ are the immersed variables and $z_0, w_0$ are the holonomy variables.  The terms $(uv)^i$ are contributed from constant polygons with corners being the immersed sectors $U$ and $V$.
	
	$$W_{\scr{L}_0}=W_{\textup{Rie}}$$
	by the coordinate change
	
	$$T\cdot u = \frac{p_{1,3}}{p_{2,3}},\,\, T^{-1} \cdot v =  \frac{p_{2,4}}{p_{1,4}},\,\, T^2 \cdot w_0 =  \frac{p_{1,4}}{p_{3,4}},\,\, T^2 \cdot z_0 =  \frac{p_{2,3}}{p_{3,4}},\,\, q = T^4.$$
	We found a combinatorial formula which is explained in Section \ref{sec:WRie=WL}.
	(Recall the Plucker relation 
	$p_{1,2} p_{3,4} - p_{2,4} p_{1,3} + p_{1,4} p_{2,3}=0$ and $p_{i,j}$ are homogeneous coordinates.)
	
	In particular, the Lagrangian $\scr{L}_0$ is non-displaceable. 
\end{example}

We use the gluing formula to deduce the disk potential for the immersed Lagrangians in $ \Gr(2,\C^n)$.  This method works in general for immersed spheres in other manifolds.

\begin{remark}
	$\Gr(2,\C^n)$ can be understood as a smoothing of a toric variety which has conifold singularities.  Smoothing of a local conifold singularity can be understood via Minkowski decomposition \cite{Altmann}.  See the works of Gross \cite{Gross-eg} and the last author \cite{L14} for Lagrangian fibrations and wall-crossing in a local smoothing.
\end{remark}

The mirror construction for flag varieties beyond type-A has not been well-understood.
They still admit Gelfand-Cetlin systems serving as Lagrangian fibrations.  However we again encounter the same problem that Lagrangian spheres sit over  the boundary of the base polytope.  We need to construct new Lagrangian fibrations so that the possibly non-displaceable spheres over the boundary are pushed into the interior of the base.  
In this paper we carry this out for the type-B flag manifold $\mathrm{OG}(1,\C^5)$.

\begin{figure}[h]
	\begin{center}
		\includegraphics[scale=0.4]{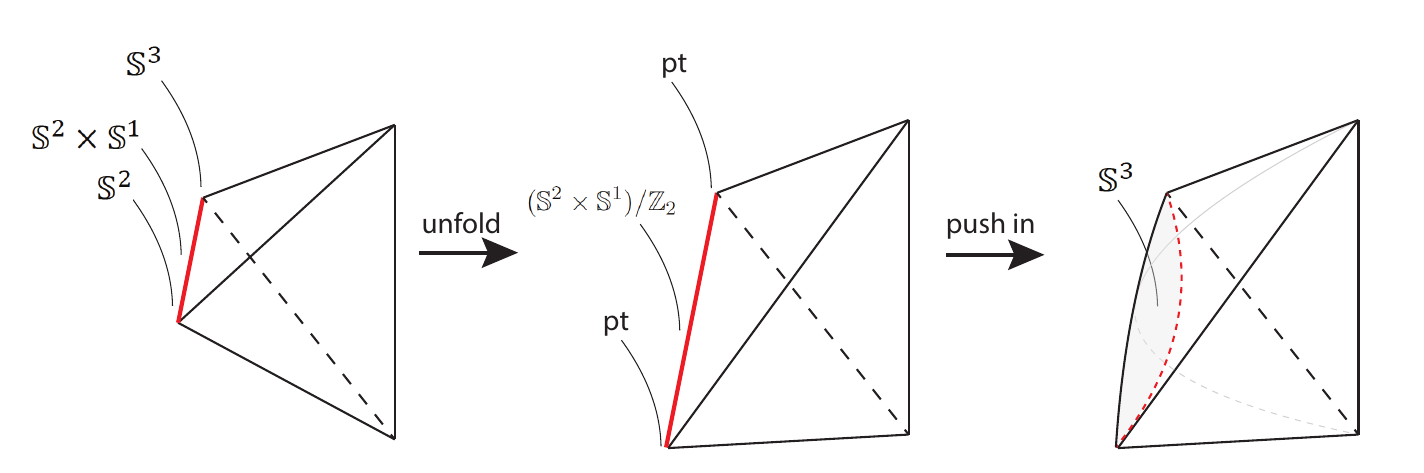}
		\caption{The Gelfand-Cetlin polytope of $\OG(1,\C^5)$ (left) and the base of the Lagrangian fibration that we use (right).}
		\label{Fig_OG1-5-unfold}
	\end{center}
\end{figure}

The Gelfand-Cetlin system of $\OG(1,\C^5)$ studied in \cite{NNU12} has Lagrangian spheres $\bS^3$ contained in the boundary of the base polytope (see the left of Figure \ref{Fig_OG1-5-unfold}).  The corresponding cone in the fan picture is generated by $(1,1,0),(-1,1,0),(0,0,1)$ which has determinant 2.  
In contrast to the case of $\Gr(2,\C^4)$, The triangle $\mathrm{Conv}\{(1,1,0),(-1,1,0),(0,0,1)\}$ does not have a Minkowski decomposition to describe the smoothing (of the corresponding toric variety). To remedy, we turn to a different local model to push in the Lagrangian spheres.

As in Theorem \ref{thm:Gr24}, we derive the LG mirror of an immersed Lagrangian which agrees with the result of Pech-Rietsch-Williams \cite{PRW}. 

\begin{thmx}[Theorem \ref{theorem_PRWrecover}]
The disk potential of an immersed Lagrangian $\scr{L}_0$ in $\mathrm{OG}(1,\C^5)$ is equals to 
$$
W_{\scr{L}_0} (u,v, z_0) = T\cdot \left(v + vz_0 + \frac{u^2}{z_0 (uv-1)}\right).
$$
In particular, the Lagrangian $\scr{L}_0$ is non-displaceable. 

Furthermore, the glued mirror equals to the Rietsch's LG mirror $(\check{X}, W_q)$ of $\OG(1,\C^5)$ in \cite{PRW}  where
$$
\check{X}= \CP^{3} \backslash \scr{D},$$ 
${\scr{D}} := \{p_{0} \cdot (p_1p_2 - p_0p_3) \cdot p_{3} = 0\}$, $[p_{i}]$ are the homogeneous coordinates of $\CP^3$,
$$
W_q = \frac{p_1}{p_0} + \frac{p_2^2}{p_1p_2 - p_0 p_3} + q \frac{p_1}{p_{3}}\colon \check{X} \to \C.
$$
\end{thmx}

As an application of our mirror construction, we prove homological mirror symmetry for $\OG(1,\C^5)$ with the help of the generation result in \cite{Sheridan-Fano}. We can also prove homological mirror symmetry for $\Gr(2,\C^4)$, yet  assuming generation of Lagrangians (which should be able to deduce from the work of Evans-Lekili \cite{EL19} or the announced work of Abouzaid-Fukaya-Oh-Ohta-Ono). As we will see in Section \ref{sec:HMS}, the mirror LG models for both cases are semi-simple, and hence it is enough to show that the associated functor gives isomorphisms of endomorphisms spaces for finitely many mirror pairs of objects.  (See \cite{CHLtoric} for the toric case.)

\begin{thmx}\label{prop:HMSpart}
Let $W$ be the Rietsch mirror (over $\Lambda$) of $X=\OG(1,\C^5)$ or $\Gr(2,\C^4)$. 
Then for each critical value $\lambda$, $D^b \mathcal{A}_\lambda \simeq D^b \MF (W-\lambda)$ where $\mathcal{A}_\lambda$ is the sub-Fukaya category of $X$ generated by the objects with the same potential value $\lambda$. Furthermore, we have $\oplus_\lambda D^b \mathcal{A}_\lambda = D^b \Fuk (X)$ for $X=\OG(1,\C^5)$, and the homological mirror symmetry for $\OG(1,\C^5)$ follows.\end{thmx}

\subsection*{Notations}

We need to work over the Novikov ring.  In general the mirrors are rigid analytic space over the Novikov ring (for instance see \cite{FOOO-T} in the toric case).  We use the following notations.
\begin{eqnarray*}
	\Lambda_+ &=& \left\{\sum_{i=0}^\infty a_i T^{A_i} \mid A_i > 0 \textrm{ increase to } +\infty, a_i \in \C \right\},\\
	\Lambda_0 &=&\left\{\sum_{i=0}^\infty a_i T^{A_i} \mid A_i \geq 0 \textrm{ increase to } +\infty, a_i \in \C \right\},\\
	\Lambda &=&\left\{\sum_{i=0}^\infty a_i T^{A_i} \mid A_i \textrm{ increase to } +\infty, a_i \in \C \right\},\\
	\Lambda_{\rm U} &=& \, \C^\times \oplus \Lambda_+.
\end{eqnarray*}
They are equipped with the valuation function 
$$\mathrm{val}:\sum_{i=0}^\infty a_i T^{A_i} \mapsto A_0 \textrm{ and } \mathrm{val}(0) = +\infty.$$
$\Lambda_{\rm U}$ is the valuation-zero subset of $\Lambda$ which forms a multiplicative group.


\subsection*{Acknowledgment}
The authors express our deep gratitude to Cheol-Hyun Cho and Dongwook Jwa.
The third named author is very grateful to Naichung Conan Leung for telling him the idea of pushing in the Lagrangian spheres at corners of Gelfand-Cetlin system in 2010.
We thank to Dmitry Tonkonog for introducing his work with T. Ekholm and G. Dimitroglou Rizell on immersed Lagrangian surfaces using Chekanov-Eliashberg algebra in May 2018.  Their approach can directly work over $\C$ and has very close relations with Legendrians. 
The second named author thanks Yuichi Nohara for explaining his works with Kazushi Ueda on wall-crossings and cluster transformations on Grassmannians of 2-planes, and Yunhyung Cho and Yong-Geun Oh for drawing his attention to Gelfand-Cetlin systems.
The first named author would like to thank Bong Lian and Shing-Tung Yau for constant support and
encouragement. His work is substantially supported by Simons Collaboration Grant on Homological Mirror Symmetry. The work of the third named author is supported by Simons Collaboration Grant for Mathematicians.

\section{Mirror functor and gluing construction}

We quickly review the basics of Lagrangian Floer theory in this section, and recall the construction of mirrors via Maurer-Cartan deformations of a Lagrangian. 

Consider a Lagrangian $\mathbb{L}$ in a symplectic manifold $X$ (not necessarily compact). A degree odd element $b \in CF(\mathbb{L},\mathbb{L})$ is called \emph{Maurer-Cartan} if it satisfies
\begin{equation}\label{eqn:defmc}
 m_0^b = m_0(1) + m_1 (b) +  m_2 (b,b) + \cdots = 0.
\end{equation}
Here, $b$ is usually assumed to have coefficient in $\Lambda_+$ to guarantee the convergence of \eqref{eqn:defmc}. In actual application later, we will sometimes weaken this condition by carefully studying  convergence of \eqref{eqn:defmc}. $b$ is called \emph{weak Maurer-Cartan} if the left hand side of \eqref{eqn:defmc} is a constant multiple $W(b)$ of the unit class in $CF(\mathbb{L},\mathbb{L})$.

For given a (weak) Maurer-Cartan element $b$, one can deform $\mathbb{L}$ as an object of $\Fuk(X)$ by taking boundary deformation of the $A_\infty$-structure on $\mathbb{L}$ in the sense of Fukaya-Oh-Ohta-Ono \cite{FOOO}. We will denote the resulting object by $(\mathbb{L},b)$. For instance, the Floer differential on $CF((\mathbb{L},b),(\mathbb{L},b))$ is given by
$$ m_1^b (x) = \sum m_k (b, \cdots, b, x, b, \cdots,b) $$
Higher operations as well as $A_\infty$-operations involving not only $(\mathbb{L},b)$ but also other Lagrangians are deformed in a similar way (by inserting arbitrary many $b$'s in all possible positions).
For more details, see for e.g., \cite{FOOO}. Notice that $m_0^b = 0$ for a Maurer-Cartan element $b$.  In such a case $(\mathbb{L},b)$ is unobstructed, whose Floer cohomology is well-defined. Throughout the article, $b$ will be either a linear combination of immersed generators (for an immersed Lagrangian) or a flat connections on a Lagrangian torus.

The collection of solutions $b$ to \eqref{eqn:defmc} can be thought of as a (unobstructed) deformation space of the Lagrangian $\mathbb{L}$ which will give a local chart of the mirror we are going to construct. We denote this space by $\mathcal{MC} (\mathbb{L})$. (To be more precise, one should take the quotient by gauge equivalences, but we will use canonical models later so that the gauge equivalence is trivial.) We will write $\mathcal{MC}_{weak} (\mathbb{L})$ for the space of weak Maurer-Cartan elements of $\mathbb{L}$. In this case, we additionally have a function $W$ on $\mathcal{MC}_{weak} (\mathbb{L})$ called the Lagrangian Floer disk potential of $\mathbb{L}$, and $(\mathcal{MC}_{weak} (\mathbb{L}),W)$ will serve as a local chart of the mirror LG model. We have the following local mirror functor by the construction in \cite{CHLnc}.

\begin{proposition}[\cite{CHLabc, CHLnc}]\label{prop:localfunc}
There exists a natural $A_\infty$-functor 
$$\mathcal{F}^{\mathbb{L}} \colon \Fuk (X) \to  D^b \,{\rm Coh} (\mathcal{MC} (\mathbb{L})),$$ 
or 
$$\mathcal{F}^{\mathbb{L}} \colon \Fuk (X) \to \MF (\mathcal{MC}_{weak} (\mathbb{L}),W)$$ in the weak Maurer-Cartan case, which sends $L \in \Fuk (X)$ to $\left( CF((\mathbb{L},b),L),m_1^b \right)$. When $b$ varies, it gives a bundle on the Maurer-Cartan space with $\left( m_1^b \right)^2 = 0$ (i.e. a cochain complex of bundles), or a bundle on the weak Maurer-Cartan space with $\left( m_1^b \right)^2 = W(b) \cdot id$ (that is, a matrix factorization of $W$). 
\end{proposition}

Higher components of $\mathcal{F}^{\mathbb{L}}$ are defined making use of deformed higher $m_k$'s analogously to $m_1^b$. See \cite{CHLabc} and \cite{CHLtoric} for more details. 

\subsection{Gluing mirror charts}\label{subsec:gluecharts}
We next explain how to glue these local deformation spaces from various Lagrangians. Consider two Lagrangians $\mathbb{L}_1$ and $\mathbb{L}_2$ in $X$, which intersect cleanly with each other. We will glue two deformation spaces $\mathcal{MC} (\mathbb{L}_1)$ and $\mathcal{MC}(\mathbb{L}_2)$ (or $\mathcal{MC}_{weak} (\mathbb{L}_1)$ and $\mathcal{MC}_{weak} (\mathbb{L}_2)$ depending on the situation) by making a choice of \emph{quasi-isomorphism} between two objects $(\mathbb{L}_1,b_1)$ and $(\mathbb{L}_2,b_2)$ in the Fukaya category of $X$.

\begin{definition}\label{def:qisomLag}
Two objects $L$ and $L'$ in $\Fuk (X)$ are said to be \emph{quasi-isomorphic} if there exist $\alpha \in CF^0 (L,L')$ and $\beta \in CF^0 (L',L)$ such that 
$$m_1 ( \alpha) = m_1 (\beta) = 0 \quad m_2 (\alpha,\beta) = \one_{L} + m_1 (\gamma) \quad m_2(\beta,\alpha) =\one_{L'} + m_1 (\delta)  $$
for some $\gamma \in CF (L,L)$ and $\delta \in CF(L',L')$.
Such $\alpha$ and $\beta$ are called quasi-isomorphisms in this case. They are called strict isomorphisms if $m_2$ between them are strictly the units.
\end{definition}

\begin{remark}
One can easily check that quasi-isomorphisms give rise to an equivalence relation in the set of objects in $\Fuk(X)$. Also, note that quasi-isomorphisms between graded Lagrangians are automatically strict due to degree reason.
\end{remark}

Now suppose that there exists a function 
$$ f \colon U \subset \mathcal{MC} (\mathbb{L}_1) \to \mathcal{MC} (\mathbb{L}_2)$$
together with a fixed quasi-isomorphism $\alpha : (\mathbb{L}_1,b_1) \stackrel{\simeq}{\to} (\mathbb{L}_2, b_2=f(b_1))$ for $b_1 \in U$ which does not depend on $b_1$. In this case, we glue two spaces by taking
\begin{equation}\label{eqn:gluedmc}
\mathcal{MC} (\mathbb{L}_1) \cup_f \mathcal{MC} (\mathbb{L}_2) = \mathcal{MC} (\mathbb{L}_1) \cup \mathcal{MC} (\mathbb{L}_2) / b_1 \sim f(b_1).
\end{equation}
Namely, we identify $U \subset \mathcal{MC} (\mathbb{L}_1)$ and $f(U) \subset \mathcal{MC}(\mathbb{L}_2)$. More detailed explanations on the topology of the resulting space will be given in later sections in which we actually apply the construction to our main examples.

The gluing of weak Maurer-Cartan spaces can be performed in a similar way. Let $W_1$ and $W_2$ be potentials for $\mathbb{L}_1$ and $\mathbb{L}_2$. Then the existence of an isomorphism implies $W_1|_U =f^\ast W_2|_{f(U)}$ since
\begin{equation}\label{equ_m0m0}
0 = m_1^2 (\alpha) = \left( W_1 (b_1) - W_2 (f(b_1)) \right) \cdot \alpha
\end{equation}
by $A_\infty$-relations.

\begin{proposition}[\cite{CHLgl}]
There exists a natural $A_\infty$-functor
$$\mathcal{F} \colon \Fuk (X) \to  \MF (\mathcal{MC}_{weak} (\BL_1) , W_1) \times^{h} \MF (\mathcal{MC}_{weak} (\BL_2) , W_2)$$
where the right hand side is a homotopy fiber product of two dg-categories over their common intersection $\MF (\mathcal{MC}_{weak} (\BL_1) \cap \mathcal{MC}_{weak} (\BL_2), W_1| = W_2|)$. (Similar statement holds for $D^b \, {\rm Coh}(\mathcal{MC} (\BL_i))$'s in the case of (strict) Maurer-Cartan deformation.) 
\end{proposition}

We do not give a precise definition of the homotopy fiber product above as it will not be used in the paper. Roughly speaking, it consists of local objects together with their gluing data which are isomorphisms up to homotopy. (See for e.g., \cite{BBB13} for details about the homotopy fiber product of dg-categories.)
On the object level, the functor $\mathcal{F}$ sends an object $L \in \Fuk (X)$ to a tuple
$ (\mathcal{F}^{\mathbb{L}_1} (L), \mathcal{F}^{\mathbb{L}_1} (L) \stackrel{m_2 (\alpha,\cdot)}{\longrightarrow} \mathcal{F}^{\mathbb{L}_2} (L) ,\mathcal{F}^{\mathbb{L}_2} (L) )$.

In what follows, the local functor given in Proposition \ref{prop:localfunc} will be enough to derive   a homological mirror symmetry statement since the mirror potentials in our case are Morse and the matrix factorization categories localizes at the critical points each of which is contained in a single chart. We will revisit this point in Section \ref{sec:HMS}.

\subsection{Review of immersed Lagrangian Floer theory}

We briefly review Floer theory for immersed Lagrangians. 
Let $\iota \colon L \to M$ be a compact immersed Lagrangian in a symplectic manifold $(M, \omega)$. We first assume that $\iota^{-1}(p)$ is at most two points for each $p \in \iota(L)$, and two branches of $\iota (L)$ intersects at $p$ transversally.

Consider the fiber product
$$
R:=L \times_\iota L = \{(p,q) \in L \times L ~|~ \iota(p) = \iota(q) \},
$$
which consists of several connected components. We label components of $R$ as follows. Obviously $R$ contains the diagonal $\{(p,p) \in L \times L\}$ as one of connected components, which will be denoted by $R_0$. The other components are isolated points of the form 
$$
(p_-, p_+) \in L \times L 
$$ 
with $\iota(p_-) = \iota(p_+) = p$ and $p_- \neq p_+$. 
For later use, define an involution $\sigma$ by  
\begin{equation}\label{eqn:sigmaR}
\sigma \colon R \to R, \quad
\begin{cases}
id \,\,\, \mbox{on} \,\, R_0, \\
(p_-, p_+) \leftrightarrow (p_+, p_-).
\end{cases}
\end{equation}
by fixing the whole $R_0$ and swapping $(p_-, p_+)$ to $(p_+, p_-)$. 
The Floer complex of $L$ is generated by $C^\ast (L)$ (or $H^\ast (L)$) and the elements of $R \setminus R_0$. The latter will be called the \emph{immersed generators} of $CF(L,L)$. Hence each self-intersection point of $\iota: L \to M$ gives to two generators of $CF(L,L)$, which intuitively will describe corners of a holomorphic polygon jumping between two branches meeting at the corresponding self-intersection.

Fix $k \in \Z_{\geq 0}$ and $\beta \in H_2(X,  \iota(L))$. Let $J$ be an $\omega$-compatible almost complex structure. 
$A_\infty$-structure on $CF(L,L)$ is defined by counting stable $J$-holomorphic maps $\varphi \colon (\Sigma, \pa \Sigma) \to (M, \iota(L))$
from a genus $0$ bordered Riemann surface $\Sigma$ together with mutually distinct marked points $z_0, \cdots, z_k$.

As mentioned, a marked point for holomorphic maps contributing to $CF(L,L)$ may map to a corner at some self-intersection point of $\iota : L \to M$. 
In order to record these corners, we additionally attach a map $\alpha$ to each holomorphic disks for $CF(L,L)$
$$
\alpha \colon \{0, \cdots, k\} \to R.
$$
such that the following properties are satisfied$\colon$
\begin{enumerate}
\item The marked points respect the counter-clockwise orientation,
\item None of marked points are nodes, 
\item The map $\alpha$ determines which immersed sectors marked points pass through. Namely, for all $i$,  
$$
\begin{cases}
\displaystyle (\varphi(z_i), \varphi(z_i)) \in R &\quad \mbox{if $\alpha(i) \in R_0$}, \\
\displaystyle \left( \lim_{z \circlearrowleft z_i} (\iota^{-1} \circ \varphi) (z),  \lim_{z \circlearrowright z_i} (\iota^{-1} \circ \varphi) (z) \right) = \alpha(i) \in R  &\quad \mbox{if $\alpha(i) \notin R_0$}.
\end{cases}
$$
\end{enumerate}
We denote the moduli spaces of such maps with additional data $\alpha$ by $\mcal{M}_{k+1}(\alpha, \beta, J)$.

We finally define the evaluation maps
$$
\begin{cases}
\textup{ev}_i \colon {\mcal{M}}_{k+1} (\alpha, \beta, J) \to L \coprod R  \\
\textup{ev}_0 \colon {\mcal{M}}_{k+1} (\alpha, \beta, J) \to L \coprod R 
\end{cases}
$$
by
$$
\textup{ev}_i \left([\Sigma, \overrightarrow{z}, u, \ell, \overline{u}] \right) :=
\begin{cases}
\overline{u}(\zeta_i) \in L, &\mbox{$\alpha(i) \in R_0$},\\
\alpha(i) \in R, &\mbox{$\alpha(i) \notin R_0$}
\end{cases}
$$
for $i = 1, \cdots, k$ and
$$
\textup{ev}_0 \left([\Sigma, \overrightarrow{z}, u, \ell, \overline{u}] \right) :=
\begin{cases}
\overline{u}(z_0) \in L, &\mbox{$\alpha(0) \in R_0$},\\
\sigma \circ \alpha(0) \in R, &\mbox{$\alpha(0) \notin R_0$}.
\end{cases}
$$
where $\bar{u}$ is the lifting of $\varphi|_{\partial D^2}$ to $R$. Having this, $A_\infty$-operation on $CF(L,L)$ can be defined in a usual way. One can generalize the definition into the case where $L$ has a clean self-intersection by counting similar types of holomorphic disks attached with Morse flow lines in self-intersection loci. In this case, critical points of the Morse functions on the self-intersection loci serve as immersed generators of the Floer complex.

We will use linear combinations of degree one immersed generators as our Maurer-Cartan deformation parameter $b$. Here, the degree of an immersed generator is determined in a usual way by two local Lagrangians branches intersecting at the corresponding self-intersection. The Maurer-Cartan equation for such $b$ is contributed by holomorphic polygons with degree one corners at self-intersections of $\iota(L)$.
 See \cite{AJ} and \cite{PWim} for more details.

\subsection{Gauge hypertori}\label{subsec:fdconn}
When a Lagrangian is equipped with a $\Lambda_{\rm U}$-line bundle together with a flat connection $\nabla$, one can deform $A_\infty$-operations on $CF(L,L)$ by measuring the parallel transport with respect to $\nabla$ along a boundary of a holomorphic disk. More precisely, a summand  $m_{k,\beta}$ of $m_k$ contributed by disks of class $\beta$ is multiplied by $hol_\nabla (\beta)$, where we regard the flat connection $\nabla$ as an element $hol_\nabla \in \hom (\pi_1 (L), \Lambda_{\rm U})$. For $A_\infty$-operations among several different Lagrangians, we similarly measure the parallel transport along each boundary component of a disk with respect to the connection on the Lagrangian over which this boundary component lies. Such a deformation gives rise to $(\Lambda_{\rm U})^{h_1}$ as a part of the Maurer-Cartan deformation space $\mathcal{MC} (L)$ (or $\mathcal{MC}_{weak} (L)$) when, for e.g., $L$ itself is already (weakly) unobstructed without further deformation. Here, $h_1$ is the first Betti number of $L$.

Throughout the paper, we will choose a particular type of a connection to represent the equivalent class of $\Lambda_{\rm U}$ flat line bundles, for which nontrivial effect of parallel transport is concentrated near chosen cycles in $H^{n-1} (L)$ (here $\dim_\R L = n$). We give more detailed explanation when $L$ is a torus, as is the case in most of our application below. For $L \cong \R^n / \Z^n$, we fix (oriented) codimension 1 tori $H_i = \epsilon_i + \R \langle e_i \rangle$ for $\epsilon_i \in \R/\Z$, which are called the {\em gauge hypertori}. Then the parallel transport over a path $\gamma$ is given by multiplying $z_i^{\pm}$ whenever $\gamma$ runs across $H_i$. Here, the sign in the exponent is determined by the parity of the intersection $\gamma \cap H_i$, and $z_i$ is the holonomy of $\nabla$ along the 1 cycle $PD [H_i]$. See \cite{CHLtoric} for more details.

\section{Completing SYZ mirror by immersed Lagrangians}\label{sec_wallcrossingimmersedLag}
In this section, we study a local model $\C^2 \backslash \{ab=\varepsilon\}$ which will be crucially used in the mirror construction of Grassmannians. The mirror construction for $\C^2 \backslash \{ab=\varepsilon\}$ via wall-crossing has been studied by Auroux \cite{auroux07} and many others. In \cite[Section 5]{HL18}, the first and the third named author briefly described the (partial) compactification of the mirror by gluing  the deformation space of an immersed sphere. In this section, we provide the detailed statements and their proofs.

We begin with the conic fibration $\Pi \colon \C^2 \to \C $ defined by $(a,b) \mapsto ab$. The fibration is equipped with the Hamiltonian $\mathbb{S}^1$-action on $\C^2$ given by
\begin{equation}\label{equ_s1actiononconicfib}
\mathbb{S}^1 \times \C^2 \to \C^2, \quad (\theta, (a,b)) \mapsto (e^{-i \theta} a, e^{i \theta} b),
\end{equation}
whose moment map for this $\mathbb{S}^1$-action is $\mu_\theta (a,b) = |b|^2 - |a|^2$. 
Fix a positive number $\varepsilon$ and a simply connected curve $\gamma$ in the base $\C \backslash \{ \varepsilon \}$. 
The union of $\mathbb{S}^1$-orbits satisfying $|a| = |b|$ in the conic fibers over $\gamma$ form a Lagrangian torus, denoted by
\begin{equation}\label{equ_Lagrangiantorusovergammai}
\mathbb{T}_{\gamma} := \{(a,b) \in \Pi^{-1}(\gamma) ~:~ |a| = |b| \}.
\end{equation}

The most standard way to get the mirror is to take concentric circles around $\varepsilon$ which gives the SYZ fibration (with one singular torus fiber). In other words, those Lagrangians form a special Lagrangian torus fibration. 

\begin{proposition}[Proposition 5.2 in \cite{auroux07}]
	Let $\gamma(r)$ be the circle around $\varepsilon$ with radus $r$. 
	The Lagrangian torus $\mathbb{T}_{\gamma(r)}$ is special Lagrangian with respect to the holomorphic volume form $\Omega = {da \wedge db}/{(ab-\varepsilon)}$. Hence, $\mathbb{T}_\gamma$ in~\eqref{equ_Lagrangiantorusovergammai} is a graded Lagrangian.
\end{proposition}

\begin{figure}[h]
	\begin{center}
		\includegraphics[scale=0.37]{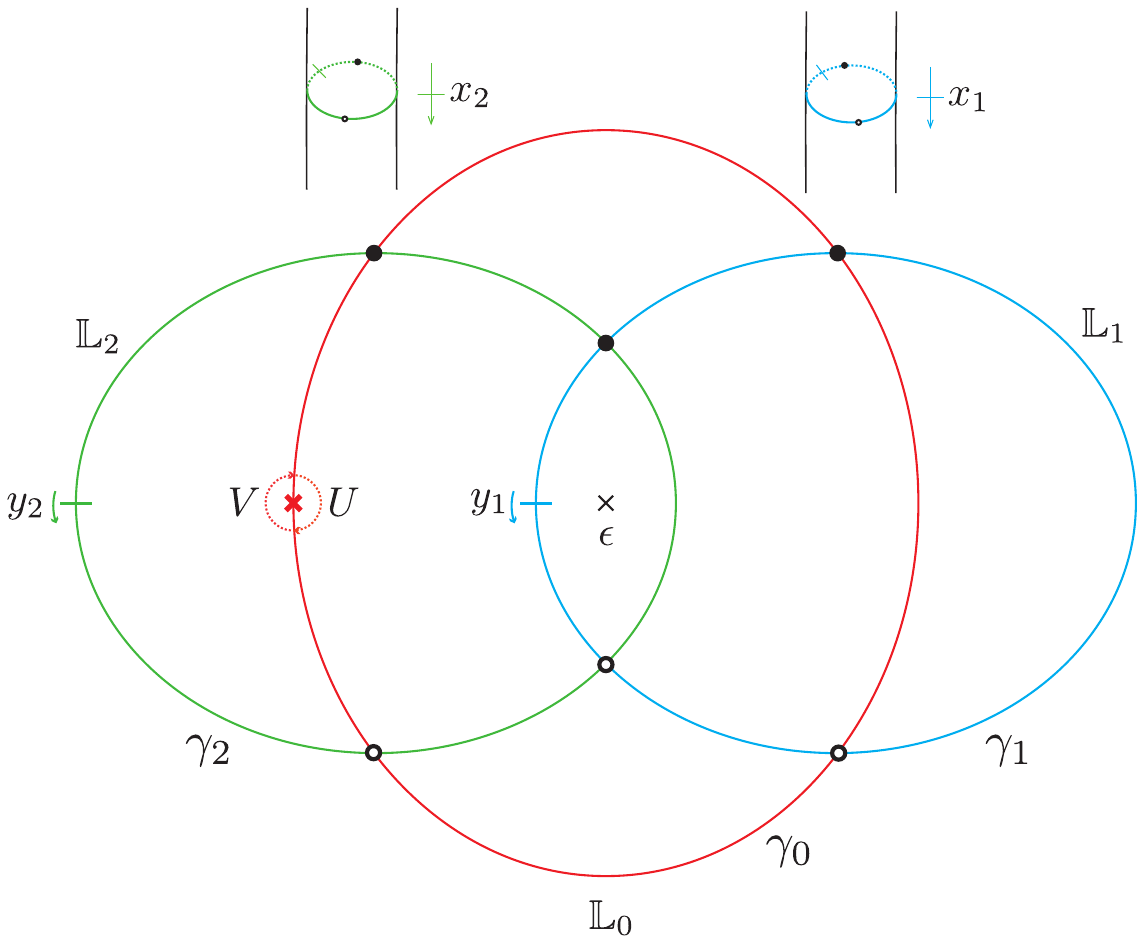}
		\caption{Simple closed curves that are base paths for the Lagrangians}
		\label{Fig_base}
	\end{center}
\end{figure}

However, in our mirror construction via Maurer-Cartan deformation spaces and their gluing, we need isomorphisms among Lagrangians (see Section~\ref{subsec:gluecharts}).  In particular we need the chosen Lagrangians to intersect.  Torus fibers do not intersect with each other, and so we need to make an alternative choice of of the following three simple closed curves$\colon$
\begin{enumerate}
\item $\gamma_{0}$ contains $\varepsilon$ and passes through $0$,
\item $\gamma_1$ contains $\varepsilon$ in its interior and $0$ in its exterior,
\item $\gamma_2$ contains $0$ and $\varepsilon$ in its interior,
\end{enumerate}
as in Figure~\ref{Fig_base}. We additionally assume that the areas bounded by $\gamma_i$'s are all the same. $\mathbb{T}_{\gamma_1},$ and $\mathbb{T}_{\gamma_2}$ will play the roles similar to those of torus fibers below/after the wall.  $\mathbb{T}_{\gamma_0}$ can be taken to be Hamiltonian isotopic to the singular torus fiber in the SYZ fibration explained above. From now on, we will write $\BL_i$ for $\mathbb{T}_{\gamma_i}$ for $i = 0, 1, 2$ for notational simplicity.

Any pair of the Lagrangians $\BL_0, \BL_1,$ and $\BL_2$ intersect cleanly with each other along disjoint circles. Later, we will use Morse model for Floer complexes among these objects, and for that we choose a generic perfect Morse function on each circle in the intersection loci. Each of them gives rise to two generators in the corresponding Floer complex with degree difference 1, which are simply its minimum and maximum. Alternatively,
one can perturb these circles by Morse functions so that they intersect transversely.

Equip the Lagrangian tori $\BL_i$ for $i = 1, 2$ with flat $\Lambda_{\rm U}$-line bundles and fix the gauge by choosing hypertori in $\BL_i$ as in Section~\ref{subsec:fdconn}. 
As depicted in Figure~\ref{Fig_base}, the parallel transport is then given by multiplying $y_i$ when a path goes across the corresponding hypertorus in $\BL_i$, which is simply the circle component of $\BL_i$ in the conic fiber over $y_i$. For the holonomy $x_i$ in the complementary direction in $\BL_i$ are roughly given in Figure \ref{Fig_base}. Since the choice of $x_i$-hypertorus is a bit delicate, we will explain this in a separate section, see Section~\ref{subsec:supple}. 
In fact, the identification of charts associated to Lagrangians $\BL_i$'s depends on this choice. (See the proof of Theorem \ref{thm:q-iso} and the discussion below the proof.)

\subsection{Unobstructedness}\label{subsec_Unobstructedness}
We will use flat connections for $\BL_1$ and $\BL_2$ and linear combinations of immersed generators $U$ and $V$ for $\BL_0$ to deform (Floer theory of) these Lagrangians. By taking suitable gradings on $\BL_0$, both $U$ and $V$ can be made of degree 1. In \cite{Alsim}, the same generators for the immersed sphere appear to have different degrees from ours, essentially due to our choice of the logarithmic (holomorphic) volume form. The $\log$-factor in the volume form (along the base direction of the conic fibration) allows us to keep the phase unchanged when we travel along a loop which starts from one of two branches (of $\BL_0$ at the self-intersection), hits the other, and comes back to the original branch.

We first prove that any such deformations satisfy the Maurer-Cartan equation.
We begin by the unobstructedness of smooth tori $\BL_i$ for $i = 1, 2$. 

\begin{lemma} \label{lem:Li_unobs}
	The pair $(\bL_i, b_{\bL_i} := \nabla^{x_i,y_i})$ is unobstructed for any $i = 1, 2$. 
\end{lemma}
\begin{proof}
	Since the torus $\bL_i$ does not bound any non-constant holomorphic disks, $m_0^{b_{\bL_i}}=0$.
\end{proof}

The immersed sphere $\BL_0$ is equipped with a perfect Morse function on $\mathbb{S}^2$ with two critical points (which are away from the immersed points). The unobstructedness of $\BL_0$ is irrelevant to a choice of Morse functions as we will see below. However, we will choose somewhat specific Morse function in Section~\ref{subsec:glueform1} mainly to obtain a simple wall-crossing formula. (See the proof of Theorem \ref{thm:q-iso}.)

The Floer cochains are spanned by the two critical points $\alpha^{\mathbb{L}_0}_0$ and $\beta_2^{\mathbb{L}_0}$ of degree $0$ and degree $2$, respectively, together with two degree $1$ immersed generators $U$ and $V$.  The generators of degree $0$ and degree $2$ are referred to as the unit class and the point class of the 2-sphere by obvious analogy.
We consider the deformation by $b_{\BL_0} = uU + vV$  of $\BL_0$, where $u,v \in \Lambda_+$.

\begin{lemma} \label{lem:L0_unobs}
	The pair $(\BL_0,b_{\BL_0} := uU + vV)$ is unobstructed.
\end{lemma}

	Since there does \emph{not} exist any non-constant holomorphic disks bounded by $\mathbb{L}_0$, the only constant stable maps are able to contribute the obstruction $m_0^b$. Namely, $m_0^b = m_{0, \beta_0}^b$ where $\beta_0$ is the trivial class in $\pi_2(\C^2 \backslash \{ ab = \varepsilon \}, \mathbb{L}_0)$, represented by a constant disk. 
Thus it should be of degree two, and hence proportional to the point class $\beta_2^{\mathbb{L}_0}$. Moreover, because of branch shifting, it is merely contributed by constant polygons supported at the self intersection point of $\BL_0$ whose corners are $(U,V,\ldots,U,V)$ or $(V,U,\ldots,V,U)$. Therefore, to verify Lemma~\ref{lem:L0_unobs}, it suffices to show that the contributions from the constant maps with corners $(U,V,\ldots,U,V)$ and with $(V,U,\ldots,V,U)$ cancel out. 
	  
	  One easy way to see such a cancellation is to consider the anti-symplectic involution $\tau \colon \C^2 \to \C$ given by $(a,b) \mapsto (\overline{b}, \overline{a})$ .
The map $\tau$ swaps the immersed sectors, i.e., $\tau$ induces the map $\tau_*$ on $R$, which coincides with the involution $\sigma$ defined in \eqref{eqn:sigmaR}. Note that $\tau$ is anti-holomorphic with respect to the standard complex structure on $\C^2$. The involution induces
the map 
\begin{equation}\label{equ_taustart}
\tau_* \colon \mcal{M}_{k+1} (\alpha, \beta_0) \to \mcal{M}^\textup{clock}_{k+1} ({\alpha}, \beta_0)
\,\, \mbox{by\, }
\begin{cases}
\tau_* (\varphi) (z) :=  (\tau \circ \varphi)(\overline{z}), \\ 
(z_0, z_1, z_2, \cdots, z_{k-1}, z_{k}) \mapsto (\overline{z}_0, \overline{z}_1, \overline{z}_2, \cdots, \overline{z}_{k-1}, \overline{z}_{k}).
\end{cases}
\end{equation}
Note that $\alpha$ does \emph{not} change under the map $\tau_*$ since both the complex conjugate and the involution swap the immersed sectors. By pulling back the class of (relative) spin structure via $\tau$ and by taking  perturbations respecting $\tau$, the orientation changes as follows, see \cite[Section 3 and 4]{FOOOanti}, which analyzed the signs of the moduli spaces of bordered stable maps bounded by a smooth Lagrangian submanifold under the $\Z/2$-action given by the anti-symplectic involution $\tau$.

\begin{lemma}\label{lemma_antiorderpre}
As a map between oriented spaces, the map $\tau_*$ in~\eqref{equ_taustart} is described by
$$
\tau_* \left(\mcal{M}_{k+1} (\alpha, \beta_0) \right) = (-1)^{\spadesuit_{\beta,k}} \mcal{M}^\textup{clock}_{k+1} ({\alpha}, \beta_0).
$$
where $\spadesuit_{\beta,k} := \frac{1}{2} \mu(\beta) + k + 1$. 
\end{lemma}

\begin{lemma}[Lemma 8.4.3 in \cite{FOOO}]\label{lemma_swapingmarkings}
Let $\sigma$ be the transposition swapping $i$ with $i+1$. 
With orientations taken into account,
The action of $\sigma$ on the moduli space $\mcal{M}_{k+1}(\beta, P_1, \cdots, P_k)$ by changing the order of marked points is described by
$$
\sigma \left( \mcal{M}_{k+1}(\beta, P_1, \cdots, P_i, P_{i+1}, \cdots P_k) \right) = (-1)^{\clubsuit_{i,i+1}} \mcal{M}_{k+1}(\beta, P_1, \cdots, P_{i+1}, P_{i}, \cdots P_k) 
$$
where $\clubsuit_{i,j} := (\deg P_i + 1)(\deg P_{j} + 1)$.
\end{lemma}

Having these lemmas in mind, Lemma~\ref{lem:L0_unobs} can be proved as follows.

\begin{proof}[Proof of Lemma~\ref{lem:L0_unobs}]
As explained above, it only remains to analyze the moduli spaces when
$$
\alpha \colon \{0, 1, \cdots, 2k\} \to R \mbox{ is given by } 
\alpha(j) = 
\begin{cases}
0 &\quad \mbox{if $j = 0$}, \\
U &\quad \mbox{if $j = 1, 3, \cdots, 2k - 1$}, \\
V &\quad \mbox{if $j = 2, 4, \cdots, 2k$},
\end{cases}
$$
and $\beta$ is the trivial class $\beta_0$. Note that $\alpha$ is preserved under the map $\tau_*$ since the complex conjugate and the involution swap the immersed sectors. 
However, since $\tau_*$ does not respect clockwise orientation for $k \geq 1$, we have to rearrange marked points by the map
$$
\sigma (0, 1, 2, \cdots, 2k-1, 2k) \mapsto (0, 2k, 2k-1,  \cdots, 2, 1),
$$ 
where $\sigma \in \frak{S}_{2k+1}$. As a result, we obtain a map $\tau^\textup{main}_* \colon \mcal{M}_{2k+1} ({\alpha}, \beta_0) \to \mcal{M}_{2k+1} (\sigma_* {(\alpha)}, \beta_0)$ defined by
$$
\begin{cases}
\tau_* (\varphi) (z) :=  (\tau \circ \varphi)(\overline{z}), \\ 
(z_0, z_1, z_2, \cdots, z_{2k-1}, z_{2k}) \mapsto (\overline{z}_0, \overline{z}_{2k}, \overline{z}_{2k-1}, \cdots, \overline{z}_2, \overline{z}_1),
\end{cases}
$$
which induces
$$
\tau^\textup{main}_* \colon \mcal{M}_{2k+1} ({\alpha}, \beta_0; U, V, \cdots, U, V) \to \mcal{M}_{2k+1} ({\sigma}_* {(\alpha)}, \beta_0; V, U, \cdots, V, U).
$$
By Lemma~\ref{lemma_antiorderpre} and Lemma~\ref{lemma_swapingmarkings}, the map $\tau^\textup{main}_*$ is orientation reversing because
$$
\spadesuit_{\beta_0, 2k} + \sum_{1 \leq i < j} \clubsuit_{i,j} \equiv 1 \quad (\textup{mod} \, 2).
$$

Recall that the orientation of a moduli space of pearl trajectories is determined by that of the fiber product between moduli spaces of disks and (un)stable submanifolds of the chosen Morse function. 
Thus, with the fixed Morse function on $\mathbb{S}^2$, the relation on the orientations of the moduli spaces of constant disks induces the same relation on those of the corresponding moduli spaces of pearl trajectories. Therefore,  
$$
\begin{cases}
m_{2k, \beta_0}(U,V, \cdots, U, V) +  m_{2k, \tau_* \beta_0}(V, U, \cdots, V, U) = 0,\\ 
m_{2k, \beta_0}(V, U, \cdots, V, U) + m_{2k, \tau_* \beta_0}(U,V, \cdots, U, V) = 0
\end{cases}
$$
and hence
\begin{align*}
m^{b_{\mathbb{L}_0}}_0(1) &= m^{b_{\mathbb{L}_0}}_{0, \beta_0}(1) = \frac{1}{2} \left( m^{b_{\mathbb{L}_0}}_{0, \beta_0}(1) + m^{b_{\mathbb{L}_0}}_{0, \tau_* \beta_0}(1) \right) = 0.
\end{align*}
\end{proof}

For instance, the constant triangles contribute as
$$ 
\begin{cases}
m_2 (uU,vV) = vu \cdot \beta_2^{\mathbb{L}_0},\\ 
m_2(vV,uU) = -uv \cdot \beta_2^{\mathbb{L}_0},
\end{cases}
$$
and hence they cancel with each other. 
In fact, the $A_\infty$-algebra $CF(\BL_0,\BL_0)$ was shown to be quasi-isomorphic to an exterior algebra with two generators in \cite[Section 11]{Seinote}.

In our situation, one can extend the deformation space $\Lambda_+ \times \Lambda_+$ for $\mathbb{L}_0$ to $\Lambda_0 \times \Lambda_+ \cup \Lambda_+ \times \Lambda_0$.
It is because $\BL_0$ only bounds constant holomorphic polygons whose corners must be in an alternating pattern of $U,V$ as we have seen in the proof of Lemma~\ref{lem:L0_unobs}. Thus, $uv$ still lies in $\Lambda_+$ even if one of $u,v$ is in $\Lambda_0$. Therefore
we obtain the following lemma.

\begin{lemma}
	The $A_\infty$-operations on Floer cochain complexes involving $(\BL_0,b_{\BL_0})$ converge for any $(u,v) \in \Lambda_0 \times \Lambda_+ \cup \Lambda_+ \times \Lambda_0$.  
\end{lemma}

By the same reason, ($A_\infty$-structures on) Floer complexes between $(\BL_0,b_{\BL_0})$ and other objects in the Fukaya category are well-defined for $(u,v) \in \Lambda_0 \times \Lambda_+ \cup \Lambda_+ \times \Lambda_0$.

\subsection{Gluing of local charts $\mathcal{MC}(\BL_i)$}\label{subsec:glueform1}
From the discussion in the previous section, we see that the Maurer-Cartan deformation of  $\mathbb{L}_i$ gives rise to the following local charts for the mirror space
\begin{equation*}
\mathcal{MC} (\mathbb{L}_i) = \left\{
\begin{array}{ll}
\Lambda_{\rm{U}} \times \Lambda_{\rm{U}} & \mbox{for} \,\, i=1,2 \\
\Lambda_0 \times \Lambda_+ \cup \Lambda_+ \times \Lambda_0 & \mbox{for} \,\, i=0
\end{array}\right. .
\end{equation*} 
We will use the coordinates $(x_i,y_i)$ for $\mathcal{MC}(\mathbb{L}_i)$ $i=1,2$, and $(u,v)$ for for $\mathcal{MC}(\mathbb{L}_0)$.
We next explain how to glue these charts to give a global mirror with help of isomorphisms among these three objects.

\begin{figure}[h]
	\begin{center}
		\includegraphics[scale=0.37]{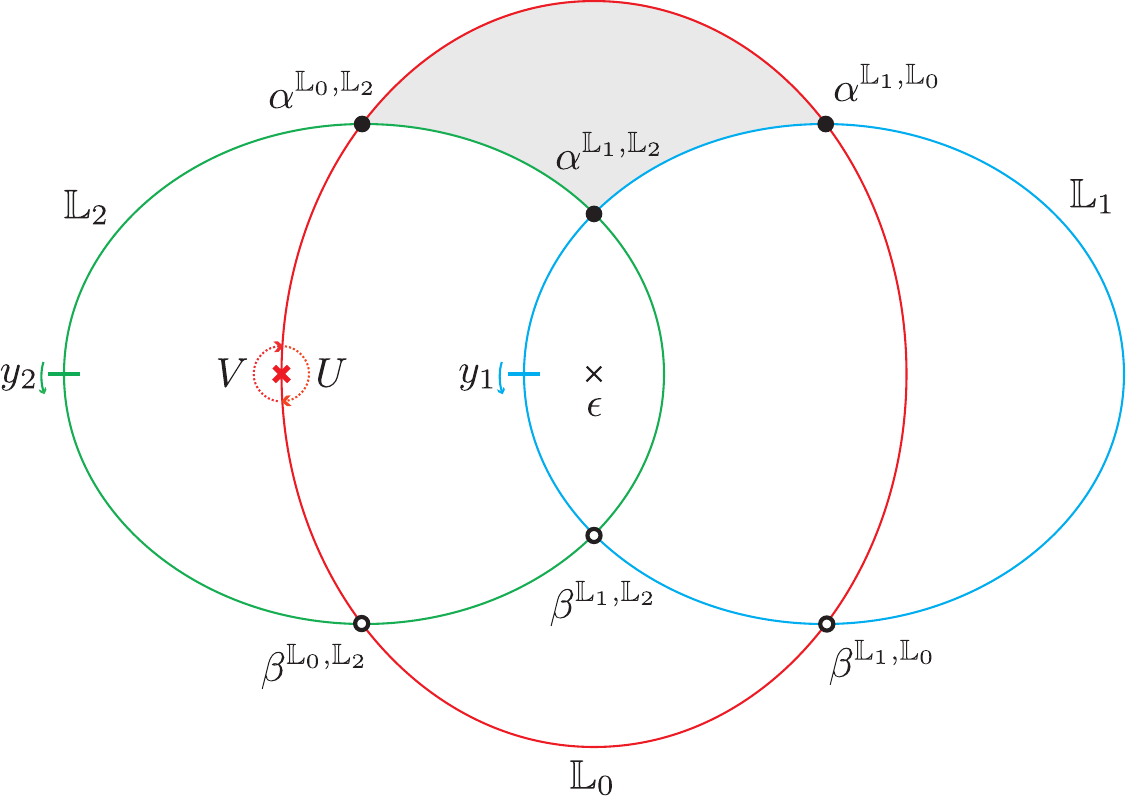}
		\caption{Base paths for $\BL_1$,$\BL_2$ and $\BL_0$ and their intersections}
		\label{fig:gluing}
	\end{center}
\end{figure}

Recall that the Lagrangian tori $\BL_1$, $\BL_2$ and the immersed sphere $\BL_0$ intersect each other cleanly along two circles, and we have chosen generic (perfect) Morse functions on these circles.
Critical points of these functions give generators of $CF(\BL_i,\BL_j)$, which we denote by
$$
\alpha_0^{\BL_i,\BL_j}, \alpha_1^{\BL_i,\BL_j},\beta_1^{\BL_i,\BL_j},\beta_2^{\BL_i,\BL_j} \quad \mbox{for $(i,j)$ with $0 \leq i \neq  j \leq 2$,}
$$
Here, $\alpha_\bullet^{\BL_i,\BL_j}$ and $\beta_\bullet^{\BL_i,\BL_j}$ are respectively over the points $\alpha^{\BL_i,\BL_j}$ and $\beta^{\BL_i,\BL_j}$ in the base of the conic fibration (see Figure~\ref{fig:gluing}), and the subscripts indicate the degrees of the morphisms. For instance, $\alpha_0^{\mathbb{L}_i,\mathbb{L}_j}$ is the maximum of the Morse function on the circle lying over the point $\alpha$. The degrees of $\beta_i$'s are shifted further by $1$ from their usual Morse indices due to the degree of the intersection of base paths at $\beta$.

We are now ready to state the main statement of this section.

\begin{theorem}\label{thm:q-iso} 
	There are quasi-isomorphisms among $(\BL_i, b_{\BL_i})$ and $(\BL_0, b_{\BL_0})$ if we define coordinate transitions as follows$\colon$ 	\begin{itemize}
		\item $(\BL_1,b_{\BL_1}) \cong (\BL_0, b_{\BL_0})$ if and only if
		\begin{equation}\label{equ_standardwallcross3}
			x_1= uv-1,\quad y_1= u 
		\end{equation} 
		where $u,y \in \Lambda_{\rm U}$, $v \in \Lambda_+$, and $x_1 \in -1 + \Lambda_+$.
		\item $(\BL_2, b_{\BL_2}) \cong (\BL_0, b_{\BL_0})$ if and only if
		\begin{equation}\label{equ_standardwallcross2}
			x_2= uv-1,\quad y_2= v^{-1} 
		\end{equation}
		where $v,y_2 \in \Lambda_{\rm U}$, $u \in \Lambda_+$, and $x_2 \in -1 + \Lambda_+$.
		\item $(\BL_1, b_{\BL_1})\cong (\BL_2, b_{\BL_2})$ if and only if
		\begin{equation}\label{equ_standardwallcross}
		x_1= x_2,\quad y_1 = y_2 (1+x_2)
		\end{equation}
		where $y_i \in \Lambda_{\rm U}$ and $x_i \in k + \Lambda_+$ for $k \in \C^\times - \{-1\}$.
	\end{itemize}
\end{theorem}

The coordinate changes above \emph{do} depend on the choice of Morse functions and hypertori for flat connections.  Different choices affects the formula by a certain automorphism on an individual chart, which we shall explain at the end of the section.

\begin{remark}
	The relation between $\BL_1$ and $\BL_2$ has been studied in Seidel's lecture notes \cite{Seinote} and also used by the work of Pascaleff-Tonkonog \cite{PT}.  Here we focus on the relation between $\BL_0$ and $\BL_1$ (and that between $\BL_0$ and $\BL_2$ can be similarly derived).
\end{remark}

\begin{proof}
	
	We compute $m_1^{\textbf{b}_\mathbb{L}}:= m_1^{b_{\BL_1}, b_{\BL_0}}$ of the unique degree-zero generator $\alpha_0^{\BL_1,\BL_0}$ in the Floer complex  $\mathrm{CF}((\BL_1,b_{\BL_1}),(\BL_0, b_{\BL_0}))$
	by counting pearl trajectories.  There are two strips from $\alpha_0^{\BL_1,\BL_0}$ to $\beta_1^{\BL_1,\BL_0}$ whose projections under $\Pi$ are the shaded regions in Figure \ref{fig:dalpha} (a). These strips can be found explicitly as holomorphic sections of the conic fibration. Once we delete $\{a=0\}$ or $\{b=0\}$, the fibration becomes a trivial $\C^\times$-fibration where one can realize these strips as constant sections of the trivial fibration. See for e.g. \cite[Section 17]{Seibook} for detailed analysis on such sections. Let us denote these strips by $u_L$ and $u_R$, respectively. 
	
For our Morse function on $\BL_0$, we require that the flows starting at the self-intersection point avoid the arcs $\partial_{\BL_0} u_L$ and $\partial_{\BL_0} u_R$. 
Similarly, we choose $x_1$-hypertorus for $\BL_1$ such that it avoids $\partial_{\BL_1} u_L$ and $\partial_{\BL_1} u_R$. For other coordinate changes, one has to choose hypertori in $\BL_2$ satisfying similar properties. See the discussion at the end of the section for more details, in particular for the existence of such hypertori and Morse functions.
	
	Under this choice, we see that $u_R$ contributes without involving any variables $u,v,x_1,y_1$ since neither $\partial_{\BL_0} u_R$ can be joined by the constant disk at the self-intersection point through a Morse flow nor $\partial_{\BL_1} u_L$ intersects the $x$-hypertorus. (It does not hit the $y$-hypertorus, either, due to our choice as in Figure \ref{fig:dalpha}.) Similarly, $u_L$ gives rise to $uy_1$ (without any $x$ or $x_1$ appearing) in its contribution since it has a corner at $U$ and passes through the $y_1$-hypertorus, but no other than these.

	Recall that we have chosen $\BL_i$'s such that the two strips have the same symplectic area, say, $\Delta$. From the above discussion we conclude that
	\begin{equation}\label{eqn:l1im}
	\langle m_1^{\textbf{b}_\mathbb{L}}(\alpha_0^{\BL_1,\BL_0}),\beta_1^{\BL_1,\BL_0} \rangle = T^{\Delta}\left( 1-  u y_1^{-1} \right).
	\end{equation}
	where the strip $u_R$ contributes $T^\Delta$, while $u_L$ contributes $T^\Delta uy_1^{-1}$. We will give a brief explanation how different choices of Morse functions affect the computation later.	
	
	\begin{figure}[h]
		\begin{center}
			\includegraphics[scale=0.37]{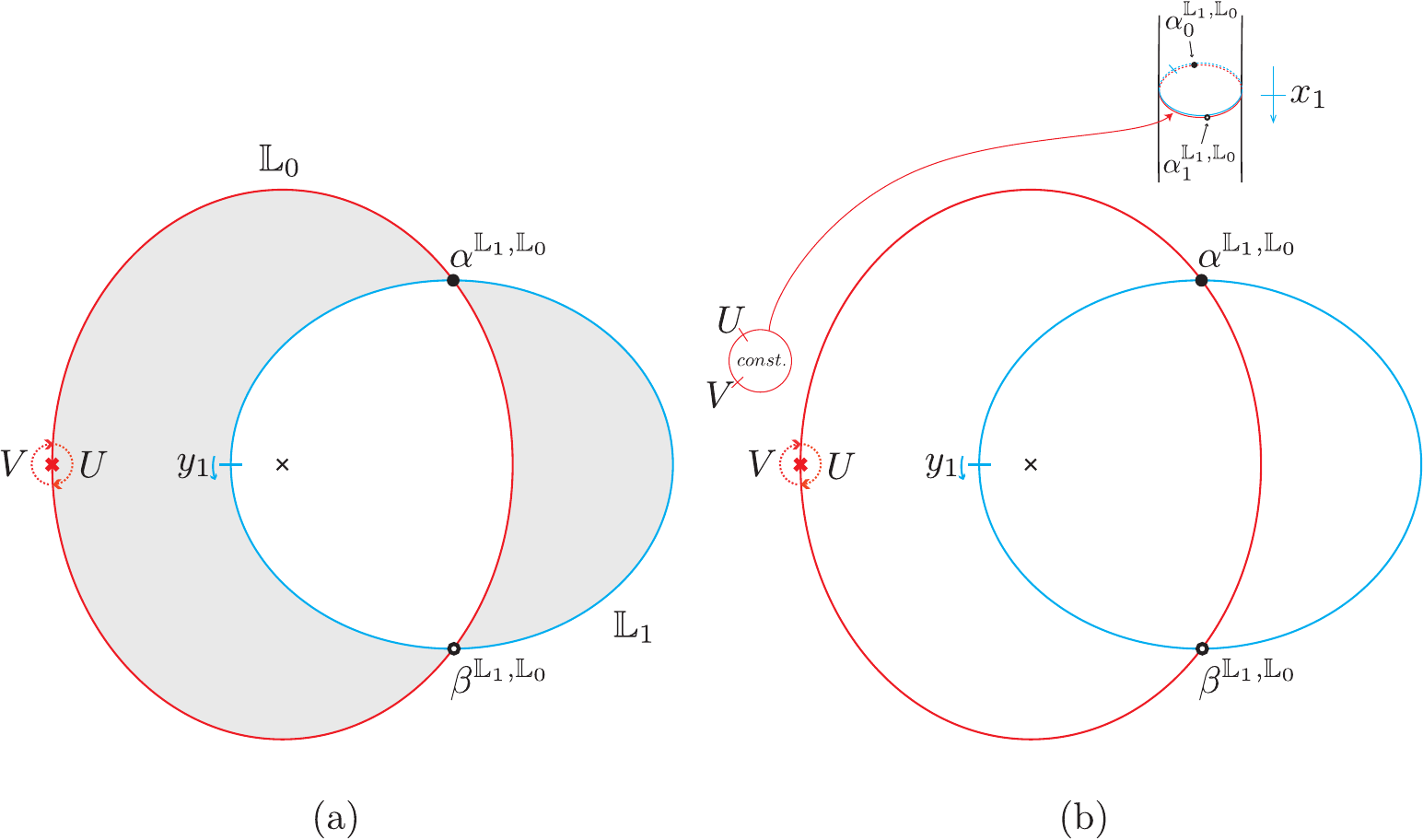}
			\caption{Pearl trajectories contributing to $m_1^{\textbf{b}_\mathbb{L}} (\alpha_0^{\BL_1,\BL_0}  )$}
			\label{fig:dalpha}
		\end{center}
	\end{figure}
	
	For $\left\langle m_1^{\textbf{b}_\mathbb{L}}(\alpha_0^{\BL_1,\BL_0}), \alpha_1^{\BL_1,\BL_0} \right\rangle$, we also have a pair of strips in $\Pi^{-1} (\alpha^{\BL_1,\BL_0})$ contributing to it.  One of the strips gives $\pm 1$, and the other gives $\pm x_1$ since one passes through the $x$-hypertorus.   In addition to those strips, we have pearl trajectories consisting of a strip, Morse flow lines in $\bL_0$ and constant polygons with corners $U,V,\ldots,U,V$.  The pearl trajectories form power series in $uv$.  Hence
	$$\left\langle m_1^{\textbf{b}_\mathbb{L}}(\alpha_0^{\BL_1,\BL_0}), \alpha_1^{\BL_1,\BL_0} \right\rangle = x_1 h(uv) + g(uv)$$
	where $g,h$ are power series with leading term $\pm 1$.  In summary
	$$ m_1^{\textbf{b}_\mathbb{L}}(\alpha_0^{\BL_1,\BL_0}) = T^{\Delta}\left( 1- uy_1^{-1} \right)\beta_1^{\BL_1,\BL_0} + (x_1 h(uv) + g(uv))\, \alpha_1^{\BL_1,\BL_0}.$$

	We see that the cocycle condition (i.e. $m_1^{\textbf{b}_\mathbb{L}}$-closedness) gives
	\begin{equation}\label{eqn:coord1im}
	\left\{
	\begin{array}{l}
	x_1 + H(uv) = 0\\
	y_1 = u
	\end{array}\right.
	\end{equation}
	where $H:=g/h$.
	One can check that the degree zero morphism in $\mathrm{CF}(\BL_0,\BL_1)$ (which is dual to $\beta_2^{\BL_1,\BL_0}$) is also closed under the corresponding Floer differential if and only if the same condition is satisfied. Moreover, one can check by similar counting that  $\alpha_0^{\BL_1,\BL}$ gives an isomorphism under the condition \eqref{eqn:coord1im} with its inverse being the dual of $\beta_2^{\BL_1,\BL_0}$.
	
	However, it is \emph{not} easy to compute the function $H(uv)$ directly as it necessarily involves a virtual perturbation due to constant bubbles and multiple covers. To compute $H(uv)$, we examine relations among three objects $\BL_0,\BL_1$, and $\BL_2$.  Namely, we consider the chain of isomorphisms
	$$(\BL_1, b_{\BL_1}) \stackrel{\alpha_0^{\BL_1,\BL_0}}{\longrightarrow} (\BL_0, b_{\BL_0}) \stackrel{\alpha_0^{\BL_0,\BL_2}}{\longrightarrow} (\BL_2, b_{\BL_2}).$$
	The composition of these two isomorphisms is given by
	$$m_2( \alpha_0^{\BL_1,\BL_0}, \alpha_0^{\BL_0,\BL_2}) = T^{\Delta'} \alpha_0^{\BL_1,\BL_2} $$
	where $\Delta'$ is the symplectic area of the triangle projecting down to the shaded region in Figure \ref{fig:gluing}, and we make similar assumption on our Morse function as before so that the flow lines from the self-intersection do not intersect the boundary of the triangle lying over $\BL_0$. This implies that $\alpha_0^{\BL_1,\BL_2}$ is also an isomorphism. Therefore, coordinate changes among three objects must be compatible (over the triple intersection). Using the relation
	$$x_1=x_2,\quad y_1 = y_2 (x_2+1)$$
	we derive
	$$ 1-H(uv) = x_2+1 =  \frac{y_1}{y_2} = \frac{u}{v^{-1}} = uv.$$
\end{proof}

\begin{figure}[h]
	\begin{center}
		\includegraphics[scale=0.45]{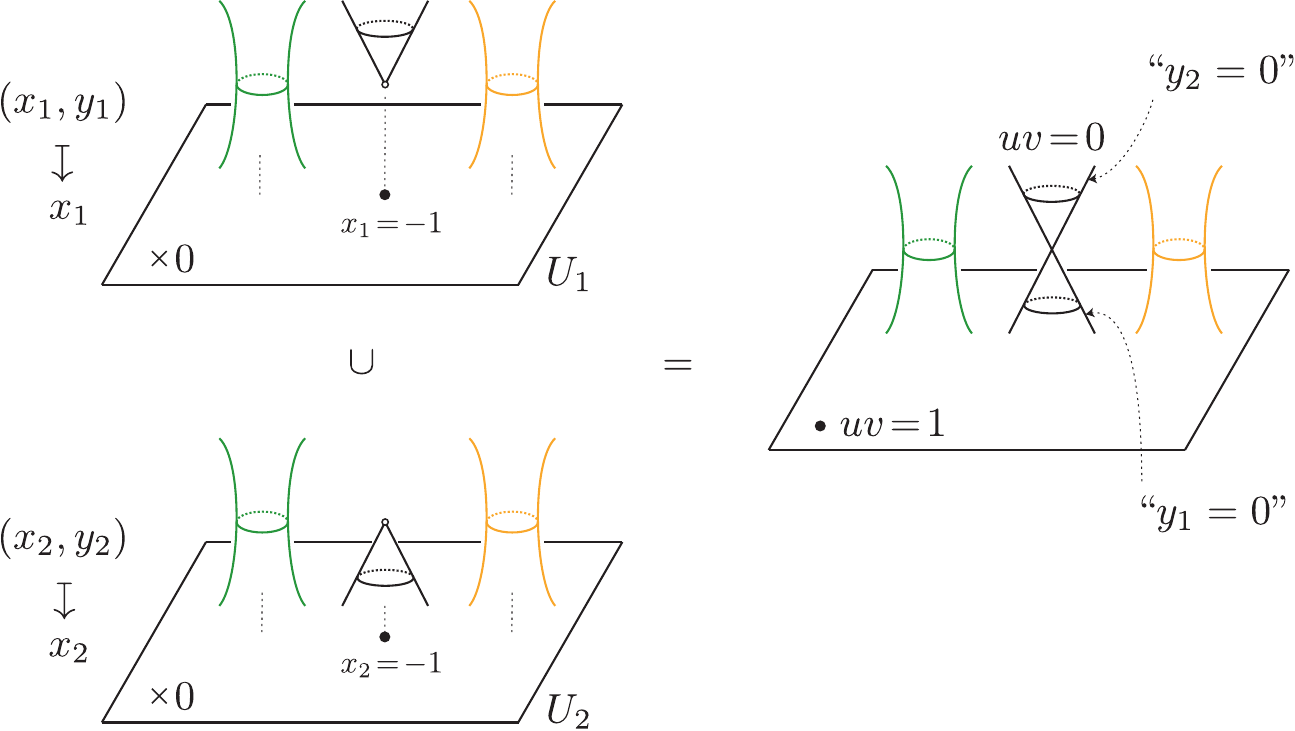}
		\caption{Illustration of gluing of two charts}
		\label{fig:glcharts}
	\end{center}
\end{figure}

The three deformation spaces are glued accordingly.  Namely we take the disjoint union of $(\Lambda_0 \times \Lambda_+) \cup (\Lambda_+ \times \Lambda_0)$ (with coordinates $(u,v)$) and two copies of $\Lambda_{\rm U} \times \Lambda_{\rm U}$ (with coordinates $(x_1,y_1)$ and $(x_2,y_2)$ respectively), and take the quotient according to the above relations. 
The gluing is illustrated in Figure \ref{fig:glcharts}.  The upper component of the degenerate cone $\{uv=0\}$ is a (partial) compactification of $\{x_1=-1\} \cong \C^\times$ by adding the point $\{x_1=-1, y_1=0\}$, and similar for the lower component.  The resulting moduli space is
$$\{(u,v) \in \Lambda_0 \times \Lambda_0: uv \not\in 1 + \Lambda_+\}.$$

In the proof, we made the simplest choice of our Morse function, for which the flow lines through the self-intersection have no intersections with contributing holomorphic polygons. If we choose arbitrary Morse function, intersection between the flows lines and holomorphic polygons will create more pearl trajectory contributions (see Section~\ref{subsec:interpretps} for instance). 
In fact one can check by following the same argument in the proof, but with different Morse function that the wall crossing formula can get affected by multiplying $(1-uv)^k$. Also, choosing different $x_i$-hypertorus can create additional $x_i^l$ in the wall crossing formula. However, the resulting mirror space will remain the same up to an equivalence, and essentially this change amounts to applying the following automorphism to $\mathcal{MC} (\BL_0) = \{ (u,v) \}$ and gluing:
$$ (u,v) \mapsto \left( (1-uv)^k u, (1-uv)^{-k} v \right).$$
See Remark \ref{rmk:intrinsicmc} for the related discussion, especially for a more intrinsic description of the Maurer-Cartan deformation space in this case. 

The choice we made here can also be interpreted as that of directions in SYZ base, along which we smooth out $\mathbb{L}_0$. We will give a brief explanation on this point at the end of the section.

\subsection{Choices of hypertori and Morse functions}\label{subsec:supple}
We supplement the proof of Theorem \ref{thm:q-iso}, describing precisely the hypertori and the Morse functions we have chosen in the proof.

	\subsubsection*{Choice of hypertori}
	
	The $y_i$-hypertorus (or $y_i$-circle) for $\bL_i$ for $i=1,2$ are canonical, namely they are taken to be vanishing circles.  However for the $x_i$-hypertorus we need to make a choice.  If we make an arbitrary choice, the wall-crossing equation would be $y'_1=x_2^k y_2 (1+ x_2)$ for some $k$.  In the following we make a suitable choice such that $k=0$.
	
	The conic fibration is non-trivial over the disk bounded by the base circle of $\bL_2$.  Let's trivialize the conic fibration by deleting $\{b=0\}\subset\C^2 - \{ab=\varepsilon\}$.   Then the $x_i$-circle of $\bL_i$ is simply taken to be a trivial section of the base circle.  
	
	Recall that we have fixed a perfect Morse function on each of the two clean intersections between $\bL_1$ and $\bL_2$ (which are circles), and denoted by $\alpha_0^{\BL_1,\BL_2}$, $\beta_1^{\BL_1,\BL_2}$ the maximum point in the upper clean intersection, and the minimum point in the lower clean intersection respectively.  Consider the holomorphic sections bounded by $(\bL_1,\bL_2)$.  There is a unique section $s_0$ passing through $\beta_1^{\BL_1,\BL_2}$ over the right region.  Over the left region, there are two sections $s_1,s_2$ passing through $\beta_1^{\BL_1,\BL_2}$, and exactly one of them, say $s_1$, is contained in the trivialization.  
	
	By considering their classes, we see that $s_0$ contributes $x_1^a (x_2)^{-a} y_2 (y_1)^{-1}$ for some $a$, and $s_1$ contributes $(x_2)^b x_1^{-b}$ for some $b$. Morse flow line in clean intersection contributes to $m_1^{\textbf{b}_\mathbb{L}}  (\alpha_0^{\BL_1,\BL_2})$ by $(x_1 - x_2) \alpha_1^{\BL_1,\BL_2}$, so we impose the relation $x_1=x_2$ (for two Lagrangians to be isomorphic). Under this relations, $s_0$ and $s_1$ give the terms $y_2 (y_1)^{-1}$ and $1$ respectively. That is, these two strips contribute to the $\beta_1^{\BL_1,\BL_2}$-component of $m_1^{\textbf{b}_\mathbb{L}}  (\alpha_0^{\BL_1,\BL_2})$ as $y_2 (y_1)^{-1}$ and $1$.
	
	For the section $s_2$, its boundary class equals to the sum of the boundary class of $s_1$ and the vanishing circle.  From the above the boundary of $s_1$ has trivial holonomy.  The vanishing circle intersects gauge $x_2$-circle once and has holonomy $x_2$. Hence the strip $s_2$ contributes $x_2$ to the $\beta_1^{\BL_1,\BL_2}$-component of $m_1^{\textbf{b}_\mathbb{L}}  (\alpha_0^{\BL_1,\BL_2})$. 
	  In total we have the wall-crossing formula $y_2 (y_1)^{-1}=1+x_2 = 1+x_1$.
	        
	\subsubsection*{Choice of Morse function of the immersed sphere}
	
		Now consider the strips bounded by $(\bL_0,\bL_i)$ for $i=1,2$.  Again we choose perfect Morse functions on the clean intersections between $\bL_0$ and $\bL_i$.  Denote by $\alpha_{0}^{\BL_0,\BL_i}$, $\beta_{1}^{\BL_0,\BL_i} $ the maximum point in the upper clean intersection, and the minimum point in the lower clean intersection respectively.
		
		Below we choose the Morse function on the normalization $S\cong\bS^2$ of $\bL_0$ such that the flow lines from $q_1,q_2$ (preimages of the immersed point) to the minimum point do not intersect the boundary of the strips in the complement of $q_1,q_2$.   This ensures that there is no constant polygons at the immersed point contributing to the $\beta_{1}^{\BL_0,\BL_i} $-component of $m_1^{\textbf{b}_\mathbb{L}}  (\alpha_0^{\BL_0,\BL_i})$, so that we have the gluing formula $u=y_1$, $v=y_2$.  (Note that in general the gluing is $h(uv) u = g(uv) y_1$ if the above flow lines intersect with boundary of strips, where $h$ and $g$ are certain series contributed from the constant polygons.)
	
		The sections bounded by $(\bL_0,\bL_i)$ can be explicitly solved by taking square root.  Namely we can set $a=e^{i\theta}b$ and solve $ab=\zeta$ where $\zeta$ is the domain variable for the base strip (on the left or on the right).   Requiring each section to pass through the minimum point of the lower clean intersection gives a unique solution.
		
		For the sections bounded by $(\bL_0,\bL_i)$ over the right base strip (which does not pass through $ab=0$), their boundaries in $\bL_0$ coincide with each other, giving a curve segment in $S$ connecting the two clean intersections.
		
		For the sections bounded by $(\bL_0,\bL_i)$ ($i=1,2$) over the left base strip, by considering the different branch cuts in taking square root, we see that their boundaries in $\bL_0$ are parallel to each other.  Namely the argument of $u$ for the two sections differs by $\pi$.  Their boundaries give two curve segments in $S$ connecting $q_i$ to the clean intersections.  They do not intersect with each other (in the complement of $q_1,q_2$).
		
\begin{figure}[h]
			\begin{center}
				\includegraphics[scale=0.37]{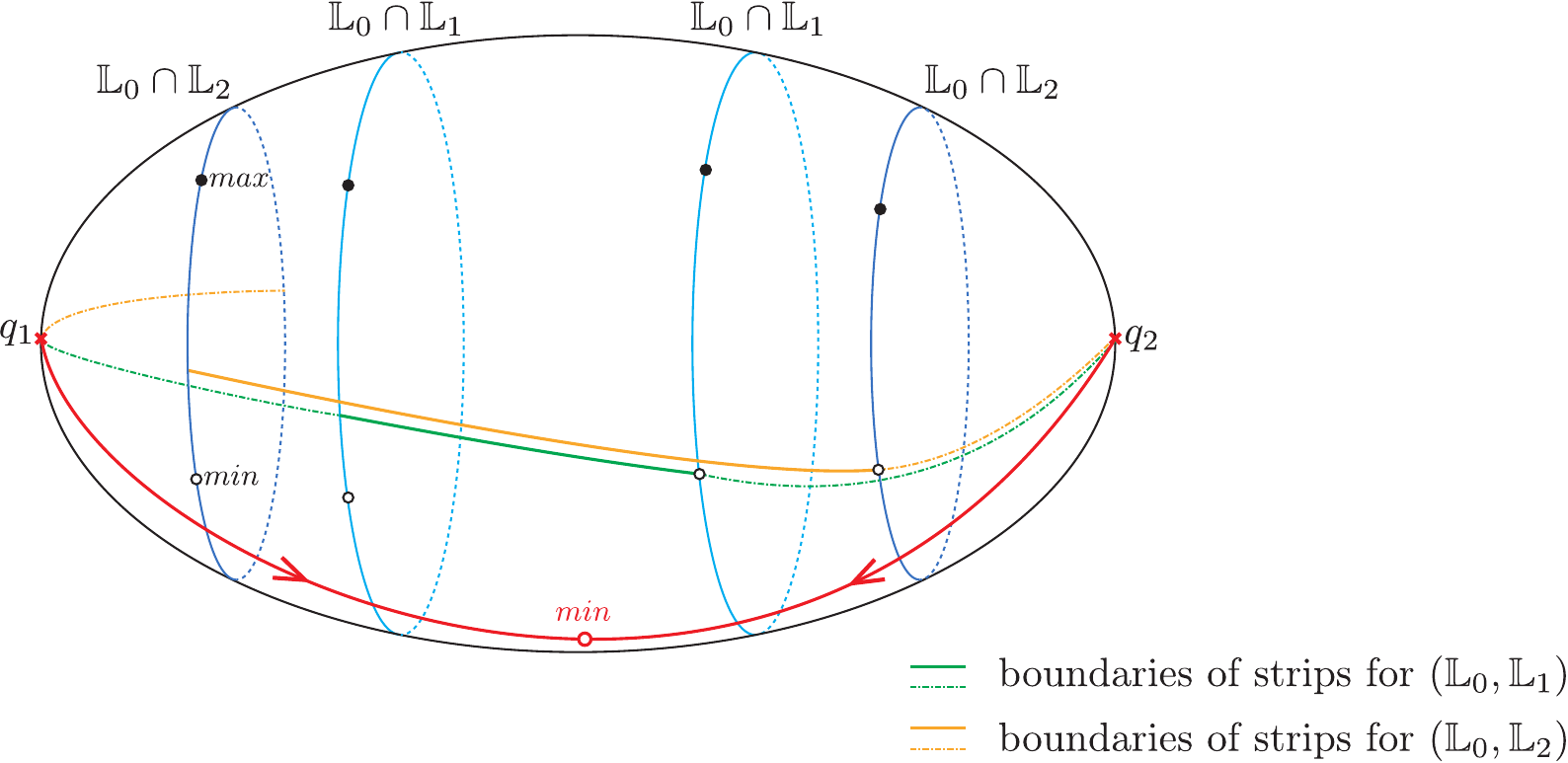}
				\caption{Choosing a Morse function for the gluing formula.}
				\label{Fig_Morse}
			\end{center}
		\end{figure}

		From the above configuration, it follows that there exists an arc in $S$ connecting $q_1,q_2$ which does not intersect with any of the above curve segments.  We take a Morse function on $S$ such that the union of the flow lines from $q_i$ to the minimum is such an arc.  See Figure \ref{Fig_Morse}. 
		
\vspace{-0.3cm}

	\subsubsection*{Relation with choice of smoothing}
	
\begin{figure}[h]
			\begin{center}
				\includegraphics[scale=0.55]{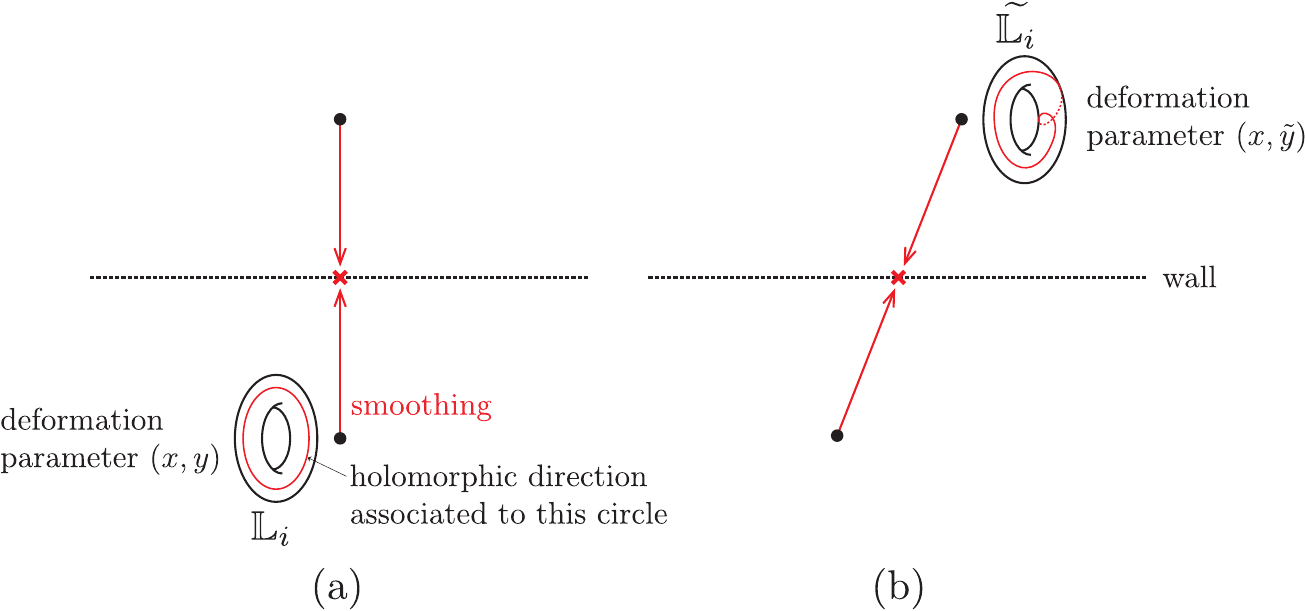}
				\caption{Two different choices of smoothings}
				\label{Fig_choicesm}
			\end{center}
		\end{figure}

	We have made a choice of smoothing $\mathbb{L}_i$ of $\mathbb{L}_0$ using the local model. In the SYZ fibration picture, it can be understood as smoothing the singular fiber along the direction perpendicular to the wall (as in (a) of Figure \ref{Fig_choicesm}). However, there is no preferred way of smoothing since we do not have a preferred complex structure in symplectic geometry. 
If we choose another directions of smoothing $\widetilde{\mathbb{L}_i}$, this gives rise to different cycles in $\mathbb{L}_i$, which we take as a new $x_i$-circle as in (b) of Figure \ref{Fig_choicesm}. 
Here, $\widetilde{\mathbb{L}_i}$ and $\mathbb{L}_i$ are isotopic through  a Lagrangian isotopy which does not cross the wall.
By the same argument as above, this induces the same wall-crossing formula, yet with respect to  new coordinates,
	$$\tilde{y_1} = \tilde{y_2} (1+x_2)$$
	where $\tilde{y_1} = x_1^k y_1$ and $\tilde{y_2} = x_2^k y_2'$ (the relation $x_2 = x_1$ remains the same). 
 
For the immersed sphere $\mathbb{L}_0$, take Morse arcs such that they do not intersect strip boundaries for $\mathbb{L}_i$ or $\widetilde{\mathbb{L}_i}$. These two different choices of Morse arcs give
$$ u=y,\quad \mbox{and} \quad \tilde{u} = \tilde{y}$$
respectively. Thus we see that  they are related by 
$$ \tilde{u} = x_1^k u = (uv -1)^k u .$$

\section{Review of Grassmannians}

In this section, we briefly recall some features of $\Gr(2,\C^n)$.  We shall construct their mirrors in Section \ref{sec_LagfibonF3} and \ref{section_mirrorcongr2n}.

\subsection{Grassmannians of 2-planes}\label{subsec:vardefgr}

The Grassmannian $\mathrm{Gr}(2, \C^n)$ of 2-planes in $\C^n$ is a partial flag manifold parametrizing the complex subspaces of two dimension in $\C^n$. Below are a few different descriptions of $\mathrm{Gr}(2, \C^n)$.

\vspace{0.3cm}

\noindent {\bf (i)} Let $V_2( \C^n)$ be the set of orthonormal frames of two dimensional subspaces in $\C^n$. An element in $V_2( \C^n)$ can be regarded as the set of columns of $\Xi$ in 
\begin{equation}\label{equ_framess}
\left\{
\Xi := 
\begin{bmatrix}
    \zeta_{1} & \xi_{1}  \\
    \zeta_{2} & \xi_{2}  \\
    \vdots & \vdots  \\
    \zeta_{n} & \xi_{n} 
\end{bmatrix}
\in M_{n \times 2}(\C)
~:~ 
\overline{\Xi}^T \Xi = I_{2}
\right\}.
\end{equation}
The Grassmannian is then the set of subspaces spanned by frames in $V_2( \C^n)$, and hence $\mathrm{Gr}(2, \C^n) \simeq V_2(\C^n) / \mathrm{U}(2)$. 

\vspace{0.2cm}

\noindent {\bf (ii)}
The unitary group $\mathrm{U}(n)$ acts linearly and transitively on the space of complex two-dimensional subspaces in $\C^n$ with the isotropy subgroup $\mathrm{U}(2) \times \mathrm{U}(n-2)$. $\mathrm{Gr}(2, \C^n)$ is a homogeneous space$\colon$ 
\begin{equation}\label{equ_u2homogeneousgr24}
\mathrm{Gr}(2, \C^n) \simeq {\mathrm{U}(n)}/{(\mathrm{U}(2) \times \mathrm{U}(n-2))}.
\end{equation}
$\mathrm{Gr}(2, \C^n)$ can be also realized as $G/P$ where $G = \mathrm{SL}(n, \C)$ and $P = \{[a_{i,j}]\}$ is a parabolic subgroup consisting of matrices in $G$ such that all entries $a_{i,j}$ with $i \geq 3$ and $j = 1, 2$ vanish. 

\vspace{0.2cm}

\noindent {\bf (iii)}
Setting $Z_{i,j} = \zeta_i \xi_j - \zeta_j \xi_i$ for $\zeta_i$ and $\xi_j$ in \eqref{equ_framess} defines an embedding of $\mathrm{Gr}(2, \C^n)$ into $\CP^{n(n-1)/2 - 1}$. $\mathrm{Gr}(2, \C^n)$ is then the subvariety of $\CP^{n(n-1)/2 - 1}$ cut out by the Pl{\"u}cker relations 
$$\left\{ (i,j,k,\ell) = 0 ~\colon~ 1 \leq i < j < k < \ell \leq n \right\}$$
where 
\begin{equation}\label{equ_pluckker}
(i,j,k,\ell) := Z_{i,j} Z_{k,\ell} - Z_{i,k} Z_{j,\ell} + Z_{i,\ell} Z_{j,k}.
\end{equation}

\vspace{0.2cm}

\noindent {\bf (iv)}
The Grassmannian $\mathrm{Gr}(2, \C^n)$ is diffeomorphic to a co-adjoint orbit.
Identifying the dual Lie algebra of $\mathrm{U}(n)$ with the set of $(n \times n)$-Hermitian matrices, the co-adjoint orbit $\mcal{O}_\lambda$ can be defined as the orbit of the diagonal matrix $I_\lambda$ under the conjugate $U(n)$-action where
\begin{equation}\label{equ_choiceofeigen}
\lambda := (n-2, n-2, \underbrace{-2,\cdots,-2}_{(n-2) \mbox{ times}} )
\end{equation}
and 
$$
I_\lambda := \textup{diag}(n-2,n-2,-2, \cdots, -2) \in M_{n \times n}(\C).
$$ 
Note that $\mathrm{U}(n)$ acts transitively on $\mcal{O}_\lambda$ by conjugation and the isotropy group of $I_\lambda$ is $\mathrm{U}(n-2) \times \mathrm{U}(2)$. 
Thus, $\mcal{O}_\lambda \simeq \mathrm{U}(n) / (\mathrm{U}(n-2) \times \mathrm{U}(2))$,
and hence $\mcal{O}_\lambda$ can be identified with $\mathrm{Gr}(2, \C^n)$ by~\eqref{equ_u2homogeneousgr24}. 

It is equipped with a $\mathrm{U}(n)$-invariant symplectic form, which is called a \emph{Kirillov-Kostant-Souriau symplectic form} $\omega_\lambda$.  
Moreover, by Guillemin-Sternberg \cite{GS83}, the Grassmannian $\mathrm{Gr}(2, \C^n)$ equipped with the pull-backed Fubini-Study form on $\CP^{n(n-1)/2 - 1}$ is isomorphic to the co-adjoint orbit $(\mcal{O}_\lambda, \omega_\lambda)$.

\subsection{Toric degenerations and Gelfand-Cetlin toric varieties}\label{Sec_toricdegg}

Although partial flag varieties are \emph{not} toric varieties in general, it is known that they admit degenerations into toric varieties (see \cite{GoLa, Cal02, AlBr} for instance) which are very useful to study partial flag varieties themselves.
For our purpose later, we take one particular toric degeneration of the Grassmannian $\mathrm{Gr}(2,\C^n)$, which appeared in \cite{BCKV, KoMi, NNU}.

To describe the toric degeneration of $\mathrm{Gr}(2, \C^n)$, we view $\mathrm{Gr}(2, \C^n)$ as a subvariety~\eqref{equ_pluckker} of $\CP^{n(n-1)/2 - 1}$. 
The toric degeneration is then defined by 
\begin{equation}\label{equ_toricdege}
\mcal{X}_\varepsilon = \left\{ \left( \left[ Z_{ij} \right], \varepsilon \right) \in \CP^{\frac{n(n-1)}{2} - 1} \times \C ~\colon~ (\varepsilon; i,j,k,\ell) =  \varepsilon Z_{ij} Z_{k\ell} - Z_{ik} Z_{j\ell} + Z_{i\ell} Z_{jk} = 0 \right\}
\end{equation}
The toric variety $\mcal{X}_{\varepsilon = 0}$, defined by binomial relations, is called a \emph{Gelfand-Cetlin (GC) toric variety}. 
Setting $u_{1,n-1} := n-2$ and $u_{3,1} := -2$, 
the associated lattice polytope in $M_\R = \R^{2(n-2)} = \{ (u_{i,j}) ~\colon~ i= 1, 2, \,\, j = 1, \cdots, n-2 \}$ is given by the following set of inequalities$\colon$
\begin{align}
\label{equ_verfirst} &\left\{ u_{1,j+1} - u_{1,j} \geq 0 \right\} \quad \mbox{for $j = 1, \cdots, n-2$}, \\
\label{equ_versec} &\left\{ u_{2,j+1} - u_{2,j} \geq 0 \right\} \quad \mbox{for $j = 1, \cdots, n-3$}, \\
\label{equ_hor} &\left\{ u_{1,j} - u_{2,j} \geq 0 \right\} \quad \mbox{for $j = 1, \cdots, n-2$}, \\
\label{equ_horsec} &\left\{ u_{2,1} - u_{3,1} \geq 0 \right\}.
\end{align}
The polytope is called the \emph{Gelfand-Cetlin (GC) polytope} and denoted by $\Delta_{(2,n)}$.
In order to describe the face structure of $\Delta_{(2,n)}$, it is convenient to use the ladder diagram $\Gamma_{(2,n)}$. 
The ladder diagram $\Gamma_{(2,n)}$ is the induced subgraph in $\R \times \Z \cup \Z \times \R$ whose vertex set is 
$$
V_{(2,n)} = \left\{ (i,j) \in \Z^2 ~\colon~ 0 \leq i \leq 2,\,\, 0 \leq j \leq n-2 \right\}.
$$
Thus the diagram $\Gamma_{(2,n)}$ is a rectangular net of size $(2 \times (n-2))$.

\begin{definition}\label{def_poneadd}
Let $\Gamma_{(2,n)}$ be the ladder diagram associated with $\Delta_{(2,n)}$ as above.
\begin{itemize}
\item A \emph{positive path} is a path with \emph{minimal length} from the bottom-left vertex to the top-right vertex in $\Gamma_{(2,n)}$. 
\item An \emph{admissible diagram} of $\Gamma_{(2,n)}$ is a subgraph that is expressed as a union of positive paths. 
\item The \emph{dimension} of an admissible diagram $\Gamma$ is defined by the number of bounded regions of $\Gamma$.\end{itemize}
\end{definition}

Let $\square^{(i,j)}$ be the unit box whose upper right vertex is located at $(i,j)$. 
Putting $u_{i,j}$ in the box $\square^{(i,j)}$, the inequalities in~\eqref{equ_verfirst},~\eqref{equ_versec},~\eqref{equ_hor} and~\eqref{equ_horsec} can be  respectively assigned to  the edges in $\{(s,j) ~\colon~ 0 \leq s \leq 1\}$, $\{(s,j) ~\colon~ 1 \leq s \leq 2\}$, $\{ (1,s) ~\colon~  j - 1 \leq s \leq j \}$ and $\{ (2,s) ~\colon~ 0 \leq s \leq 1 \}$ in an obvious manner.

Suppose we are given an admissible diagram $\Gamma$ of $\Gamma_{(2,n)}$, and consider the facets of $\Delta_{(2,n)}$ which are obtained by turning the inequalities assigned to edges in $\Gamma_{(2,n)} \setminus \Gamma$ into equalities. Intersecting all such facets, one can associate a unique face of $\Delta_{(2,n)}$ with each admissible diagram $\Gamma$. In fact, this defines an one-to-one correspondence as follows.

\begin{proposition}[\cite{ACK}]\label{prosposition_ACK}
There is an order-preserving one-to-one correspondence 
\begin{equation}\label{equ_facediagramcorrespondence}
\left\{ \mbox{admissible diagrams of $\Gamma_{(2,n)}$} \right\} \to \left\{ \mbox{faces of $\Delta_{(2,n)}$} \right\}
\end{equation}
where the set-theoretical inclusion gives the set of admissible diagrams a partial order.

Moreover, the dimension of an admissible diagram equals that of the corresponding face. 
\end{proposition}

\begin{example}\label{example_gr24diag}
We examine the case where $n = 4$. The ladder diagram $\Gamma_{(2,4)}$ is given as in Figure~\ref{fig_ladder4}. 
The Gelfand-Cetlin polytope $\Delta_{(2,4)}$ has six facets whose corresponding admissible diagrams are $\Gamma(f_1), \cdots, \Gamma(f_6)$, each of which contains three bounded regions. 
For instance, $\Gamma(f_1)$ maps into the facet supported by $u_{1,2} = 2$ under the correspondence~\eqref{equ_facediagramcorrespondence}. The face $g$ corresponding to $\Gamma(g)$ in Figure~\ref{fig_ladder4} is one-dimensional because $\Gamma(g)$ consists of one bounded region. 
It is contained in the affine line defined by $u_{1,1} = u_{1,2} = u_{2,1} = u_{2,2}$. 
Since $\Gamma(g) \subset \Gamma(f_i)$, the face $g$ is contained in the facet $f_i$ for $i = 2,3,4,5$.
Finally, the face $h$ corresponding to $\Gamma(h)$ is the interior of $\Delta_{(2,4)}$. 

\begin{figure}[h]
	\begin{center}
		\includegraphics[scale=0.8]{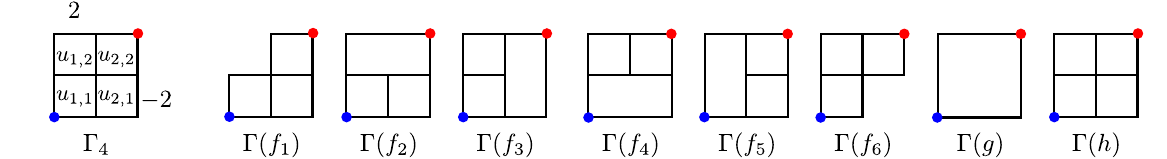}
	\caption{\label{fig_ladder4} The ladder diagram $\Gamma_{(2,4)}$ and its admissible diagrams.}	
	\end{center}
\end{figure}
\end{example}

For each $\varepsilon$ and $1 \leq k \leq n$, the variety $\mcal{X}_\varepsilon$ carries the Hamiltonian $\mathbb{S}^1$-action $\theta_k (\in [0,2\pi]/ 0 \sim 2 \pi)$ defined by
\begin{equation}\label{equ_Taction}
(\theta_k, [Z_{ij}]) \mapsto 
\begin{cases}
[e^{\sqrt{-1} \theta_k} \cdot Z_{ij}] \quad &\mbox{if either $i = k$ or $j = k$,}\\
[Z_{ij}] \quad &\mbox{otherwise}
\end{cases}
\end{equation}
for $k =1, \cdots, n$. 
In other words, $\theta_k$ acts on the $k$-th row $[\zeta_k, \xi_k]$ by rotation and fixes the other rows in~\eqref{equ_framess}.

\subsection{Block combinatorics and Gelfand-Cetlin systems}\label{subsection_blockcombi}

The Gelfand-Cetlin system in a partial flag manifold (of Lie type A, B, and D) was constructed by Guillemin-Sternberg \cite{GS83} by using Thimm's trick.
Regarding $\mathrm{Gr}(2,\C^n)$ as the co-adjoint orbit $\mcal{O}_\lambda$ for a choice of $\lambda$ in~\eqref{equ_choiceofeigen}, the \emph{Gelfand-Cetlin (GC for short in what follows) system} on $\mcal{O}_\lambda$ is defined by
\begin{equation}\label{equ_GCsystem}
\Phi_\textup{GC} = (\Phi_{i,j} \colon i = 1, 2, j = 1, \cdots, n-2) \colon \mcal{O}_\lambda \to \R^{2(n-2)}
\end{equation}
where $\Phi_{i,j}$ is the $i$-th largest eigenvalue of the $(i+j-1)$-th leading principal minor matrix of a Hermitian matrix in $\mcal{O}_\lambda$. 
It is a complete integrable system with respect to $\omega_\lambda$, and forms of action variables. 
The image of the system coincides with the GC polytope $\Delta_{(2,n)}$ determined by the min-max principle (if $\mcal{O}_\lambda$ is equipped with the form $\omega_\lambda$ with a choice of $\lambda$ in~\eqref{equ_choiceofeigen}). The reader is referred to \cite{GS83} for more details. 

Each component generates a \emph{local} $\mathbb{S}^{1}$-Hamiltonian action. 
Specifically, the component $\Phi_{i,j}$ is smooth and periodic Hamiltonian on the dense open subset given by the inverse image of $\left\{ \textbf{u} \in \Delta_\lambda ~\colon~ u_{i-1,j+1} < u_{i,j} < u_{i+1, j-1} \right\}$.
In particular, $\mathrm{Gr}(2,\C^n)$ contains an algebraic torus $(\C^\times)^{2(n-2)}$ consisting of torus orbits generated by $\Phi_{i,j}$'s.
Although each component may \emph{not} be smooth on the whole $\mcal{O}_\lambda$, their certain combinations can still be smooth. 
For instance, the combination $\Phi_{1, j+1} + \Phi_{2,j} - \Phi_{1,j} - \Phi_{1,j-1} -2$, the difference of sums of components in  consecutive anti-diagonals, is a moment map for $\theta_{j+1}$ in~\eqref{equ_Taction}. 

Since the action of the algebraic torus (in full dimension) does \emph{not} extend to $\mathrm{Gr}(2,\C^n)$ for $n \geq 4$, non-torus Lagrangian fibers of $\Phi_\textup{GC}$ appear at lower dimensional strata of $\Delta_{(2,n)}$. 
Such Lagrangian fibers are classified in \cite{CKO}, see also \cite{BMZ}. 
In the case of Grassmannian of 2-planes, it turns out that fibers are diffeomorphic to a product of 3-spheres or tori. 

\begin{figure}[h]
	\begin{center}
		\includegraphics[scale=0.8]{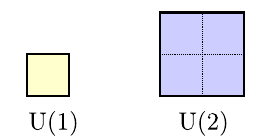}
	\caption{\label{fig_blocks} $\mathrm{U}(1)$-block and $\mathrm{U}(2)$-block}	
	\end{center}
\end{figure}

There is a simple combinatorial rule using two types of square blocks in the ladder diagram as shown in Figure~\ref{fig_blocks},
which tells us whether a fiber is a Lagrangian or not and reveals its diffeomorphic type.
\begin{itemize}
\item ($\mathrm{U}(1)$-block) is a $(1 \times 1)$-block, which amounts to $\mathbb{S}^1$,
\item ($\mathrm{U}(2)$-block) is a $(2 \times 2)$-block, which amounts to $\mathbb{S}^3 \times \mathbb{S}^1$.
\end{itemize}
\vspace{0.2cm}

For any given face $f$ of $\Delta_{(2,n)}$, consider the admissible diagram $\Gamma$ corresponding to $f$ in Proposition~\ref{prosposition_ACK}. 
The diagram $\Gamma$ divides the rectangular region $(n-2) \times 2$ into several pieces. 
The fiber over a point in the relative interior of $f$ is a Lagrangian if and only if any divided piece is either a $\mathrm{U}(1)$-block or a $\mathrm{U}(2)$-block and, in addition, is bounded by $\Gamma$. 
If this is the case, both $f$ and $\Gamma$ are said to be \emph{Lagrangian}.
Equivalently, a Lagrangian admissible diagram $\Gamma$ in a ladder diagram $\Gamma_{(2,n)}$ satisfies 
\begin{equation}\label{equ_numberofu1u2}
n_1 + 4 n_2 = 2(n-2)
\end{equation}
and vice versa, where $n_i$ is the number of $\mathrm{U}(i)$-blocks bounded by $\Gamma$ for $i = 1, 2$.
In summary, we have the following theorem. 

\begin{theorem}[\cite{CKO}]\label{thm_CKOaa}
Suppose that a face $f$ corresponds to $\Gamma$ via the correspondence~\eqref{equ_facediagramcorrespondence}. 
The admissible diagram $\Gamma$ is Lagrangian if and only if the fiber over each point in the relative interior of $f$ is a Lagrangian. 
Moreover, the fiber is diffeomorphic to $(\mathbb{S}^1)^{n_1 + n_2} \times (\mathbb{S}^3)^{n_2}$ where $n_i$ is the number of $\mathrm{U}(i)$-blocks bounded by $\Gamma$ for $i = 1, 2$. 
\end{theorem}

Furthermore, a Lagrangian fiber is \emph{monotone} if and only if its position satisfies the so-called Bohr-Sommerfeld condition as below.
\begin{theorem}[\cite{CK-mono}]\label{thm_CKmo}
Suppose that a face $f$ corresponds to $\Gamma$ is Lagrangian.
The Lagrangian fiber over a point $\textbf{\textup{u}}$  in the relative interior of $f$ is monotone if and only if 
\begin{equation}\label{equation_positionofmonotone}
\begin{cases}
u_{1,i} = u_{1,i+1} = u_{2,i} = u_{2,i+1} = i - 1  &\mbox{for $\square^{(1,i)}, \square^{(1,i+1)}, \square^{(2,i)}, \square^{(2,i+1)}$ in $\mathrm{U}(2)$-block,} \\
u_{i,j} = j - i  &\mbox{for $\square^{(i,j)}$ in $\mathrm{U}(1)$-block.}
\end{cases}
\end{equation}
\end{theorem}

\begin{example}\label{exa_fiberofgr24}
Let $\mcal{O}_\lambda \simeq \mathrm{Gr}(2,\C^4)$ for $\lambda = (2, 2, -2, -2)$. There are two Lagrangian faces of $\Delta_{(2,4)}$, which come from $\Gamma(g)$ and $\Gamma(h)$ in Figure~\ref{fig_ladder4}.
On the other hand, $\Gamma(f_1)$ is \emph{not} Lagrangian. Although the whole ($2 \times 2$)-block is cut into four $\mathrm{U}(1)$-blocks in this case, one of the four is not bounded by $\Gamma(f_1)$. 
Also, $\Gamma(f_2)$ are \emph{not} Lagrangian since one of the divided pieces is neither $\mathrm{U}(1)$-block nor $\mathrm{U}(2)$-block.
The fibers over the relative interiors of $g$ and $h$ are respectively $\mathrm{U}(2)$ and $\mathbb{T}^4$ by Theorem~\ref{thm_CKOaa}.
\end{example}

\begin{example}\label{exa_fiberpofgr26}
The Grassmannian $\mathrm{Gr}(2,\C^6) \simeq \mcal{O}_\lambda$ for $\lambda = (4,4,-2,-2,-2,-2)$ contains five Lagrangian faces, see Figure~\ref{fig_gr26Lagfaces}. 

\begin{figure}[h]
	\begin{center}
		\includegraphics[scale=0.8]{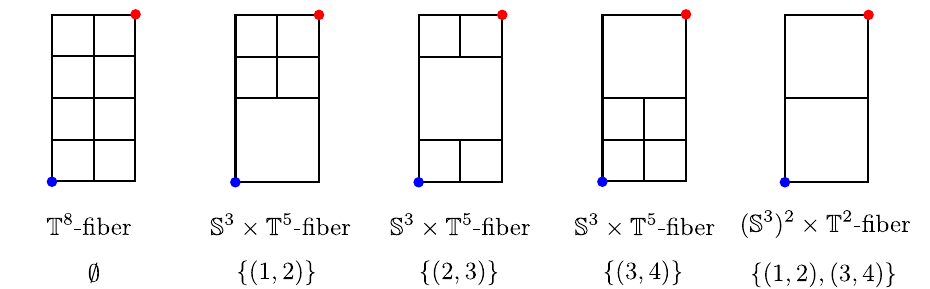}
	\caption{\label{fig_gr26Lagfaces} Lagrangian faces in $\mathrm{Gr}(2,\C^6) \simeq \mcal{O}_\lambda$ for $\lambda = (4,4,-2,-2,-2,-2)$.}	
	\end{center}
\end{figure}
\end{example}

For later use, we introduce combinatorial objects parametrizing Lagrangian faces. 
Let $\scr{P}_n$ be the set of pairs consisting of two consecutive integers in $\{1, 2, \cdots, n-2\}$, that is,
$$
\scr{P}_n := \left\{ (1,2), (2,3), \cdots, (n-3, n-2) \right\}.
$$
For a subset $\mcal{I} \subset \scr{P}_n$, the number of pairs in $\mcal{I} $ is denoted by $| \mcal{I}|$.
Let $\scr{I}_n$ be the collection of subsets of $\scr{P}_n$ whose distinct elements do \emph{not} share any common integers. 
Namely, for $\mcal{I} \in \scr{I}_n$ and for two pairs $(i, i+ 1)$ and $(j, j+1) \in \mcal{I}$, we have $i \neq j-1$ and $i \neq j + 1$. 
Equip $\scr{I}_n$ with a partial order given by
\begin{equation}\label{equ_partialorder}
\mcal{I}_1 < \mcal{I}_2 \, \Leftrightarrow \, (i, i+1) \in \mcal{I}_2 \mbox{ whenever } (i, i+1) \in \mcal{I}_1.
\end{equation}
Denote by $\scr{I}^\textup{max}_n$ the set of maximal elements of $\scr{I}_n$ with respect to the partial order~\eqref{equ_partialorder}. 
Note that $\scr{I}_n$ parametrizes monotone Lagrangians. 

\begin{remark}
We shall construct an immersed Lagrangian associated with each element in $\scr{I}_n$. 
We will see, however, that immersed Lagrangians associated with maximal elements in $\scr{I}^\textup{max}_n$ are sufficient for the purpose of recovering the Rietsch's mirror (see Lemma~\ref{lemma_enoughRIe}).
\end{remark}

\begin{example}
When $n = 6$, the index sets are 
$$
\begin{cases}
\scr{I}_6 = \{\emptyset, \{(1,2)\}, \{(2,3)\}, \{(3,4)\}, \{(1,2), (3,4)\} \} \\
\scr{I}^\textup{max}_6 = \{\{ (1,2), (3,4) \}, \{ (2,3) \}\}.
\end{cases}
$$
As in Figure~\ref{fig_gr26Lagfaces}, $\scr{I}_6$ parametrizes the set of faces in $\Delta_{(2,6)}$ whose relative interior has Lagrangian fibers. 
In other words, it parametrizes the monotone Lagrangian fibers of $\mathrm{Gr}(2,\C^6)$.

By applying the criteria in \cite[Section 9]{CKO}, the monotone non-toric Lagrangians corresponding to $ \{(1,2)\}$ or $\{(3,4)\}$ are displaceable. Thus they are trivial objects in the Fukaya category, and we do not need to take them into account.  
The other two monotone Lagrangians will be replaced by immersed Lagrangians to recover the Rietsch's mirror later on.
\end{example}

To compute the Floer potential (introduced by Fukaya-Oh-Ohta-Ono \cite{FOOO}) of a Lagrangian GC torus $\scr{T} \subset \mcal{O}_\lambda$, Nishinou-Nohara-Ueda \cite{NNU} constructed a degeneration of the GC system on $\mcal{X}_1$ into the toric moment map on $\mcal{X}_0$. 
The toric degeneration in~\eqref{equ_toricdege} is a family of Fano varieties and the central fiber $\mcal{X}_0$ admits a small resolution of the singular loci. Using these properties, they showed that a holomorphic disk has Maslov index strictly greater than two whenever it  intersects codimension $\geq 2$ strata. Therefore, the resulting superpotential is equal to the Givental-Hori-Vafa potential, analogously to the toric Fano case in Cho-Oh \cite{CO-T}.

\begin{theorem}[\cite{NNU}]\label{theorem_NNUmain}
The (Floer) potential function of $\scr{T}\subset \mathrm{Gr}(2,\C^n)$ is  
\begin{equation}\label{equ_potentialGC}
W_{\scr{T}}(\textbf{\textup{z}}) = z_{2,1} +  \frac{T^n}{z_{1,n-2}} + \sum_{j=1}^{n-3}
\left( \frac{z_{1,j+1}}{z_{1,j}} +  \frac{z_{2,j+1}}{z_{2,j}} \right) + \sum_{j=1}^{n-2} \frac{z_{1,j}}{z_{2,j}} 
\end{equation}
where $z_{i,j}$ is the exponential variable corresponding to the loop generated by $u_{i,j}$. 
\end{theorem}

\begin{remark}
Here, we only presented the potential function of a GC torus fiber in the case of $\mathrm{Gr}(2,\C^n)$.
It is also proved in \cite{NNU} that the Floer potential of a GC torus fiber agrees with the one introduced by Givental-Kim and Eguchi-Hori-Xiong in \cite{GK95, EHX97} for every partial flag varieties of Lie type $A$.
\end{remark}

\section{Floer theoretical strategy to construct SYZ mirrors for partial flag varieties}

In this section, we outline a program to construct open mirror symmetry for partial flag varieties.

\subsection{Strategy to construct mirrors for partial flag varieties}\label{sec_general}
Let $X$ be a partial flag variety.  The most distinguished torus fiber in the GC system of $X$ is the monotone one, which is denoted by $\scr{T}$ 
in $X$.   We can start with the local mirror $(\scr{U}_\scr{T}, W_\scr{T})$ of the torus $\scr{T}$, where $\scr{U}_\scr{T} \simeq (\Lambda_\textup{U})^{\dim_{\C} X}$ is its Maurer-Cartan deformation space given as in Section~\ref{subsec:fdconn}, and $W_\scr{T}$ is the Floer potential defined as the count of Maslov two holomorphic disks bounded by $\scr{T}$. 
The pair $(\scr{U}_\scr{T}, W_\scr{T})$ can be computed with the aid of a toric degeneration, see \cite{NNU, CO-T, FOOO-T}.  

We would like to construct the mirror $\check{X}$ of the complement of an anti-canonical divisor in $X$.  Then the Landau-Ginzburg mirror of $X$ would be given by $(\check{X},W)$, where $W$ is the disk potential counting holomorphic disks emanated from the anti-canonical divisor.  This is in line with the general program of Auroux \cite{auroux07}.  However, a big difference from the toric case is that, the inverse image of the boundary of the GC polytope is \emph{not} a divisor.  The reason is that the GC system has non-toric Lagrangian fibers over some of the faces.

In other words, the GC system does not restrict to a torus fibration on the anti-canonical divisor complement.  We need to take into account of monotone non-toric Lagrangian fibers in the GC system.  

In the cotangent bundle of each non-toric Lagrangian, viewed as a local neighborhood of the Lagrangian, we construct an immersed Lagrangian $\scr{L}$ whose certain Maurer-Cartan deformations are quasi-isomorphic to $\scr{T}$.  We use the quasi-isomorphisms (Definition \ref{def:qisomLag}) to derive the gluing data between their Maurer-Cartan deformation spaces (local mirrors). 
Using the compatibility between the quasi-isomorphism and Floer potentials, the disk potential of $\scr{L}$ can be explicitly computed from $W_\scr{T}$.
Moreover, the local mirrors glued together and give a partial compactification of the torus chart $\scr{U}_\scr{T}$.

For partial flag varieties of Lie type $A$, every GC fiber is diffeomorphic to the total space of an iterated bundle whose fibers are either points or products of odd-dimensional spheres, see \cite{CKO} for details.
We shall study cotangent bundles of odd dimensional spheres, and find suitable immersed Lagrangians in the local models.  
Notice that partial flag varieties are locally defined by quadratic equations of the form
$$
X_{1} Y_{1} + X_{2} Y_{2} + \cdots + X_{n} Y_{n} + 1 = 0,
$$
each of which exactly defines the cotangent bundle of an odd-dimensional sphere. 

In this article, we shall construct an immersed Lagrangian in the local chart $T^* \mathbb{S}^3$ and $T^* \mathrm{U}(2) \simeq T^* \mathbb{S}^3 \times T^* \mathbb{S}^1$ and derive the relationship between the constructed Lagrangians and the Lagrangian torus $\scr{T}$. 
The first two examples of flags using these two as local models are the complete flag variety $\mcal{F}\ell(1,2;3)$ and the Grassmannian $\mathrm{Gr}(2,\C^4)$ respectively. Indeed, they are compactifications of a single $T^* \mathbb{S}^3$ in $\CP^2 \times \CP^2$ and a single $T^* \mathrm{U}(2)$ in $\C\bP^5$ respectively. 

Indeed the LG mirror constructed from the torus $\scr{T}$ is good enough to study mirror symmetry of complete flag varieties, since it already contains enough critical points.
Thus we will first focus on $\mathrm{Gr}(2,\C^4)$ for which the Jacobian ring of the mirror from $\scr{T}$ is not isomorphic to the quantum cohomology of $\mathrm{Gr}(2,\C^4)$. 
We will discuss how to construct its partially compactified mirror in details in Section~\ref{sec_LagfibonF3}.
Then we will move on to the Grassmannian $\mathrm{Gr}(2,\C^n)$, which has charts being products of $T^* \mathrm{U}(2)$ and $T^* \mathbb{S}^1$.  Employing local models as building blocks, 
we will complete a mirror construction for $\mathrm{Gr}(2,\C^n)$ in Section~\ref{section_mirrorcongr2n}.

Following the scheme of \cite{CHLabc,CHLtoric,CHLnc, CHLgl}, the mirror construction can be upgraded into the categorical level. 
Namely, each reference Lagrangian has an associated localized mirror functor.
By gluing the mirror functors through the quasi-isomorphisms, we obtain a glued mirror functor from the Fukaya category to the category of Matrix factorizations. We exhibit homological mirror symmetry in some examples in Section~\ref{sec:HMS} although it does not crucially use the glued functor due to isolatedness of the critical loci of the mirror potential. 

Our program does \emph{not} require constructing a  special Lagrangian fibration, which is one of the main difficulties in the SYZ program.  Instead we construct a finite collection of Lagrangian immersions of $X$ in place of singular SYZ fibers.  This is good enough for mirror symmetry of anti-canonical divisor complement and also the Fano variety of $X$ itself.

In a forthcoming work with Cheol-Hyun Cho, the authors will study other local models and apply them to construct mirrors and derive homological mirror symmetry for more general classes of partial flag varieties. 

\subsection{Immersed charts for Rietsch's mirrors}\label{subimmersedre}

Inspired by Peterson's lecture, Rietsch \cite{Rie} constructed Lie-theoretical LG models of general partial flag manifolds. Later on, Marsh and Rietsch \cite{MR} described the Lie-theoretical mirror of  $ \mathrm{Gr}(k, \C^n)$ in terms of the dual Pl{\"u}cker coordinates and investigated the cluster structure on it. 
It turns out that a Floer theoretically derived SYZ mirror matches up with the Lie-theoretical mirror.
After recalling their work briefly on $ \mathrm{Gr}(2, \C^n)$, we introduce immersed charts for the Rietsch's mirror. 
Those charts that will be geometrically derived will be used to fill up the strata that cannot be covered by the cluster charts.

We write $X:= \mathrm{Gr}(2, \C^n)$ in this section. Regarding $X$ as $G/P$ (see Section~\ref{subsec:vardefgr}), the Rietsch's mirror consists of the complement of an anti-canonical divisor in its Langlands dual Grassmannian $\check{G}/\check{P} = \mathrm{Gr}(n-2, \C^n)$ together with a superpotential.
The map
\begin{equation}\label{equ_pijpihatjhat}
p_{\hat{i}, \hat{j}} := p_{1, \cdots, i-1, i+1, \cdots, j-1, j+1, \cdots, n} \mapsto p_{\vphantom{\hat{i}}i, \vphantom{\hat{j}}j}
\end{equation}
provides an identification between $\mathrm{Gr}(n-2, \C^n)$ with $\mathrm{Gr}(2, \C^n)$. 

There is an one-to-one correspondence between the Pl{\"u}cker variables and the Young diagrams in the $\left(2 \times (n-2)\right)$-rectangle (or block). 
Consider a \emph{negative} path, a path from the top-right vertex to the bottom-left vertex with the shortest distance, lying in the $\left(2 \times (n-2)\right)$-block. Its \emph{horizontal} steps determine the indices of Pl{\"u}cker variable. 
On the other hand, cutting the $\left(2 \times (n-2)\right)$-block along the given negative path, the left-upper corner one will be taken as the corresponding Young diagram. 

Under the reflection with respect to the line through the right-top vertex with the slope $-1$, the Young diagrams in the $\left(2 \times (n-2)\right)$-rectangle map into those in $\left((n-2) \times 2\right)$-rectangle.
The identification~\eqref{equ_pijpihatjhat} is compatible with the reflection on Young diagrams when recording the horizontal steps of a reflected positive path in the $\left((n-2) \times 2\right)$-rectangle as the index of a variable.

Under the identification~\eqref{equ_pijpihatjhat}, the Rietsch's mirror of $\mathrm{Gr}(2, \C^n)$ is the LG model $(\check{X} ,W_{\textup{Rie}})$ defined as 
\begin{itemize}
\item $\check{X} := \mathrm{Gr}(2, \C^n) \backslash \scr{D}$ where $\scr{D} = \left\{ p_{{1},{2}} \cdot p_{2,3} \cdots p_{n-1,n} \cdot p_{1,n} = 0 \right\}$,
\begin{equation}\label{equ_Rietschsmir}
W_{\textup{Rie}}\left( \left[ p_{i,j} \right] \right) := q \frac{p_{2,n}}{p_{1,2}} + 
\sum_{j=2}^{n-1} \frac{p_{j-1,j+1}}{p_{j,j+1}} + \frac{p_{1,n-1}}{p_{1,n}}.
\end{equation}
\end{itemize}

In \cite{MR}, the authors found a cluster structure, which was invented by Fomin-Zelevinsky \cite{FZ02}, on the Rietsch's mirror $\left(\check{X}, W_{\textup{Rie}}\right)$ in~\eqref{equ_Rietschsmir}.
The superpotential~\eqref{equ_Rietschsmir}  has the variables $p_{j,j+1}$ $(1 \leq j \leq n)$ in the denominators, which are called \emph{frozen} variables. 
The other Pl{\"u}cker variables are called \emph{cluster} variables.
Observe that the variables in $\scr{D}$ are frozen variables, and hence $W_{\textup{Rie}}$ is a regular function on $\check{X}$. 

In fact, one can find that the variables in each monomial of the superpotential~\eqref{equ_Rietschsmir} are arranged by the following combinatorial rule described in terms of the associated Young diagrams.
For $j \neq 1$, the Young diagram associated with the variable in the numerator is made out of that in its denominator by adding one unit box. 
When $j = 1$, the numerator corresponds to the $\left((n-3) \times 1\right)$-block obtained by removing the bottom row and the second column from $\left((n-2) \times 2\right)$-block. We refer readers to \cite[Section 6]{MR} for details.

\begin{remark}
By the quantum Pieri rule, a variable in any numerator of superpotential~\eqref{equ_Rietschsmir} corresponds to the Schubert cycle obtained by multiplying the cycle corresponding to the unit box to the cycle corresponding to the frozen variable in its denominator.
\end{remark}

\begin{remark}
Nohara and Ueda constructed generalized GC systems in \cite{NU14} and obtained LG mirrors by counting holomorphic disks with boundaries on their torus fibers. In \cite{NUclu} these LG mirrors respect the cluster structure, and the resulting glued mirror partially covers $\left(\check{X}, W_{\textup{Rie}}\right)$. 
\end{remark}

\begin{example}\label{exam_rietschgr24}
For $\lambda = (2, 2, -2, -2)$, we have
$\mcal{O}_\lambda \simeq \mathrm{Gr}(2,\C^4)$. The Rietsch's mirror consists of $\check{X} := \mathrm{Gr}(2, \C^4) \backslash \{p_{1,2} \cdot p_{2,3} \cdot p_{3,4} \cdot p_{1,4} = 0\}$ and
\begin{equation}\label{equ_rietschpogr24}
W_{\textup{Rie}} \left( \left[ p_{i,j} \right] \right) = q \, \frac{p_{2,4}}{p_{1,2}} + \frac{p_{1,3}}{p_{2,3}} + \frac{p_{2,4}}{p_{3,4}} + \frac{p_{1,3}}{p_{1,4}} \colon \check{X} \to \C.
\end{equation}
\end{example}

We now introduce \emph{immersed charts} of $\left(\check{X}, W_{\textup{Rie}}\right)$, which will later turn out to be Maurer-Cartan deformation spaces of immersed Lagrangians (see Section \ref{section_mirrorcongr2n}). 
While the cluster charts can recover the Rietsch's mirror up to codimension $\geq 2$ only, we will see that the immersed charts fill out the missing parts lying outside any of cluster charts. 

\begin{definition}
For $\mcal{I} \in \scr{I}_n$, we define an \emph{immersed chart} associated with $\mcal{I}$ of $\check{X}$ by
\begin{equation}\label{equ_immersedchar}
\mcal{U}_\mcal{I} = \left\{ [p_{i,j}] \in \check{{X}} ~\colon~ p_{n-i-1,n} \neq 0 \mbox{ for $i$ with $(i, i+1) \notin \mcal{I}$} \right\}.
\end{equation}
\end{definition}

There are several different ways of describing this cluster structure.
For example, it is governed by the Postnikov diagram or the plabic graph as explained in \cite{MR}. 
For $\mathrm{Gr}(2,\C^n)$, the triangulations of an $n$-gon also provides a cluster structure.
 
\begin{figure}[h]
	\begin{center}
		\includegraphics[scale=0.7]{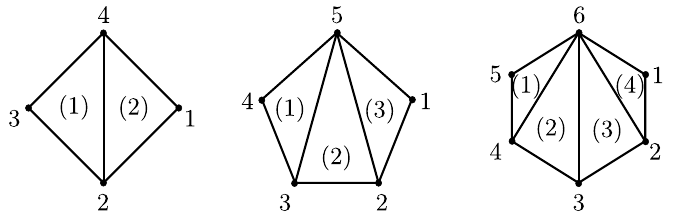}
	\caption{\label{fig_referenceclus} The reference chart $\mcal{U}_\emptyset$ for $\mathrm{Gr}(2,\C^n)$ for $n = 4,5,6$.}	
	\end{center}
\end{figure}

We start from the $n$-gon and enumerate the vertices of the $n$-gon respecting the clockwise orientation. 
Notice that $\mcal{U}_\emptyset$ is given by $\{p_{2,n} \cdots p_{n-2, n} \neq 0\}$. 
We call $\mcal{U}_\emptyset$ the \emph{reference chart} of $\check{X}$.
We then associate the chart with the triangulation of the $n$-gon by drawing the edges connecting the $n$-th vertex and the other vertices.
This particular triangulation divides the $n$-gon into $(n-2)$-triangles. We assign the integers $(i)$ to the triangle whose vertices are $n, n-i$, and $n-i-1$, see Figure~\ref{fig_referenceclus}. 

Each combinatorial input $\mcal{I} \in \scr{I}_n$ indeed tells us how the triangulation for the reference chart $\mcal{U}_\emptyset$ mutates into another triangulation. For each pair $(i, i+1) \in \mcal{I}$, consider the rectangle whose vertices are $n, n-i, n-i-1$, and $n-i-2$ and perform the flip move, changing the edge $(n,n-i-1)$ to the edge $(n-i,n-i-2)$. 
Immersed chart arises through the mutation process. 
We assign the immersed chart $\mcal{U}_\mcal{I}$ to the subdivision of $n$-gon by removing the edges $(n, n-i-1)$ for all $(i,i+1) \in \mcal{I}$. 
Those subdivision attached to the immersed charts will be used to produce the identification between the immersed charts and the Maurer-Cartan deformation spaces of immersed Lagrangians, see Section~\ref{completesyzmirrgr2n}. 

\begin{example}
As in Figure~\ref{fig_coord25}, there are five algebraic torus charts associated with five triangulations. 
For a choice $\mcal{I} = \{ (1,2)\}$, consider two triangles $(5,4,3)$ and $(5,3,2)$ and then obtain triangles $(5,4,2)$ and $(4,3,2)$ after the flip move corresponding to $\mcal{I}$. 
The immersed chart $\mcal{U}_{(1,2)}$ has the partition of the $5$-gon into one $4$-gon $(5,4,3,2)$ and $3$-gon $(1,2,5)$. 

\vspace{-0.2cm}
\begin{figure}[h]
	\begin{center}
		\includegraphics[scale=0.6]{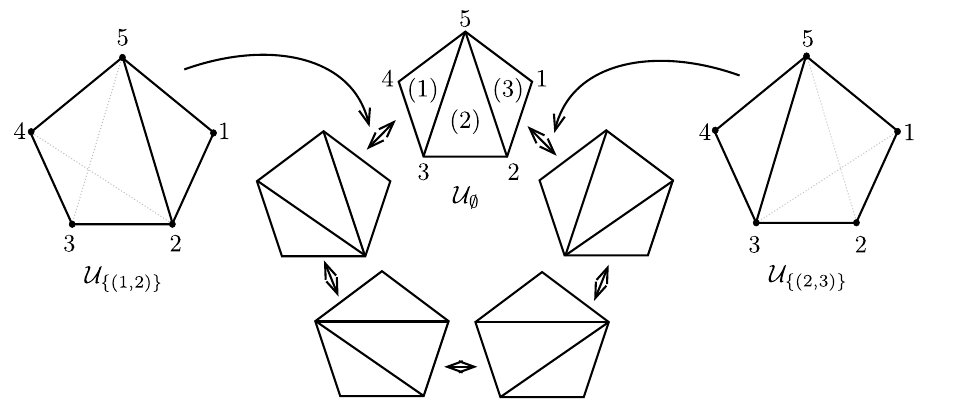}
	\caption{\label{fig_coord25} Immersed charts and cluster chats of $\check{X}$ for $\mathrm{Gr}(2,\C^5)$.}	
	\end{center}
\end{figure}
\vspace{-0.2cm}
\end{example}

The following lemma characterizes immersed charts that cover the Rietsch's mirror of $\mathrm{Gr}(2,\C^n)$.

\begin{lemma}\label{lemma_enoughRIe}
Let $(\check{X}, W_{\textup{Rie}})$ be the Rietsch's mirror of $X = \mathrm{Gr}(2,\C^n)$. The mirror space is covered by the set of complex charts $\{ \mcal{U}_\mcal{I} ~\colon~ \mcal{I} \in \scr{I}_n^\textup{max} \}$ where $\scr{I}_n^\textup{max}$ is given as in Section~\ref{subsection_blockcombi}. Namely, $\check{{X}} = \mathrm{Gr}(2, \C^{n}) \backslash \scr{D}$
is contained in
$$
\bigcup_{\mcal{I} \in \scr{I}^\textup{max}_n} \left( \check{{X}} \cap \mcal{U}_\mcal{I} \right).
$$
\end{lemma}

\begin{proof}
Assume that both $p_{i,n}$ and $p_{i+1,n}$ vanish.
Due to the Pl{\"u}cker relation 
\begin{equation*}
p_{i,i+1} \cdot p_{i+2,n} - p_{i, i+2}  \cdot p_{i+1, n} + p_{i, n}  \cdot p_{i+1, i+2} = 0,
\end{equation*}
$p_{i,n} = p_{i+1,n} = 0$ forces $p_{i+2, n} = 0$. 
Proceeding inductively, we get $p_{n-1,n} =0$. 
Thus, any point satisfying $p_{i, n} = p_{i+1, n} = 0$ lies in $\scr{D}$. 

For any subset $\mcal{J}$ of $\scr{P}_n$, set $\mcal{U}_\mcal{J} := \left\{ [p_{i,n}] \in \check{X} ~\colon~ p_{n-i-1,n} \neq 0 \mbox{ for $i$ with $(i,i+1) \notin \mcal{J}$} \right\}.$ 
If $\mcal{J}$  contains both $(i,i+1)$ and $(i+1,i+2)$, the above observation yields
\begin{equation}\label{equ_consecu}
\mcal{U}_\mcal{J} \subset \mcal{U}_{\mcal{J} \backslash \{(i,i+1)\}} \cup \mcal{U}_{\mcal{J} \backslash \{(i+1,i+2)\}}.
\end{equation}
Also, observe that
\begin{equation}\label{equ_incluss}
\mcal{I} \subset \mcal{I}^\prime \, \Rightarrow \, \mcal{U}_{\mcal{I}^\prime} \subset \mcal{U}_{\mcal{I}}.
\end{equation}

Combining~\eqref{equ_consecu} and~\eqref{equ_incluss}, we establish
$$
\check{X} = \bigcup_{\mcal{I} \subset \scr{P}_n} (\check{X}  \cap \mcal{U}_{{\mcal{I}}}) = \bigcup_{\mcal{I} \in \scr{I}_n} (\check{X}  \cap \mcal{U}_{{\mcal{I}}}) = \bigcup_{\mcal{I} \in \scr{I}^\textup{max}_n} (\check{X}  \cap \mcal{U}_{{\mcal{I}}}).
$$
This completes the proof. 
\end{proof}

We next discuss how to adorn the superpotential on the immersed chart $\mcal{U}_\mcal{I}$. 
The reference chart $\mcal{U}_\emptyset$ consists of the frozen variables and the set $\{ p_{2,n}, p_{3,n}, \cdots, p_{n-2,n} \}$ of cluster variables. 
By using the Pl{\"u}cker relation~\eqref{equ_pluckker}, 
$W_{\textup{Rie}}$ restricted on $\mcal{U}_\emptyset$ can be written in terms of the above set of variables. 
It agrees with $W_\scr{T}$ in~\eqref{equ_potentialGC} up to a suitable coordinate change. 

As observed in \cite[Proposition 6.10]{MR}, the LG model $(\mcal{U}_\emptyset, W_{\textup{Rie}} |_{\mcal{U}_\emptyset})$ can be identified with the model in \cite{EHX97, BCKV, NNU}. 
In other words, the LG mirror $\left((\C^\times)^{n(n-2)}, W_{\scr{T}}\right)$ in~\eqref{equ_potentialGC} provides one cluster chart of $\left(\check{X}, W_{\textup{Rie}}\right)$. More specifically, under the coordinate change 
\begin{equation}\label{equ_coordinateidentifi}
z_{1,j} \mapsto \frac{p_{n-j-1,n-j}}{p_{n-j,n}}, \quad z_{2,j} \mapsto \frac{p_{n-j-1,n}}{p_{n-1,n}},
\end{equation} 
the algebraic torus $(\C^\times)^{n(n-2)}$ embeds into $\check{X}$, and $W_{\textup{Rie}}$ restricted to this embedded torus matches  $W_{\scr{T}}$. 

Consider the immersed chart $\mcal{U}_\mcal{I}$ for $\mcal{I} \in \scr{I}_n$. 
Suppose that $ \mcal{I} = \{(i_1, i_1+1), \cdots, (i_r, i_r +1) \}$.  
Then the superpotential on $\mcal{U}_\mcal{I}$ is uniquely determined by clearing the variables $p_{n-1-i_s,n}$ $(s=1, \cdots, r)$ in the denominator of each Laurent monomial (containing $p_{n-1-i_s,n}$) in $ W_{\textup{Rie}} |_{\mcal{U}_\emptyset}$
by the relation
\begin{equation}\label{plui2i21is2}
(i_s, i_s+1, i_s+2, n) \colon p_{i_s, i_s+1} \cdot p_{i_s+2, n} - p_{i_s,i_s +2} \cdot p_{i_s+1, n} + p_{i_s,n} \cdot p_{i_s +1, i_s +2} = 0. 
\end{equation}
Once we do not have any $p_{n-1-i_s,n}$ in any denominators of monomials, the expression extends to $\mcal{U}_\mcal{I}$. 

\begin{remark}
To write down cluster transformations and cluster charts explicitly, one needs to assign the dual Pl{\"u}cker variable $p_{i,j}$ to the edge $(i,j)$ of a triangulation of the $n$-gon. 
By connecting the barycenters of the edges that forms each triangle and giving an orientation, the quiver can be associated to a triangulation. 
Then quiver mutations give rise to cluster transformations, see \cite{FWZ1} for instance. 
We will not recall the process because we are mainly concerned with the specific types of immersed charts in Lemma~\ref{lemma_enoughRIe}. 
\end{remark}

\begin{example}
Consider $\mcal{O}_\lambda \simeq \mathrm{Gr}(2,\C^6)$ where $\lambda = (4,4,-2,-2,-2,-2)$. 
In this case, the superpotential function associated with $\emptyset$ is given by 
\begin{align*}
W_\emptyset ([p_{i,j}]) &= q \, \frac{p_{26}}{p_{12}} + \frac{p_{12} p_{36}}{p_{26} p_{23}} + \frac{p_{16}}{p_{26}} + \frac{p_{23}p_{46}}{p_{36}p_{34}} + \frac{p_{26}}{p_{36}} + \frac{p_{34} p_{56}}{p_{46} p_{45}} + \frac{p_{36}}{p_{46}} \\
&+ \frac{p_{12} p_{56}}{p_{26} p_{16}} + \frac{p_{23} p_{56}}{p_{36} p_{26}} + \frac{p_{34} p_{56}}{p_{46} p_{36}} +  \frac{p_{45}}{p_{46}} + \frac{p_{46}}{p_{56}},
\end{align*}
which is the potential function of $\scr{T}$ (after setting $q = T^6$). 
As examples, consider $\mcal{I} := \{(2,3)\}$ and $\mcal{I}^\prime := \{(1,2), (3,4)\}$ (see Figure~\ref{fig_gr26Lagfaces}). 
The superpotential functions restricted on the immersed charts $\mcal{U}_\mcal{I}$ and $\mcal{U}_{\mcal{I}^\prime}$ are respectively
\begin{align}\label{align_gri2}
W_{\mcal{I}} = q \, \frac{p_{26}}{p_{12}} + \frac{p_{12} p_{36}}{p_{26} p_{23}} + \frac{p_{16}}{p_{26}} + \frac{p_{24}}{p_{34}} + \frac{p_{34} p_{56}}{p_{46} p_{45}} + \frac{p_{36}}{p_{46}} + \frac{p_{12} p_{56}}{p_{26} p_{16}} + \frac{p_{24} p_{56}}{p_{26} p_{46}} + \frac{p_{45}}{p_{46}} + \frac{p_{46}}{p_{56}}.
\end{align}
and 
\begin{align}\label{align_gri}
W_{\mcal{I}^\prime} = q \, \frac{p_{26}}{p_{12}} + \frac{p_{13}}{p_{23}} + \frac{p_{23}p_{46}}{p_{36}p_{34}} + \frac{p_{26}}{p_{36}} + \frac{p_{35}}{p_{45}} + \frac{p_{13}p_{56}}{p_{16}p_{36}} + \frac{p_{35}}{p_{36}} + \frac{p_{46}}{p_{56}}. 
\end{align}
Two LG models $(\mcal{U}_\mcal{I}, W_{\mcal{I}})$ and $(\mcal{U}_{\mcal{I}^\prime}, W_{\mcal{I}^\prime})$ cover $\left(\check{X}, W_{\textup{Rie}}\right)$. In Section~\ref{section_mirrorcongr2n}, they will be realized as Maurer-Cartan deformation spaces of immersed Lagrangians.
\end{example}

\begin{remark}
In order to deal with a convergence issue of $A_\infty$-operations, we need to work over \emph{Novikov rings} in Lagrangian Floer theory.
Precisely speaking, our first goal is to construct the Rietsch's mirror over the Novikov rings using Floer theory. 
To match up two mirrors, we let the formal variable $T^\bullet$ in~\eqref{theorem_NNUmain} play a role of the K{\"a}hler parameter $q$ in~\eqref{equ_rietschpogr24}.  
\end{remark}

\section{Complete SYZ mirror of $\mathrm{Gr}(2,\C^4)$}\label{sec_LagfibonF3}

The aim of this section is to construct an immersed Lagrangian in the model $T^* \mathbb{S}^3$ and $T^* \mathrm{U}(2)$ as building blocks for a completion of mirrors of partial flag varieties. 
We analyze the quasi-isomorphisms between the immersed Lagrangian and Lagrangian tori therein.
As an example of the completion process, the Grassmannian $\mathrm{Gr}(2,\C^4)$ will be studied in details.
The glued mirror via quasi-isomorphisms recovers the Rietsch's mirror.

\subsection{Construction of immersed Lagrangian in $T^* \mathbb{S}^3$ and $T^* \mathrm{U}(2)$}\label{sec_consimmergr24}

As in~\eqref{equ_pluckker}, the Grassmannian $\mathrm{Gr}(2, \C^4)$ embeds into $\CP^5$ as a hypersurface defined by 
\begin{equation}\label{equ_pluckerrel}
Z_{12} Z_{34} - Z_{13} Z_{24} + Z_{14} Z_{23} = 0.
\end{equation}
The toric degeneration in~\eqref{equ_toricdege} is then a family of hypersurfaces 
\begin{equation}\label{equ_pluckertori}
\mcal{X}_\varepsilon := \{ [Z_{ij}] \in \CP^5 ~\colon~ \varepsilon \cdot Z_{12} Z_{34} - Z_{13} Z_{24} + Z_{14} Z_{23} = 0 \}.
\end{equation}

\begin{figure}[h]
	\begin{center}
		\includegraphics[scale=0.25]{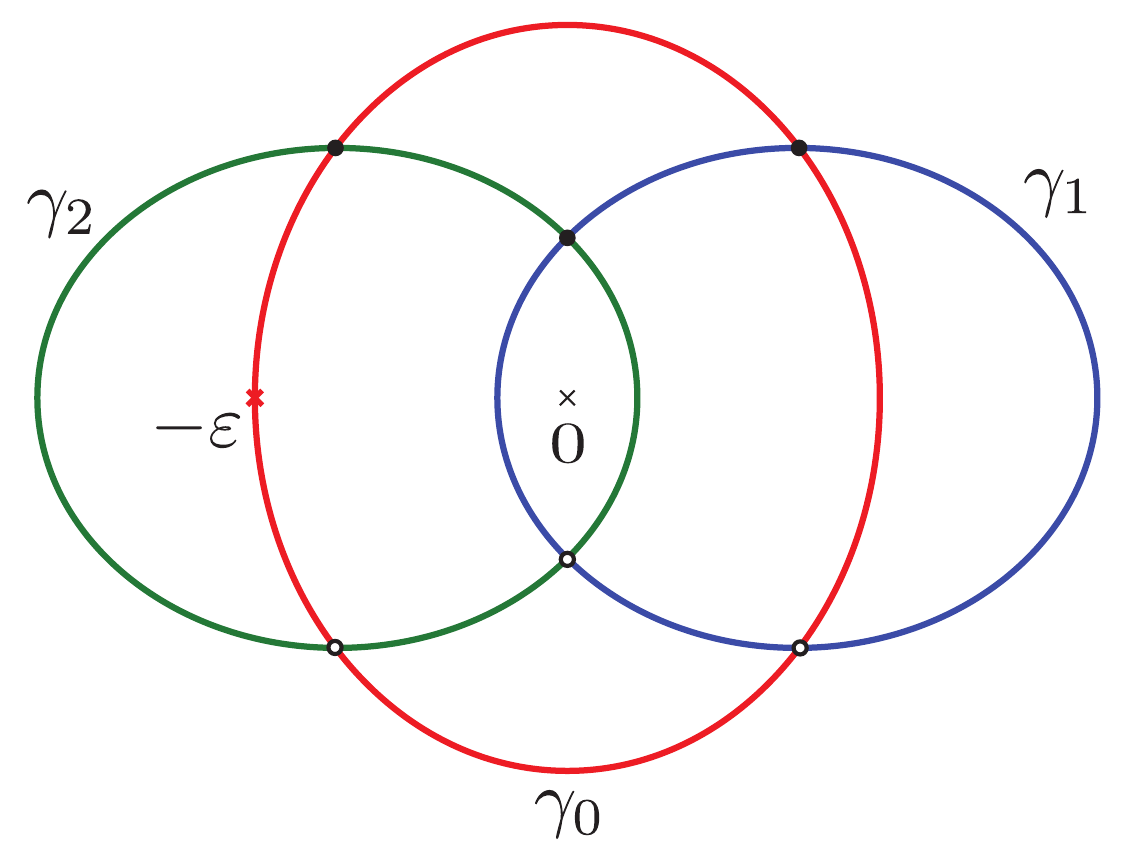}
		\caption{Simple closed curves $\gamma_0, \gamma_1$, and $\gamma_2$ for $\scr{L}_0, \scr{L}_1,$ and $\scr{L}_2$ on the complex plane $\C_Z$.}
		\label{Fig_baseforGr24}
	\end{center}
\end{figure}

By passing it to the affine chart given by $\{ Z_{12} \neq 0\}$,~\eqref{equ_pluckertori} can be written as
\begin{equation}\label{equ_affinecharteqq}
\varepsilon \cdot \frac{Z_{34}}{Z_{12}} - \frac{Z_{13}}{Z_{12}} \cdot \frac{Z_{24}}{Z_{12}} + \frac{Z_{14}}{Z_{12}} \cdot \frac{Z_{23}}{Z_{12}} = 0. 
\end{equation}
We then take a symplectic reduction of the open set $Z_{12} \cdot Z_{34} \neq 0$ of $\mathrm{Gr}(2,\C^4)$ by the $\mathbb{S}^1$-action $\theta_4$ in~\eqref{equ_Taction}. 
The reduced space is identified with a smoothing of a conifold, which is given by
\begin{equation}\label{equ_localconifoldsmo}
\varepsilon - X_{1} \cdot Y_{1} + X_{2} \cdot Y_{2} = 0
\end{equation}
where $X_{1} = {Z_{13}}/{Z_{12}}, Y_{1} = {Z_{24}}/{Z_{34}}, X_{2} = {Z_{14}}/{Z_{34}}$, and $Y_{2} = {Z_{23}}/{Z_{12}}$. 
The equation~\eqref{equ_localconifoldsmo} is diffeomorphic to the cotangent bundle of $\mathbb{S}^3$. 
Regarding~\eqref{equ_localconifoldsmo} as the intersection
$$
\begin{cases}
Z = X_1 \cdot Y_1 + \varepsilon, \\
Z = X_2 \cdot Y_2,
\end{cases}
$$
the $Z$-projection defines a double conic fibration. 
It carries the Hamiltonian $\mathbb{T}^2$-action by
$$
\left( \left(\eta_1, \eta_2\right), (X_1, Y_1, X_2, Y_2) \right) \mapsto \left( e^{- \sqrt{-1} \eta_1} X_1, e^{\sqrt{-1} \eta_1} Y_1,   e^{- \sqrt{-1} \eta_2} X_2, e^{\sqrt{-1} \eta_2} Y_2 \right).
$$
For $i = 0, 1$, and $2$, by taking $\mathbb{T}^2$-orbits over $\gamma_i$ as depicted in Figure~\ref{Fig_baseforGr24}, we produce two Lagrangian tori and one immersed Lagrangian. 
The immersed Lagrangian can be employed to (partially) compactify a mirror of a Lagrangian torus in $T^*\mathbb{S}^3$ that degenerates into a toric fiber of the conifold. 

\vspace{-0.2cm}
\vspace{-0.2cm}

\begin{remark}
By following the process, (after the reduction) one obtains a Lagrangian fibration on a smoothing of the conifold. 
It is an example of a Lagrangian fibration constructed by Gross \cite{Gross-eg} via Minkowski decomposition. 

There are many works related to this fibration. Chan-Pomerleano-Ueda \cite{CPU} studied mirror symmetry of the conifold using this fibration.  The third named author studied SYZ via Minkowski decompositions \cite{L14} and mirror symmetry of an orbifolded conifold jointly with Kanazawa \cite{KL}.
The first named and the third named authors with Fan and Yau in \cite{FHLY} studied the mirror geometry of the Atiyah flop with help of this model.
The second named author with Cho and Oh \cite{CKO} studied non-displaceable Lagrangians in this model. 
\end{remark}

In order to take the reduced $\C^\times$-factor into consideration, we define the map
\begin{equation}
\mcal{X}_\varepsilon \backslash \{ Z_{12} \cdot Z_{34} \neq 0\} \to \left\{ (Z_0, X_1, Y_1, X_2, Y_2) \in (\C^\times) \times \C^4 ~\colon~ \varepsilon - X_{1} \cdot Y_{1} + X_{2} \cdot Y_{2} = 0
 \right\}
\end{equation}
given by $[Z_{ij}] \mapsto \left({Z_{34}}/{Z_{12}}, {Z_{13}}/{Z_{12}}, {Z_{24}}/{Z_{34}}, {Z_{14}}/{Z_{34}}, {Z_{23}}/{Z_{12}} \right)$. 
It is diffeomorphic to $T^*\mathrm{U}(2)$.
The projection $Z = X_2 Y_2$ carries the fiber-wise Hamiltonian $\mathbb{T}^3$-action generated by $\theta_1 \circ \theta_4, \theta_2 \circ \theta_4$, and $\theta_4$ where $\theta_\bullet$'s are in~\eqref{equ_Taction}. Take the Hamiltonian
\begin{equation}\label{equ_hamiltonianfucntiont3}
(\Phi_{1,1} - \Phi_{2,2} , \Phi_{1,2} + \Phi_{2,1} - \Phi_{1,1} - \Phi_{2,2}, - \Phi_{2,2})
\end{equation}
The projection is a symplectic fibration.
For $i = 0, 1$, and $2$, by parallel transporting $\mathbb{T}^3$-orbits over $\gamma_i$ in Figure~\ref{Fig_baseforGr24}, we obtain two Lagrangian tori and one immersed Lagrangian by collecting orbits. 

We require our curves $\gamma_i$ to satisfy the condition that holomorphic sections over the regions bounded by $\gamma_0, \gamma_1$, and $\gamma_2$ have area $1$. 
In addition to the condition, we collect the following specific orbits$\colon$
\begin{enumerate}
\item The orbits of the $\mathbb{S}^1$-action $\theta_1 \circ \theta_4$ are at the level of $\mu_{\theta_1} + \mu_{\theta_4} = \Phi_{1,1} - \Phi_{2,2} = 0$,
\item The orbits of the $\mathbb{S}^1$-action $\theta_2 \circ \theta_4$ are at the level of $\mu_{\theta_2} + \mu_{\theta_4} = (\Phi_{1,2} + \Phi_{2,1}) - \Phi_{1,1} - \Phi_{2,2} = 0$,
\item The orbits of the $\mathbb{S}^1$-action $\theta_4$ are at the level of $\mu_{\theta_4} = - \Phi_{2,2} = 0$.
\end{enumerate}
As in Section \ref{sec_wallcrossingimmersedLag}, those conditions enable us to construct quasi-isomorphisms among constructed Lagrangians later on. 

Let us denote by $\scr{L}_i$ the Lagrangian consisting of such orbits over $\gamma_i$.
Then $\scr{L}_1$ and $\scr{L}_2$ are monotone Lagrangian tori because $c_1(T \mathrm{Gr}(2, \C^4)) = [\omega]$, which follows from our choice of $\lambda$ in~\eqref{equ_choiceofeigen}. Observe that $\scr{L}_0$ is an immersed Lagrangian intersecting with the non-toric monotone Lagrangian $\mathbb{S}^3 \times \mathbb{S}^1$ over $[-\varepsilon, 0]$ and the monotone Lagrangians $\scr{L}_1$ and $\scr{L}_2$ as well. The immersed Lagrangian $\scr{L}_0$ will play a crucial role in constructing the complete mirror in Section~\ref{sec_chartsmirrorGr(2c4)}.

\begin{figure}[h]
	\begin{center}
		\includegraphics[scale=0.85]{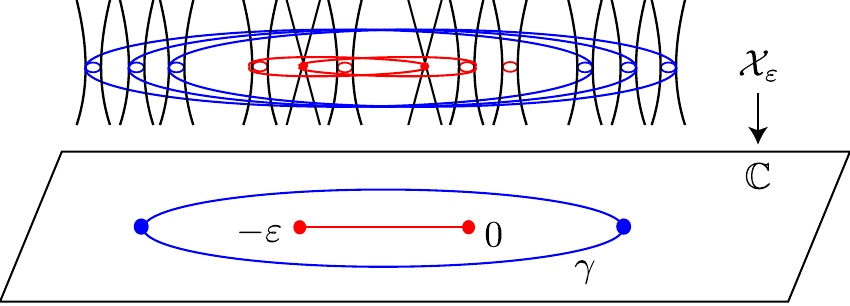}
	\caption{\label{Fig_tripleconic} Construction of Lagrangians $\scr{L}_\gamma$}	
	\end{center}
\end{figure}

\begin{remark}
There is a one-parameter family of Lagrangian submanifolds diffeomorphic to $\mathbb{S}^3 \times \mathbb{S}^1$ over the line segment $[-\varepsilon, 0]$, as the union of two vanishing cycles produces $\mathbb{S}^3$, and a choice of $\mathbb{S}^1$-orbit generated by $\theta_4$ provides the $\mathbb{S}^1$-factor, see Figure~\ref{Fig_tripleconic}. The family can be regarded as the counterpart of the GC fibers $\simeq \mathbb{S}^3 \times \mathbb{S}^1$ over the edge $g$ in Example~\ref{example_gr24diag}.
\end{remark}

\begin{remark}\label{remark_redcons}
We take the symplectic reduction of $\mathrm{Gr}(2,\C^4) \backslash \{ Z_{12} \cdot Z_{23} \cdot Z_{34} \cdot Z_{14} = 0 \}$ with respect to the $\mathbb{T}^3$-action generated by~\eqref{equ_hamiltonianfucntiont3}.
Regarding the $\C$-plane in Figure~\ref{Fig_tripleconic} as the one-dimensional reduced space, the Lagrangian $\scr{L}_i$ can be alternatively constructed by collecting the $\mathbb{T}^3$-orbit over a Lagrangian $\gamma_i$ in the reduced space. As being Lagrangian is a closed condition, $\scr{L}_0$ is an immersed Lagrangian even though $\gamma_0$ passes through $- \varepsilon$ where the $\mathbb{S}^1$-orbit of $\theta_2 \circ \theta_4$ degenerates at the zero level.
\end{remark}

\subsection{Local mirrors of the Lagrangians $\scr{L}_0, \scr{L}_1$, and $\scr{L}_2$}\label{localmirrorsecti}

To equip our Lagrangians with a grading, we take the meromorphic volume form
\begin{equation}\label{equ_Omegazgrad}
\Omega :=  d \log \left( \frac{Z_{23}}{Z_{12}} \right) \wedge d \log \left( \frac{Z_{34}}{Z_{13}} \right) \wedge d \log \left( \frac{Z_{13}}{Z_{12}} \right)  \wedge d \log \left( \frac{Z_{14}}{Z_{12}} \right)
\end{equation}
on~\eqref{equ_affinecharteqq}, which is a holomorphic volume form on the complement of the anti-canonical divisor $\scr{D} = \{ Z_{12} \cdot Z_{23} \cdot Z_{34} \cdot Z_{14} = 0 \}$.

We trivialize the complement by the map
\begin{equation}\label{equ_identificationofremovedd}
\mcal{X}_\varepsilon \backslash \scr{D} \to \left( \C^2 \backslash \{A_1 A_2 = \varepsilon \} \right) \times (\C^\times)^2, \quad [Z_{ij}] \mapsto \left(\frac{Z_{13}}{Z_{12}}, \frac{Z_{24}}{Z_{34}} , \frac{Z_{34}}{Z_{12}}, \frac{Z_{14}}{Z_{34}} \right).
\end{equation}
Equip $\left( \C^2 \backslash \{A_1 A_2 = \varepsilon \} \right) \times (\C^\times)^2$ with the symplectic form inherited from $\mcal{X}_\varepsilon \backslash \scr{D}$ via the map~\eqref{equ_identificationofremovedd}. 
Consider the $\mathbb{T}^3$-action 
$$
\left((\eta_1, \eta_2, \eta_3), (A_1,A_2,A_3,A_4)\right) \mapsto \left(e^{-\sqrt{-1} \eta_1} A_1, e^{\sqrt{-1} \eta_1} A_2, e^{\sqrt{-1} \eta_2} A_3, e^{\sqrt{-1} \eta_3} A_4 \right)
$$
on $\left( \C^2 \backslash \{A_1 A_2 = \varepsilon \} \right) \times (\C^\times)^2$. 
The map~\eqref{equ_identificationofremovedd} is $\mathbb{T}^3$-equivariant $(\theta_2 \circ \theta_4, \theta_4, \theta_1 \circ \theta_4) \mapsto (\eta_1, \eta_2, \eta_3)$. 
Under the identification, the Lagrangian $\scr{L}_i$ can be regarded as the product of $\bL_i\subset \C^2 \backslash \{ab=\epsilon\}$ in Section~\ref{sec_wallcrossingimmersedLag} and $\mathbb{T}^2$-orbits generated by $\eta_2 = \theta_4$ and $\eta_3 = \theta_1 \circ \theta_4$. 
Moreover, the meromorphic volume form $\Omega$ in~\eqref{equ_Omegazgrad} is expressed as
\begin{equation}\label{equ_inducedOmegazgrad}
\Omega = \frac{dA_1 \wedge dA_2}{(A_1 A_2 - \varepsilon)} \wedge d \log A_3 \wedge d \log A_4
\end{equation}
(up to a nonzero scalar).  Then it is obvious that
\begin{lemma}
	The Lagrangians $\scr{L}_i$ for $i=0,1,2$ are graded with respect to the above volume form.
\end{lemma}

The above lemma is useful to find the Maslov index of holomorphic disks bounded by $\scr{L}_i$.

We equip the torus factor $\mathbb{T}^2 \subset (\C^\times)^2$ with flat connections $\nabla^{z,w}$ where $z,w \in \Lambda_{\rm U}$ are the holonomies along two orbits generated by $\theta_1, - \theta_4$ respectively. 
The immersed generators $U$ and $V$ in $\mathbb{L}_0 \subset \C^2 \backslash \{ ab = \varepsilon \}$ induce the degree one generators in $\scr{L}_0$ together with unit classes in the other factors. 
By abuse of notation, they will be still denoted by $U$ and $V$. 
Holonomy variables $x_i$ and $y_i$ $(i = 1, 2)$ are analogous to the same variables appearing in Section \ref{sec_wallcrossingimmersedLag}. 
We then have the \emph{formally} deformed Lagrangians$\colon$
\begin{enumerate}
\item $(\scr{L}_i, b_{\scr{L}_i} := \nabla^{x_i,y_i,z_i,w_i})$ for $i=1,2$, 
\item $(\scr{L}_0, b_{\scr{L}_0} := (uU+vV, \nabla^{z_0,w_0}))$
\end{enumerate}
and, from the discussion in Section \ref{sec_wallcrossingimmersedLag}, we see that their deformation spaces are given by
\begin{equation}\label{equ_deformspacess}
\begin{cases}
\scr{U}_i \simeq (\Lambda_{\rm U})^4 \quad \mbox{for $i=1,2$}, \\
\scr{U}_0 \simeq \left(\Lambda_0 \times \Lambda_+ \right) \cup \left( \Lambda_+ \times \Lambda_0\right) \times (\Lambda_{\rm U})^2.
\end{cases}
\end{equation}

Note that the monotone GC Lagrangian torus fiber in~\eqref{equ_GCsystem} is $\theta_2 \circ \theta_4, \theta_4, \theta_1 \circ \theta_4$ invariant.
Thus, it is fully contained in the level set of~\eqref{equ_hamiltonianfucntiont3}.
Taking the reduction of $\left( \C^2 \backslash \{ab = \varepsilon \} \right) \times (\C^\times)^2$ with respect to  $\eta_1, \eta_2$, and $\eta_3$, we see that the fiber is located over a simple closed curve on the reduced space, equivalently on the base of the projection. 
A Lagrangian isotopy on the base preserving the enclosed area leads to a Hamiltonian isotopy between the GC fiber and $\scr{L}_2$. 

Through the toric degeneration~\eqref{equ_pluckertori}, $\scr{L}_2$ degenerates into a Lagrangian torus, which is Hamiltonian isotopic to the toric fiber over the barycenter of the GC polytope. Since the polytope is reflexive and the toric variety admits a small resolution, Theorem~\ref{theorem_NNUmain} yields 
\begin{equation}\label{equ_uichart}  
\left(\scr{U}_2 := (\Lambda_\textup{U})^4, W_{\scr{L}_2} \right) \quad  W_{\scr{L}_2}(\textbf{\textup{z}}) = \frac{1}{z_{1,2}} + \frac{z_{1,2}}{z_{1,1}} + \frac{z_{1,2}}{z_{2,2}} + \frac{z_{1,1}}{z_{2,1}} + \frac{z_{2,2}}{z_{2,1}} + z_{2,1}.
\end{equation}
Since a  disk of Maslov index zero is bounded by the $\mathbb{S}^1$-orbit generated by $\mu_{\theta_2} + \mu_{\theta_4} = \left( \Phi_{1,2} + \Phi_{2,1} - \Phi_{1,1}\right) - \Phi_{2,2}$ at $Z = - \varepsilon$, its boundary class corresponds to the monomial 
${z_{1,2}\, z_{2,1}}/{z_{1,1}\, z_{2,2}}.$ 
Among four holomorphic disks of Maslov index two over the region bounded by $\gamma$, let us choose the disk which degenerates into the one with the boundary class ${z_{2,2}}/{z_{2,1}}.$ 
Recall that the other factors are generated by $\theta_1$ and $- \theta_4$ in~\eqref{equ_Taction} and their moment maps are respectively written by $\Phi_{1,1}$ and $\Phi_{2,2}$. 
Therefore, we have the following identification 
\begin{equation}\label{equ_relationbetweenyijandxyz}
x_2 = \frac{z_{1,2}\, z_{2,1}}{z_{1,1}\, z_{2,2}},\,\,   y_2 = \frac{z_{2,2}}{z_{2,1}},\,\, z_2 = {z_{1,1}},\, \mbox{ and } w_2 = {z_{2,2}}
\end{equation}
where $x_2,y_2$ are analogous to the same variables appearing in Section \ref{sec_wallcrossingimmersedLag}, and the other two correspond to holonomies in the additional circle factors.

The remaining part of this section will be devoted to verifying that the deformations are weakly unobstructed, and to computing the coordinate changes and the potential functions.

\begin{theorem}\label{theorem:q-iso-hat}
In $\mathrm{Gr}(2, \C^4)$, the potential functions of $\scr{L}_0, \scr{L}_1$, and $\scr{L}_2$ are 
\begin{itemize}
\item $\left( \scr{U}_0 = \left(\Lambda_0 \times \Lambda_+ \right) \cup \left( \Lambda_+ \times \Lambda_0\right) \times (\Lambda_{\rm U})^2, W_{\scr{L}_0} \right)$ where
\begin{equation}\label{equn_potentil0}
W_{\scr{L}_0} (u,v, z_0, w_0) = \frac{v}{(uv-1)z_0}  + u + \frac{uz_0}{w_0} + {v w_0}.
\end{equation}
\item $\left(\scr{U}_1 = (\Lambda_{\rm U})^4, W_{\scr{L}_1} \right)$ where
\begin{equation}\label{equn_potentil1}
W_{\scr{L}_1} (x_1, y_1, z_1, w_1) =  \frac{1}{x_1y_1z_1} +  \frac{1}{y_1z_1} + y_1 + \frac{y_1 z_1}{w_1} + \frac{x_1w_1}{y_1} + \frac{w_1}{y_1}
\end{equation}
\item $\left(\scr{U}_2 = (\Lambda_{\rm U})^4, W_{\scr{L}_2} \right)$ where
\begin{equation}\label{equn_potentil2}
W_{\scr{L}_2} (x_2, y_2, z_2, w_2) =  \frac{1}{x_2 y_2 z_2} + y_2 + x_2 y_2 + \frac{x_2y_2z_2}{w_2} + \frac{y_2 z_2}{w_2} + \frac{w_2}{y_2}.
\end{equation}
\end{itemize}
In particular, they glue to form a partially compactified mirror of $\left(\scr{U}_2 = (\Lambda_{\rm U})^4, W_{\scr{L}_2} \right)$.
\end{theorem}

\begin{remark}
Note that $W_{\scr{L}_0}$ is a rational function. It should be viewed as an analytic continuation of its Taylor expansion, each of whose term genuinely counts holomorphic disks.
See Section~\ref{subsec:interpretps} for more details.
\end{remark}

\begin{remark}\label{remark_referece}
To match up Rietsch's mirror in Section~\ref{sec_chartsmirrorGr(2c4)}, we need to adjust the valuations of coordinates so that the first term of~\eqref{equn_potentil0} and~\eqref{equn_potentil2} and the first two terms of~\eqref{equn_potentil1} are multiplied by $T^4$.
For instance, the appropriate adjustment is 
$$
x_2 \mapsto x_2,\, y_2 \mapsto T^{-1} y_2,\, z_2 \mapsto T^{-2} z_2,\, w_2 \mapsto T^{-2} w_2.  
$$
\end{remark}

\begin{corollary}
	The disk potential~\eqref{equn_potentil0} of the immersed Lagrangian $\scr{L}_0$ in $\Gr(2,\C^4)$ has a critical point, and hence $\scr{L}_0$ has non-trivial Floer cohomology. 
	In particular $\scr{L}_0$ is non-displaceable.
\end{corollary}

We begin by showing that all deformations are unobstructed.

\begin{lemma} \label{lem:prod-unobs}
	In the product $(\C^2 \backslash \{ab=\epsilon\})\times(\C^\times)^2$, for any $x_i, y_i, z_i, w_i \in \Lambda_{\rm U}$ and $u, v \in \left( \Lambda_0 \times \Lambda_+ \right) \cup \left( \Lambda_+ \times \Lambda_0 \right)$ the deformed Lagrangians $(\scr{L}_i, b_{\scr{L}_i})$ for $i=0, 1,2$ are unobstructed.
\end{lemma}
\begin{proof}
	Since each $\scr{L}_i$ does \emph{not} bound any non-constant holomorphic disks, $(\scr{L}_1, b_{\scr{L}_1})$ and $(\scr{L}_2, b_{\scr{L}_2})$ are unobstructed. It remains to check that the Lagrangian $(\scr{L}_0, b_{\scr{L}_0})$ is unobstructed. 
	Because each $b_{\scr{L}_0}$ is of degree one,
	by the dimension reason, the only possible output of the obstruction $m^{(\scr{{L}}_0,b_{\scr{L}_0})}_0$ is of degree two, and hence must be of the form
	$$
		m^{(\scr{{L}}_0,b_{\scr{L}_0})}_0 (\alpha^{\mathbb{L}_0}_0 \otimes  \alpha_0^{\mathbb{T}^2}) = n_{(2,0)} \cdot \left( \beta^{\mathbb{L}_0}_2 \otimes \alpha_0^{\mathbb{T}^2} \right)
	$$
	for some $n_{(2,0)} \in \Lambda_+$. It follows that $n_{(2,0)} = 0$ because such contributions cancel pairwise by Lemma~\ref{lem:L0_unobs}. Thus, the pair $(\scr{{L}}_0,b_{\scr{L}_0})$ is unobstructed.
	\end{proof}

Here is a higher-dimensional analogue of Theorem~\ref{thm:q-iso}.  

\begin{lemma}\label{cor:q-iso-hat} 
In the product $(\C^2 \backslash \{ab=\epsilon\})\times(\C^\times)^2$, there are quasi-isomorphisms among $(\scr{L}_i, b_{\scr{L}_i})$ for $i = 0, 1, 2$ given as follows$\colon$
	\begin{itemize}
		\item $(\scr{L}_1,b_{\scr{L}_1}) \cong (\scr{L}_0,b_{\scr{L}_0})$ if and only if
		\begin{equation}\label{equ_coordinatechange10}
			x_1= uv-1, \,\, y_1= u, \,\, z_1=z_0, \mbox{ and } \, w_1=w_0
		\end{equation} 
		where $u,y \in \Lambda_{\rm U}$, $v \in \Lambda_+$, and $x_1 \in -1 + \Lambda_+$.
		\item $(\scr{L}_0,b_{\scr{L}_0}) \cong (\scr{L}_2,b_{\scr{L}_2})$ if and only if
		\begin{equation}\label{equ_coordinatechange20}		
		uv-1 = x_2, \,\, v^{-1} = y_2, \,\, z_0=z_2, \mbox{ and } \, w_0 = w_2
		\end{equation}
		where $v,y_2 \in \Lambda_{\rm U}$, $u \in \Lambda_+$, and $x_2 \in -1 + \Lambda_+$.
		\item $(\scr{L}_1,b_{\scr{L}_1}) \cong (\scr{L}_2, b_{\scr{L}_2})$ if and only if
		\begin{equation}\label{equ_coordinatechange12}
		x_1= x_2, \,\, y_1 = y_2 (x_2+1), \,\, z_1=z_2, \mbox{ and } \, w_1=w_2
		\end{equation}
		where $y_i \in \Lambda_{\rm U}$ and $x_i \in \Lambda_{\rm U} \backslash (-1 + \Lambda_+)$.
	\end{itemize}
\end{lemma}

\begin{proof}
	By Theorem \ref{thm:q-iso}, the pair $(\alpha_0^{\BL_i,\BL_j},\beta_0^{\BL_j,\BL_i} := \beta_2^{\BL_i,\BL_j})$ of morphisms provides quasi-isomorphisms between $(\BL_i,b_{\BL_i})$ and $(\BL_j,b_{\BL_j})$ for each $(i,j) = (1,2), (1,0),$ and $(0,2)$ under suitable relations among $b_{\BL_\bullet}$.	
	We also have quasi-isomorphisms $(\alpha_{0}^{\mathbb{T}_i,\mathbb{T}_j}, \beta_{0}^{\mathbb{T}_j,\mathbb{T}_i}$) between $(\mathbb{T}_i,b_{\mathbb{T}_i})$ and $(\mathbb{T}_j,b_{\mathbb{T}_j})$ where $\alpha_{0}^{\mathbb{T}_i,\mathbb{T}_j}$ is the maximum point, and $\beta_{0}^{\mathbb{T}_j,\mathbb{T}_i} = \beta_{2}^{\mathbb{T}_i,\mathbb{T}_j}$ is the minimum point (of a certain Morse function on the cleanly intersected torus $\mathbb{T}^2$).
	
	We claim that $\alpha_0^{\BL_1,\BL_0} \otimes \alpha_{0}^{\mathbb{T}_1,\mathbb{T}_0} \in  \mathrm{CF}(\scr{L}_1, \scr{L}_0)$ is a quasi-isomorphism if and only if the bounding cochains are related by the coordinate change~\eqref{equ_coordinatechange10}. 
	Because of the maximum principle, we may analyze Floer differentials and Floer products on each summand of the product $\C^2 \backslash \{ab = \epsilon\}$ and $(\C^\times)^2$ separately. Specifically, the output of $m_1^{\textbf{b}_\scr{L}}(\alpha_0^{\BL_1,\BL_0} \otimes \alpha_{0}^{\mathbb{T}_1,\mathbb{T}_0})$ has degree one and is of the form
$$
m_1^{\textbf{b}_\scr{L}}(\alpha_0^{\BL_1,\BL_0} \otimes \alpha_{0}^{\mathbb{T}_1,\mathbb{T}_0}) = m_1^{\textbf{b}_\mathbb{L}}(\alpha_0^{\BL_1,\BL_0}) \otimes \alpha_{0}^{\mathbb{T}_1,\mathbb{T}_0} \pm \alpha_0^{\BL_1,\BL_0}\otimes m_1^{\textbf{b}_\mathbb{T}}(\alpha_{0}^{\mathbb{T}_1,\mathbb{T}_0})
$$
where $\textbf{b}_\bullet$ is meant to be a sequence of bounding cochains for $\bullet$.  
Because of the area condition on holomorphic sections over the regions bounded by $\gamma_0$ and $\gamma_1$,
	we then have
	\begin{itemize}
	\item $m_1^{\textbf{b}_\mathbb{L}}(\alpha_0^{\BL_1,\BL_0}) = 0$ if and only if $x_1= uv-1$ and $y_1= u$ 
	\item $m_1^{\textbf{b}_\mathbb{T}}(\alpha_{0}^{\mathbb{T}_1,\mathbb{T}_0}) = 0$ if and only if 
	$z_1 = z_0$ and $w_1 = w_0$. 
	\end{itemize}
	The first one follows from Theorem~\ref{thm:q-iso}. 
	
	Also, 	
$$ 
m_2^{\textbf{b}_\scr{L}}(\alpha_0^{\BL_1,\BL_0} \otimes \alpha_0^{\mathbb{T}_1,\mathbb{T}_0},\beta_0^{\BL_0,\BL_1} \otimes \beta_0^{\mathbb{T}_0,\mathbb{T}_1}) = m_2^{\textbf{b}_{\mathbb{L}}}(\alpha_0^{\BL_i,\BL_0},\beta_0^{\BL_0,\BL_i}) \otimes m_2^{\textbf{b}_{\mathbb{T}}}(\alpha_0^{\mathbb{T}_1,\mathbb{T}_0},  \beta_0^{\mathbb{T}_0,\mathbb{T}_1})
$$
where $\beta_0^{\BL_0,\BL_1} := \beta_2^{\BL_1,\BL_0}$ and $\beta_0^{\mathbb{T}_0,\mathbb{T}_1} := \beta_2^{\mathbb{T}_1,\mathbb{T}_0}$.
	Again by Theorem~\ref{thm:q-iso}, we then have
	\begin{itemize}
	\item $m_2^{\textbf{b}_{\mathbb{L}}}(\alpha_0^{\BL_1,\BL_0},\beta_0^{\BL_0,\BL_1}) = \alpha_0^{\mathbb{L}_1}$,
	\item $m_2^{\textbf{b}_{\mathbb{T}}}(\alpha_0^{\mathbb{T}_1,\mathbb{T}_0},  \beta_0^{\mathbb{T}_0,\mathbb{T}_1}) = \alpha_0^{\mathbb{T}_1}$
	\end{itemize}
where $\alpha_0^{{\bullet}}$ is the critical point of degree zero. 
In sum, the condition for $\alpha_0^{\BL_1,\BL_0} \otimes \alpha_{0}^{\mathbb{T}_1,\mathbb{T}_0}$ being a quasi-isomorphism is equivalent to~\eqref{equ_coordinatechange10}.
	The others~\eqref{equ_coordinatechange20} and~\eqref{equ_coordinatechange12} can be similarly shown. 
\end{proof}

From now on, by abuse of notation, $\scr{L}_i$ is regarded as a Lagrangian in $\Gr(2,\C^4)$.

\begin{lemma}\label{lemma_weaklyunobgr24}
	For each $i = 0, 1, 2$, the pair $(\scr{L}_i, {b}_{\scr{L}_i})$ are weakly unobstructed in $\Gr(2, \C^4)$.
\end{lemma}
\begin{proof}
	All stable disks are of Maslov index greater or equal to two because any non-constant disk intersects the anti-canonical divisor $\scr{D}$ and $\Gr(2, \C^4)$ is Fano. Since $b$ is a combination of degree one elements, disks of Maslov index two can only contribute to $m_0^{b_{\scr{L}_\bullet}}$ by the dimension reason, and hence the output must have degree zero.  As we have chosen perfect Morse functions on the Lagrangians, we have a \emph{unique} degree zero element $\alpha_0^{\scr{L}_\bullet}$ in the Morse model, which is indeed the unit of the $A_\infty$-algebra on $\scr{L}_\bullet$.  Hence $m_0^{b_{\scr{L}_\bullet}}$ is proportional to the unit.
\end{proof}

\begin{proposition}\label{proposition_quasisogr24}
	Regarding $\scr{L}_i$'s as Lagrangians in $\Gr(2,\C^4)$, 
	the statements on quasi-isomorphisms in Lemma~\ref{cor:q-iso-hat} still hold.
\end{proposition}
\begin{proof}
	Any ambient disks of Maslov index greater than or equal to two 
	 contribute to neither $m_1^{\textbf{b}_\scr{L}}(\alpha_0^{\BL_i,\BL_j} \otimes \alpha_{0}^{\mathbb{T}_i,\mathbb{T}_j})$ nor  $m_2^{\textbf{b}_\scr{L}}(\alpha_0^{\BL_i,\BL_j} \otimes \alpha_0^{\mathbb{T}_i,\mathbb{T}_j},\beta_0^{\BL_j,\BL_i} \otimes \beta_{0}^{\mathbb{T}_j,\mathbb{T}_i})$, since the inputs have degree zero.  The conditions for quasi-isomorphisms remain same after the compactification.
\end{proof}

Now, we are ready to prove Theorem~\ref{theorem:q-iso-hat}. 

\begin{proof}[Proof of Theorem~\ref{theorem:q-iso-hat}]
In~\eqref{equ_uichart} and~\eqref{equ_relationbetweenyijandxyz}, we have obtained  the mirror $\left( \scr{U}_2, W_{\scr{L}_2} \right)$ of $\scr{L}_2$. 
We can compute the LG mirrors associated to $\scr{L}_1$ and $\scr{L}_0$ from the formula of $W_{\scr{L}_2}$ for $\scr{L}_2$
with help of Lemma~\ref{lemma_weaklyunobgr24} and Proposition~\ref{proposition_quasisogr24} as follows.
By the $A_\infty$-relation, 
\begin{align*}
(m^{\textbf{b}_\scr{L}}_1 \circ m^{\textbf{b}_\scr{L}}_1) (\alpha_0^{\scr{L}_2, \scr{L}_0}) &=  m^{\textbf{b}_\scr{L}}_2 ( \alpha_0^{\scr{L}_2, \scr{L}_0}, m^{\textbf{b}_{\scr{L}}}_0 (1)) - m^{\textbf{b}_\scr{L}}_2 ( m^{\textbf{b}_{\scr{L}}}_0 (1), \alpha_0^{\scr{L}_2, \scr{L}_0}) \\ 
&= ( W_{\scr{L}_0}({b}_{\scr{L}_0}) - W_{\scr{L}_2}({b}_{\scr{L}_2})) \cdot \alpha_0^{\scr{L}_2, \scr{L}_0}
\end{align*}
where $\alpha_0^{\scr{L}_2, \scr{L}_0} := \alpha_0^{\BL_2,\BL_0} \otimes \alpha_{0}^{\mathbb{T}_2,\mathbb{T}_0}$.
Therefore, on open dense subsets of the domains of $W_{\scr{L}_0}$ and $W_{\scr{L}_2}$, two potential functions agree
$$
W_{\scr{L}_2}({b}_{\scr{L}_2}) = W_{\scr{L}_0}({b}_{\scr{L}_0}),
$$
provided that the coordinate change between ${b}_{\scr{L}_2}$ and ${b}_{\scr{L}_0}$ is in ~\eqref{equ_coordinatechange20}. By applying the coordinate change~\eqref{equ_coordinatechange20}, the expression~\eqref{equn_potentil0} is derived. 
Since it coincides with the potential function $W_{\scr{L}_0}$ on an open dense subset, and the potential is analytic, it gives the potential function globally on $\scr{U}_0$. Similarly, we can compute $W_{\scr{L}_1}$. 
\end{proof}

\begin{remark}
	The above argument works in general, and proves that disk potentials for two Lagrangians (with Maurer-Cartan deformations) agree (i.e. $W_{\scr{L}_0}(b_{\scr{L}_0})=W_{\scr{L}_1}(b_{\scr{L}_1})$), if there is a non-zero Floer cohomology class in $\mathrm{HF}((L_0,b_{\scr{L}_0}),(L_1,b_{\scr{L}_1}))$.  Similar argument appeared in \cite{PT}.
\end{remark}

\subsection{Interpretations}\label{subsec:interpretps}
\begin{figure}[h]
	\begin{center}
		\includegraphics[scale=0.55]{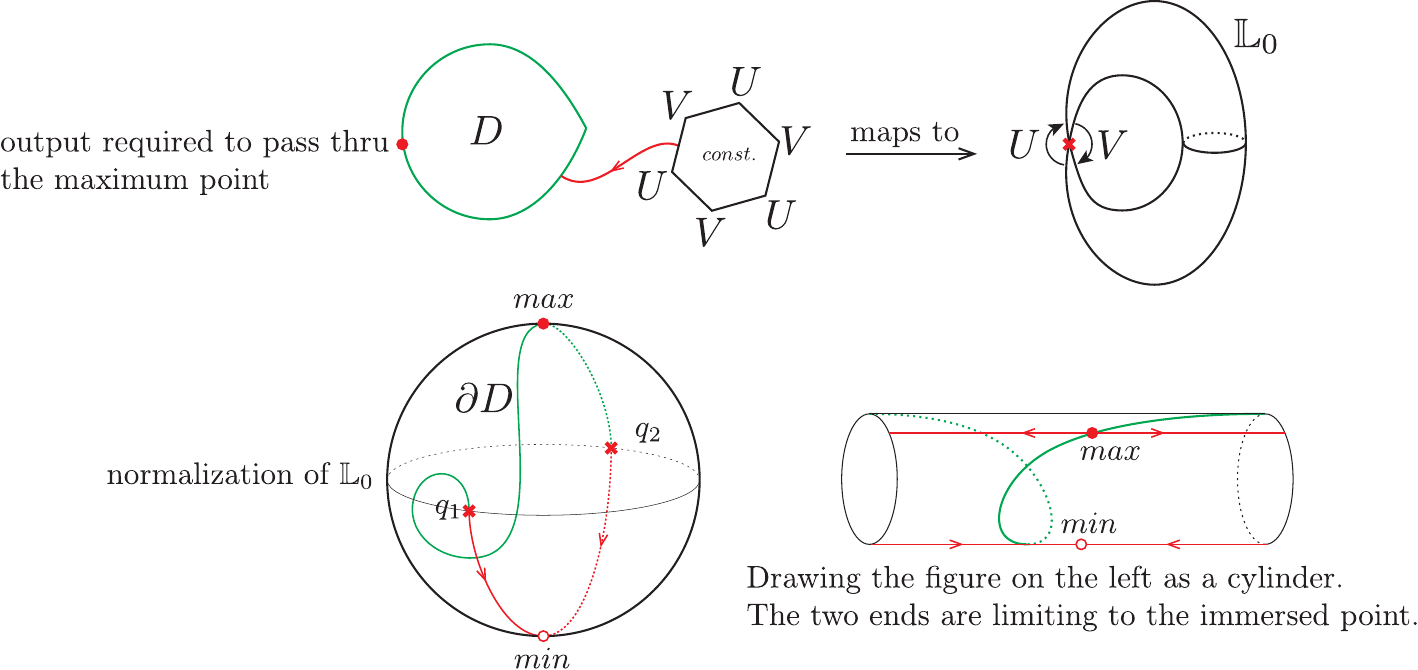}
		\caption{An example of a pearl trajectory bounded by the immersed Lagrangian sphere passing through a generic marked point (which is the maximum point of a generic Morse function).  Constant disk bubbles can be attached to a Maslov-two disk $D$ via Morse flow lines, if $\partial D$ hits a flow line from (preimages of) the immersed point to the minimum point of the Morse function.}
		\label{Fig_pearl-W-imm}
	\end{center}
\end{figure}

The disk potential $W_{\scr{L}_0}$ for our immersed Lagrangian $\scr{L}_0$ in $\Gr(2,\C^4)$ is a polynomial in $u,v,(uv-1)^{-1}$, which in turn is an infinite series in $u,v$.  They are contributed from pearl trajectories as shown in the top of Figure \ref{Fig_pearl-W-imm}.  Since constant disk bubbles can pass through the immersed generators $u,v$ arbitrarily many times without affecting the Maslov index (and hence the moduli dimension), arbitrarily high powers of $u,v$ appear in $W$.  (We only depict the factor of the immersed sphere, since the other torus factor has trivial effect on wall-crossing.)

In Equation \eqref{equn_potentil0}, the Maslov-two disks $u, u z_0 w_0^{-1}, {v}{w_0}$ do not interact with constant disk bubbles and hence do not involve the factor $(uv-1)^{-1}$, while the Maslov-two disk $vw_0$ admits disk bubbling with constant disks.  

To explain this, consider the bottom of Figure \ref{Fig_pearl-W-imm}, where a normalization of the immersed sphere is drawn (which is a sphere with two marked points $q_+,q_-$, or simply a cylinder).  Denote by $\overline{q_+q_-}$ the arc joining $q_+,q_-$ formed from the union of the two flow line from $q_+$ or $q_-$ to the minimum point of the Morse function.

A Maslov-two disk $D$ (passing through the maximum point of a generic perfect Morse function on the sphere) admits disk bubbling only if its boundary intersects with $\overline{q_+q_-}$.  This condition corresponds to that the Maslov-two disk bounded by a torus fiber has its boundary winding around the vanishing circle (degenerating to the immersed point).  This explains why the disk $vz_0^{-1}w_0^{-1}$ admits disk bubbling while the others do not.

\begin{remark}
	The contributions from constant disk bubbles are very hard to compute directly from analysis of Kuranishi structures.  On the other hand, using the disk potential of a smooth torus fiber and the quasi-isomorphism with the immersed fiber, we obtain the disk potential and conclude that the constant disk bubbles at the immersed point contribute to each disk term by multiplication of $(uv-1)^{-k}$, where $k$ is the signed intersection number of the Maslov-two disk boundary with the arc $\overline{q_+q_-}$.  This computation works for the disk potential of an immersed sphere (or its product with a torus) in general.
\end{remark}

\begin{remark} \label{rmk:t-power}
	In \cite{DET}, they define the Floer theory of the immersed sphere using the generators $u,v,t$ with the relation $uv=1+t$.  They prove that the disk potential $W$ is an invariant of the immersed sphere.
	
	Here we use a perfect Morse function and adopt the formulation of \cite{FOOO-can,AJ,CHLabc,CHLnc,HL18}.  The chosen perfect Morse function provides a minimal model which is unique up to $A_\infty$ isomorphisms.   The $t$-power in the formulation of \cite{DET} corresponds to the (signed) intersection number explained above.  
\end{remark}

\begin{remark}\label{rmk:intrinsicmc}
	Note that the expression of $W$ in $(u,v)$ depends on the choice of the Morse function.   Namely, we can pull back the function by a self-diffeomorphism of the cylinder which rotates one end by $2k\pi$ (and keeps the other end fixed).  Then $W$ would change by $(u,v)\mapsto ((uv-1)^k u,(uv-1)^{-k} v)$, since the intersection number of a disk boundary with the arc $\overline{q_+q_-}$ changes by $k$.
	
	Thus $W$ should be understood as a function on
	$$ \mathcal{M} := \left(\coprod_i \left(\Lambda_+^2\right)_{(i)}\right) \big/ \langle (u_{(i)},v_{(i)}) \sim ((uv-1)^{j-i} u_{(j)},(uv-1)^{i-j} v_{(j)}) \rangle \cong \Lambda_+^2.  $$
	This gives a more coordinate-free description of the weak Maurer-Cartan deformation space of the immersed sphere.  Dimitroglou-Ekholm-Tonkonog \cite{DET} makes sense of this space in terms of a certain path space between $q_+$ and $q_-$ called the string topology algebra.
\end{remark}

\subsection{Identification with the Rietsch's mirror of $\mathrm{Gr}(2, \C^4)$}\label{sec_chartsmirrorGr(2c4)}

In this section, we show that the glued mirror  from Maurer-Cartan deformations of $ \scr{L}_{0}, \scr{L}_1$, and $\scr{L}_2$ in the previous subsection gives the full Riestch's mirror of $\mathrm{Gr}(2, \C^4)$.

As in Example~\ref{exam_rietschgr24}.
recall that the Lie-theoretic mirror in \cite{Rie,MR} is the affine variety given by  
$\check{X} := \mathrm{Gr}(2, \C^4) \backslash \scr{D}$ where ${\scr{D}} := \{p_{12} \cdot p_{23} \cdot p_{34} \cdot p_{14} = 0\}$ is an anti-canonical divisor together with the superpotential $W_{\textup{Rie}}$ in~\eqref{equ_rietschpogr24}.

It has a cluster algebra structure with the frozen variables $\{p_{12}, p_{23}, p_{34}, p_{14}\}$ and the cluster variables $\{p_{13}, p_{24}\}$. 
We consider two cluster charts for $\check{X}$ as follows$\colon$
\begin{equation}\label{equ_algebraictori}
\begin{cases}
\mcal{U}_{1} = \{ [p_{ij}] \in \check{X} ~\colon~ p_{13} \neq 0 \} \\
\mcal{U}_{2} = \{ [p_{ij}] \in \check{X} ~\colon~ p_{24} \neq 0 \}.
\end{cases}
\end{equation}
Restrict to each chart $\mcal{U}_i$ and using $p_{13} p_{24} = p_{12} p_{34} + p_{14} p_{23}$,
\eqref{equ_rietschpogr24} can be expressed as a Laurent polynomial in terms of the frozen variables $p_{12}, p_{23}, p_{34}, p_{14}$ and one cluster variable.  Explicitly, they are given as 
\begin{equation}\label{equ_Laurentpoly}
\begin{cases}
\displaystyle W_{1} := W_{\textup{Rie}} \big{|}_{\mcal{U}_1} \colon \mcal{U}_1 \to \C, \, W_{1} = q \, \frac{p_{34}}{p_{13}} + q \, \frac{p_{14}p_{23}}{p_{12}p_{13}} + \frac{p_{13}}{p_{23}} + \frac{p_{12}}{p_{13}} + \frac{p_{14}p_{23}}{p_{34}p_{13}} + \frac{p_{13}}{p_{14}} \\ \\
\displaystyle W_{2} := W_{\textup{Rie}} \big{|}_{\mcal{U}_2} \colon \mcal{U}_2 \to \C, \, W_{2} = q \, \frac{p_{24}}{p_{12}} + \frac{p_{12} \, p_{34}}{p_{23} \, p_{24}} + \frac{p_{14}}{p_{24}} + \frac{p_{24}}{p_{34}} + \frac{p_{12} \, p_{34}}{p_{14} \, p_{24}} + \frac{p_{23}}{p_{24}}.
\end{cases}
\end{equation}
In this case, the immersed chart is the pair $(\mcal{U}_{\{(1,2)\}}, W_{\textup{Rie}})$. Set $\mcal{U}_0 := \mcal{U}_{\{(1,2)\}}$ for simplicity.

The Laurent polynomial in~\eqref{equ_Laurentpoly} only has four critical points, while the rank of the quantum cohomology ring of $\mathrm{Gr}(2, \C^4)$ is six. On the contrary, the key feature of the LG model $(\check{X}, W_{\textup{Rie}})$ in~\eqref{equ_rietschpogr24} is that its Jacobian ring has rank six, and is isomorphic to the quantum cohomology ring. 

\begin{lemma}\label{lemma_critcal24}
The superpotential \eqref{equ_rietschpogr24} has six critical points$\colon$ Let $\xi$ be the $4$-th root of unity. 
\begin{enumerate}
\item 
For $j = 0, 1, 2, 3$, 
$$
p_{12} = q,\, p_{13} = \sqrt{2} \, \xi^{-j} q^{\frac{3}{4}}, p_{14} = p_{23} =  \xi^{-2j}  q^{\frac{1}{2}}, p_{24} = \sqrt{2} \, \xi^{j} q^{\frac{1}{4}},\,\, p_{34} = 1.  
$$ 
The critical values are respectively $4 \sqrt{2} \, \xi^j \, q^{\frac{1}{4}}$.

\item For $j = 0, 1$, 
$$
p_{13} = 0,\, p_{24} = 0,\, p_{14} = \sqrt{-1} \, \xi^{2j} q^{\frac{1}{2}},\, p_{23} = - \sqrt{-1} \, \xi^{2j} {q^{\frac{1}{2}}},\, p_{34} = 1,\, p_{12} = - q.  
$$ 
The critical values are $0$. 
\end{enumerate}
\end{lemma}

Keeping the cluster structure of LG models on $\check{X}$ in mind, we now explain how to recover Rietsch's mirror using  Floer theory. 
The main theorem of this section is the following.

\begin{theorem}\label{theorem_RietschSYZ}
The Rietsch's mirror $(\check{X}, W_{\textup{Rie}})$ of the Grassmannian $X = \mathrm{Gr}(2, \C^4)$ is isomorphic to the LG model glued from Maurer-Cartan deformation spaces of the monotone Lagrangian tori $\scr{L}_1$, $\scr{L}_2$, and 
the immersed Lagrangian $\scr{L}_{0}$ via quasi-isomorphisms given in Proposition \ref{proposition_quasisogr24}.
\end{theorem}

\begin{proof}
We identify the variables corresponding to the orbits generated by $\theta_1$ and $- \theta_4$ as follows. 
\begin{equation}\label{equ_coorchangezw}
z_0 = z_1 = z_2 = z_{1,1} \mapsto  \frac{p_{23}}{p_{34}}, \quad w_0 = w_1 = w_2 = z_{2,2} \mapsto \frac{p_{14}}{p_{34}}.
\end{equation}
In addition, the identification for the other variables is taken as follows$\colon$
\begin{equation}\label{equ_coorchangeuvu}
 x_1 = x_2 = uv - 1 = \frac{z_{1,2}\, z_{2,1}}{z_{1,1}\, z_{2,2}} \mapsto \frac{p_{12}p_{34}}{p_{14}p_{23}}, \,  y_2 = v^{-1} = \frac{z_{2,2}}{z_{2,1}} \mapsto \frac{p_{14} }{p_{24}}, \, y_1 =  u \mapsto \frac{p_{13}}{p_{23}}.
\end{equation}

Under the coordinate change~\eqref{equ_coorchangezw} and~\eqref{equ_coorchangeuvu}, we express $W_{\scr{L}_i}$ for $i = 1, 2$ in terms of $p_{i,j}$. 
It is straightforward to see that two superpotentials coincide if ignoring $q$.
As we have seen in Lemma~\ref{lemma_critcal24}, the valuations of critical points $\{[p_{ij}] = [s_{ij}]\}$ in $W_{i}$'s are not zero.
We need to adjust the valuation of coordinates, see Remark~\ref{remark_referece}.
Then $W_{\scr{L}_0}$ coincides with the mirror superpotential $W_\textup{Rie}$ in~\eqref{equ_rietschpogr24} (up to some a constant multiple) as desired. 

Recall that the deformation space $\scr{U}_{0}$ for $W_{\scr{L}_0}$ was restricted to $(u,v) \in \Lambda_0 \times  \Lambda_+ \cup \Lambda_+ \times \Lambda_0$ to ensure the convergence.  
As a consequence, it cannot fully cover $\check{X}$. 
The missing part 
$$
\left\{ [p_{ij}] \in \check{X} ~\colon~  p_{13} \in \C^\times, p_{24} = 0 \right\} \cup \left\{ [p_{ij}] \in \check{X} ~\colon~  p_{13} = 0, p_{24}  \in \C^\times \right\}, 
$$
is covered by $\scr{U}_1$ and $\scr{U}_2$. Therefore, the glued mirror from $\scr{U}_0, \scr{U}_1$, and $\scr{U}_2$ covers the entire $\check{X}$.
\end{proof}

\section{Complete SYZ mirror of $\mathrm{Gr}(2,\C^n)$}\label{section_mirrorcongr2n}
Now we study the Grassmannians of 2-planes in higher dimensions, and construct their mirrors using similar techniques developed in previous sections. 
We complete the local mirror of a smooth Lagrangian torus of $\mathrm{Gr}(2,\C^n)$ by gluing in the mirror charts from suitably constructed immersed Lagrangians. 
We shall show that the glued LG mirror of $\mathrm{Gr}(2,\C^n)$ coincides with the Rietsch's mirror. The main argument is in parallel to the one for $\mathrm{Gr}(2,\C^4)$ except some complications mainly caused by high dimensionality. Another subtlety in constructing Lagrangians lies the choice of symplectic forms on a local model, which we did not have before in low dimensions (Remark \ref{remark_issuee}).

\subsection{Construction of immersed Lagrangians}\label{sec_constuofimmLag}

Let $X$ be the Grassmannian $\mathrm{Gr}(2,\C^n)$ of 2-planes. 
To produce Lagrangians in the complement of the anti-canonical divisor $\scr{D} := \{Z_{1,2} \cdot Z_{2,3} \cdots Z_{n-1,n}\cdot  Z_{1,n} = 0\}$ in $X$, we shall focus on the local chart of $\mathrm{Gr}(2,\C^n)$ around the monotone Lagrangian non-torus fiber associated with $\mcal{I} \in \scr{I}_n$ in Section~\ref{subsection_blockcombi}.
It is then realized as a product of local charts of $\mathrm{Gr}(2,\C^4)$. 
Taking a product of the Lagrangians constructed in Section~\ref{sec_consimmergr24}, we obtain  Lagrangians in $\mathrm{Gr}(2,\C^n)$ that are ingredients of our mirror construction.

\begin{remark}\label{remark_issuee}
One may try to construct Lagrangians by taking the reduction of a maximal torus action on $\mathrm{Gr}(2,\C^n)$ for $n > 4$ as in Section~\ref{sec_consimmergr24}. However, construction of Lagrangians in the reduced space, which has complex dimension $\geq 2$ and is equipped with a reduced symplectic form, is more difficult. 
\end{remark}

Details on our construction are in order.
The first task is to equip $\mathrm{Gr}(2, \C^n)$ with a meromorphic volume form. This is crucial to obtain the global mirror later. Let us restrict ourselves to the chart $(\C^\times)^{2n-4}$ given by 
\begin{equation}\label{equ_clusterchartempt}
\prod_{j=1}^{n-1} Z_{j,j+1} \times \prod_{j=3}^n Z_{1,j} \neq 0.
\end{equation}
Consider the holomorphic volume form 
\begin{equation}\label{equation_merovolgr2n}
\Omega := \bigwedge_{j = 1}^{n-2} \left( d \log \left(  \frac{Z_{j+1,j+2}}{Z_{1,j+1}} \right) \wedge  d \log \left(  \frac{Z_{1,j+2}}{Z_{1,2}} \right) \right)
\end{equation}
on the algebraic torus $(\C^\times)^{2n-4}$, which extends to a meromorphic volume form on $\mathrm{Gr}(2, \C^n)$. 
The extended form has a simple pole on the anti-canonical divisor $\scr{D}$ as observed in  \cite[Lemma 8.5]{MR}.

For each index $\mcal{I}$ in $\scr{I}_n$, consider the trivialization
\begin{equation}\label{equ_complexide}
\Upsilon_\textup{cpx} \colon 
\mathrm{Gr}(2, \C^n) \backslash \scr{D}^\mcal{I} \to
\prod_{(i,i+1) \in \mcal{I}}
(\C^2 \backslash \{A_{i,1} A_{i,2} = \varepsilon \} \times (\C^\times)^2)
\times (\C^\times)^{2n -4 - 4 \cdot |\mcal{I}|}
\end{equation}
defined by
\begin{equation}\label{equ_cpxmapaa}
[Z_{r,s}] \mapsto 
\begin{cases}
\left( \frac{Z_{1,i+2}}{Z_{1,i+1}}, \frac{Z_{i+1, i+3}}{Z_{i+2,i+3}}, \frac{Z_{i+2,i+3}}{Z_{1,i+1}}, \frac{Z_{1,i+3}}{Z_{i+2,i+3}} ~\colon~ (i,i+1) \in \mcal{I} \right) \\
\left( \frac{Z_{i+1,i+2}}{Z_{1,i+1}}, \frac{Z_{1,i+2}}{Z_{1,2}} \colon (i, i+1), (i-1, i) \notin \mcal{I}  \right).
\end{cases}
\end{equation}
where $\mathrm{Gr}(2, \C^n) \backslash \scr{D}^\mcal{I}$ is the chart of $\mathrm{Gr}(2, \C^n) \backslash \scr{D}$ on which the denominators in~\eqref{equ_cpxmapaa} and $\C^\times$-factors do not vanish.  
Here, ($A_{i,1}, A_{i,2}, A_{i,3}, A_{i,4})$ and $(B_{1,i}, B_{2,i})$ are the complex coordinates on the target space.

Consider 
\begin{equation}\label{equ_productymcali}
\mcal{Y}_\mcal{I}:= (\mathrm{Gr}(2, \C^4) \backslash \{Z_{1,2} \cdot Z_{2,3} \cdot Z_{3,4} \cdot Z_{1,4} = 0\})^{|\mcal{I}|} \times (\mathrm{Gr}(1, \C^2) \backslash \{ Z_0 \cdot Z_1 = 0 \})^{2n - 4 - 4 \cdot |\mcal{I}|},
\end{equation}
which is diffeomorphic to $\Upsilon_\textup{cpx} (\mathrm{Gr}(2, \C^n) \backslash \scr{D}^\mcal{I})$ 
because 
$$
\mathrm{Gr}(2, \C^4) \backslash \{Z_{1,2} \cdot Z_{2,3} \cdot Z_{3,4} \cdot Z_{1,4} = 0\} \simeq \C^2 \backslash \{A_{i,1} A_{i,2} = \varepsilon \} \times (\C^\times)^2.
$$
The local chart $\mcal{Y}_\mcal{I}$ admits two symplectic forms$\colon$
\begin{enumerate}
\item $\omega_{(2,n)}$ is denoted by the K{\"a}hler form inherited from $\mathrm{Gr}(2,\C^n)$,
\item $\omega_\mcal{I}$ is denoted by the product K{\"a}hler form.
\end{enumerate}
Since $\mathrm{Gr}(2, \C^n)$ and $\mathrm{Gr}(2, \C^4)^{|\mcal{I}|} \times \mathrm{Gr}(1, \C^2)^{2n - 4 - 4 \cdot |\mcal{I}|}$ are locally the cotangent bundle of $(\mathbb{S}^3)^{|\mcal{I}|} \times \mathbb{T}^{2n- 4 - 3 \cdot |\mcal{I}|}$,
$(\mcal{Y}_\mcal{I}, \omega_\mcal{I})$ can be embedded into $(\mcal{Y}_\mcal{I}, \omega_{(2,n)})$ symplectically (by rescaling symplectic form $\omega_\mcal{I}$ if necessary). 
We then extend $\omega_\mcal{I}$ into one on $\mcal{Y}_\mcal{I}$ by interpolating the pull-backed symplectic form and the symplectic form $\omega_{2,n}$ equivariantly.
The extended symplectic form on $\mcal{Y}_\mcal{I}$ is still denoted by $\omega_\mcal{I}$ by abuse of notation. 
By applying an equivariant version of the Moser argument, we obtain a $\mathbb{T}^{2n - 4 -|\mcal{I} |}$-equivariant endomorphism $\Upsilon_\textup{sym}$ on $\mcal{Y}_\mcal{I}$ such that  $\Upsilon_\textup{sym}^* \omega_{(2,n)} = \omega_\mcal{I}$. 

Applying the construction in Section~\ref{sec_LagfibonF3}, we produce Lagrangians in each local $\mathrm{Gr}(2, \C^4)$-factor.
Specifically, for $(i, i+1) \in \mcal{I}$, consider the Hamiltonian $\mathbb{T}^3$-action generated by~\eqref{equ_hamiltonianfucntiont3} on $\mathrm{Gr}(2,\C^4)$. 
Take the reduction at the level at the origin in~\eqref{equ_hamiltonianfucntiont3}.  
As in Figure~\ref{Fig_baseforGr24}, we then consider three simple closed curves in the base, say $\gamma_{i,k}$ for $k = 0, 1, 2$. 
Collecting the ${\mathbb{T}^{3}}$-orbits over $\gamma_{i,k}$ produces a Lagrangian in $\mathrm{Gr}(2, \C^4)$. Also, we take a circle whose center is the origin in each $\C^\times$-summand.
The product of the Lagrangians and circles gives rise to a Lagrangian $\scr{L}_k^{\mcal{I}}$ for $k= 0, 1, 2$ in $(\mcal{Y}_\mcal{I}, \omega_\mcal{I})$, and hence in $(\mcal{Y}_\mcal{I}, \omega_{(2,n)})$ via $\Upsilon_\textup{sym}$. 
By taking a suitable size of circles, the Lagrangians are chosen to be monotone.

In local coordinates given in~\eqref{equ_cpxmapaa}, the volume form~\eqref{equation_merovolgr2n} can be written as (up to a non-zero scalar multiple) a holomorphic volume form on $\mcal{Y}_\mcal{I}$ by
\begin{equation}\label{equ_Zgradedvolmloca}
\Omega = \left( \bigwedge_{(i,i+1) \in \mcal{I}} \frac{dA_{i,1} \wedge dA_{i,2} \wedge dA_{i,3} \wedge dA_{i,4}}{(A_{i,1} A_{i,2} - 1) A_{i,3} A_{i,4}} \right) \wedge  \left( \bigwedge_{(i-1, i), (i,i+1) \notin \mcal{I}} 
\frac{ d B_{1,i} \wedge d B_{2,i}}  {B_{1,i} B_{2,i}} \right).
\end{equation}
The Lagrangians on $(\mcal{Y}_\mcal{I}, \omega_{(2,n)})$ constructed above admit  $\Z$-gradings with respect to $\Omega$, which will be used for a local calculation later. Because the grading is topological, the $\Z$-grading is preserved through the diffeomorphism $\Upsilon_\textup{sym}$ so that it makes sense in $(\mcal{Y}_\mcal{I}, \omega_\mcal{I})$ as well.

\begin{remark}
The identification~\eqref{equ_complexide} was used to interpolate generalized GC systems in Nohara-Ueda~\cite{NUclu}.
\end{remark}

\subsection{Local SYZ mirrors of immersed Lagrangians}\label{completesyzmirrgr2n}

In this section, we investigate Floer theoretical relation between the Lagrangians $\scr{L}^{\mcal{I}}_k$ for $\mcal{I} \in \scr{I}_n$ and $k = 0 ,1 ,2$ constructed in~\ref{sec_constuofimmLag}. 

We first explain a combinatorial algorithm how to read the potential functions of immersed Lagrangian $\scr{L}^\mcal{I}_0$. 
Let us start from the potential function $W_\scr{T}$ in~\eqref{equ_potentialGC}. 
\begin{itemize}
\item (Step 1) For each $(i,i+1) \in \mcal{I}$, get rid of the following six terms from $W_\scr{T}$
$$
\frac{z_{1,i+2}}{z_{1,i+1}} + \frac{z_{1,i+1}}{z_{1,i}} + \frac{z_{2,i+1}}{z_{2,i}} + \frac{z_{1,i+1}}{z_{2,i+1}} + \frac{z_{1,i}}{z_{2,i}} + \frac{z_{2,i}}{z_{2,i-1}}. 
$$
In terms of block combinatorics, those terms correspond to the edges in the $\mathrm{U}(2)$-block associated with $(i,i+1)$, see Figure~\ref{fig_potentialimmer}.

\vspace{-0.2cm}
\begin{figure}[h]
	\begin{center}
	\includegraphics[scale=0.6]{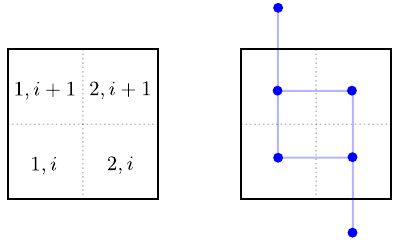}
	\vspace{-0.1cm}
	\caption{\label{fig_potentialimmer} Step 1 for the disk potential function of $\scr{L}^\mcal{I}_0$.}	
	\end{center}
\end{figure}
\vspace{-0.1cm}

\item (Step 2) For each $(i,i+1) \in \mcal{I}$, insert the following four terms to the output of the first step.
$$
u_i + u_i \frac{z_{1,i}}{z_{2,i+1}} + v_i \frac{z_{2,i+1}}{z_{2,i-1}} + \frac{v_i z_{1,i+2}}{(u_iv_i - 1)z_{1,i}}.
$$
\end{itemize}

We denote by $W_{\mcal{I}}$ the resulting potential function. 

\begin{example}\label{example_immersedpote}
Let us revisit Example~\ref{exa_fiberpofgr26} and Figure~\ref{fig_gr26Lagfaces}. 
By applying the above algorithm, $W_{(2,3)}$ and $W_{(1,2),(3,4)}$ can be written as follows. 
$$
W_{(2,3)} = u_2 + u_2 \frac{z_{12}}{z_{23}} + v_2 \frac{z_{23}}{z_{21}} + \frac{v_2 z_{14}}{(u_2v_2 - 1) z_{12}} + 
\frac{T^6}{z_{14}} + \frac{z_{14}}{z_{24}} + \frac{z_{24}}{z_{23}} + \frac{z_{12}}{z_{11}} + \frac{z_{11}}{z_{21}} + \frac{1}{z_{21}},
$$
and
$$
W_{(1,2),(3,4)} = u_1 + u_1 \frac{z_{11}}{z_{22}} + v_1 z_{22} + v_1 \frac{z_{13}}{(u_1v_1-1)z_{11}} + u_3 + u_3 \frac{z_{13}}{z_{24}} + v_3 \frac{z_{24}}{z_{22}} + v_3 \frac{T^6}{(u_3v_3 - 1)z_{13}}.
$$
\end{example}

The main theorem of this section states as follows$\colon$

\begin{theorem}\label{theroem_diskpotentialimgr26}
The superpotential $W_{\mcal{I}}$ is a Floer theoretical disk potential of the immersed Lagrangian $\scr{L}_0^\mcal{I}$. 
\end{theorem}

Fix $\mcal{I} \in \scr{I}_n$. We first compute quasi-isomorphisms between the objects $\scr{L}_0^{\mcal{I}}, \scr{L}_1^{\mcal{I}}$, and $\scr{L}_2^{\mcal{I}}$.
In order to fix coordinates on the {formal} deformation spaces of these Lagrangians, it suffices to equip the carried Lagrangians with a flat $\Lambda_\mathrm{U}$-bundle in~\eqref{equ_productymcali}.
For each $k = 0, 1 ,2$ and $(i,i+1) \in \mcal{I}$,
a flat $\Lambda_{\mathrm{U}}$-connection $\nabla^{\textbf{z}}_{2,i+1,k}$ (resp. $\nabla^{\textbf{z}}_{1,i,k}$) is adorned on $\scr{L}^{\mcal{I}}_k$
to realize the holonomy $z_{2,i+1,k}$ (resp. $z_{1,i,k}$) along the circle-factor in the plane $A_{i,3}$ (resp. $A_{i,4}$) with the counter-clockwise orientation. 
Similarly, flat $\Lambda_{\mathrm{U}}$-connections $\nabla^{\textbf{z}}_{1,i,k}$ and $\nabla^{\textbf{z}}_{2,i,k}$ are equipped for the holonomies $z_{1,i,k}$ and $z_{2,i,k}$ along the circle-factors of $B_{1,i}$ and $B_{2,i}$ with the counter-clockwise orientation. 

For each summand $\C^2 \backslash \{ A_{i,1} A_{i,2} = \varepsilon \}$, we define flat connections and holonomy variables (denoted by $x_{i,k}$ and $y_{i,k}$) as in Section~\ref{subsec:supple} that lead to the formulae in Theorem~\ref{thm:q-iso}.
Let $U_{i}$ and $V_{i}$ be the products of immersed generators lying over $\gamma_{i,0}$ and the unit classes in the other factors, both of which are of degree one. 
The corresponding component of the Maurer-Cartan deformation space is $(\Lambda_0 \times \Lambda_+) \cup (\Lambda_+ \times \Lambda_0)$, whose coordinates are taken to be $(u_i, v_i)$.
We then have the following boundary-deformed Lagrangians$\colon$ 
\begin{enumerate}
\item  two Lagrangian tori $\left( \scr{L}_k^\mcal{I}, b^\mcal{I}_k := \nabla^{\textbf{x}, \textbf{y}, \textbf{z}} \right)$ for $k=1,2$ whose associated deformation spaces are $\scr{U}_k^\mcal{I} \simeq (\Lambda_\textup{U})^{2n -4}$,
\item one immersed Lagrangian 
$$
\left( \scr{L}_0^\mcal{I}, b^\mcal{I}_0 := \left( \sum_{(i,i+1) \in \mcal{I}} u_{i} U_{i} + v_{i} V_{i}, \nabla^{\textbf{z}} \right) \right)
$$
with the deformation space $\scr{U}^\mcal{I}_0 \simeq \left((\Lambda_0 \times \Lambda_+) \cup (\Lambda_+ \times \Lambda_0)\right)^{| \mcal{I} |} \times (\Lambda_\textup{U})^{2n - 4 - 2 \cdot |\mcal{I}|}.$
\end{enumerate}

\begin{lemma}\label{lemma_localchartcomq}
For each fixed $\mcal{I} \in \scr{I}_n$, consider a monotone immersed Lagrangian $\scr{L}^\mcal{I}_0$ and monotone Lagrangian tori $\scr{L}^\mcal{I}_1, \scr{L}^\mcal{I}_2$ in a local chart. Then their deformation spaces are related by quasi-isomorphisms between $(\scr{L}^{\mcal{I}}_k, b^\mcal{I}_{k})$ for $k = 0, 1, 2$ in the following way$\colon$
	\begin{itemize}
		\item $(\scr{L}^{\mcal{I}}_1,b^\mcal{I}_{1}) \cong (\scr{L}^{\mcal{I}}_0,b^\mcal{I}_{0})$ if and only if
		\begin{equation}\label{equ_coordinatechange10gr2n}
			\begin{cases}
				x_{i,1} = u_{i}v_{i}-1, \, y_{i,1} = u_{i} &\quad \mbox{ for $(i,i+1) \in \mcal{I}$} \\
				z_{i,j,1} = z_{i,j,0} &\quad \mbox{otherwise.}
			\end{cases}
		\end{equation} 
		\item $(\scr{L}^{\mcal{I}}_0,b^\mcal{I}_{0}) \cong (\scr{L}^{\mcal{I}}_2,b^\mcal{I}_{2})$ if and only if
		\begin{equation}\label{equ_coordinatechange20gr2n}	
			\begin{cases}	
				x_{i,2} = u_{i}v_{i}-1, \, y_{i,2} = v_{i}^{-1} &\quad \mbox{ for $(i,i+1) \in \mcal{I}$} \\
				z_{i,j,2} = z_{i,j,0} &\quad \mbox{ otherwise.}
			\end{cases}
		\end{equation}
		\item $(\scr{L}^{\mcal{I}}_1,b^\mcal{I}_{1}) \cong (\scr{L}^{\mcal{I}}_2, b^\mcal{I}_{2})$ if and only if
		\begin{equation}\label{equ_coordinatechange12gr2n}
			\begin{cases}	
				x_{i,2} = x_{i,1}, \, y_{i,1} = y_{i,2} (1 + x_{i,2}) &\quad \mbox{ for $(i,i+1) \in \mcal{I}$} \\
				z_{i,j,2} = z_{i,j,1} &\quad \mbox{ otherwise.}
			\end{cases}		
		\end{equation}
	\end{itemize}
\end{lemma}

\begin{proof}
For a local calculation, we may pass Lagrangians to $\mcal{Y}_\mcal{I}$ equipped with the standard complex structure.
They are $\Z$-graded with respect to the holomorphic volume form~\eqref{equ_Zgradedvolmloca}. 
On each local $\mathrm{Gr}(2, \C^4)$ summand, there is a unique degree zero element of $\mathrm{CF}(\scr{L}_i, \scr{L}_j)$ in Section~\ref{sec_consimmergr24}.
We claim that the product of these degree zero elements is a quasi-isomorphism if the deformation spaces are glued by the relations~\eqref{equ_coordinatechange10gr2n},~\eqref{equ_coordinatechange20gr2n}, and~\eqref{equ_coordinatechange12gr2n}. 
The maximum principle enables us to analyze holomorphic curves on each piece. 
Then Lemma~\ref{cor:q-iso-hat} tells us how to identify deformation spaces in order to obtain quasi-isomorphisms in the summand.
\end{proof}

The holonomy variables $z_{i,j,k}$ for $k= 0,1,2$  are glued identically, that is, $z_{i,j,0} = z_{i,j,1} = z_{i,j,2}$. By abuse of notation, they will be denoted by $z_{i,j}$ from now on. 

Now, we need to find an isomorphism of the two charts $\scr{U}_2^{\mcal{I}}$ and $\scr{U}_2^{\mcal{I}^\prime}$ associated with the two Lagrangian tori $\scr{L}_2^\mcal{I}$ and $\scr{L}_2^{\mcal{I}^\prime}$ for different $\mcal{I}$ and $\mcal{I}^\prime$.  They are all Hamiltonian isotopic to the same monotone $\mathbb{T}^{2n-4}$-orbit, which limits to the monotone toric fiber under toric degeneration (See Section~\ref{localmirrorsecti}).

\begin{proposition}\label{prop_compii1}
For $\mcal{I}$ and $\mcal{I}^\prime$ in $\scr{I}_n$, the deformation spaces $\scr{U}_2^\mcal{I}$ and $\scr{U}_2^{\mcal{I}^\prime}$ are related by a quasi-isomorphism between $(\scr{L}^{\mcal{I}}_2, b^\mcal{I}_{2})$ and $(\scr{L}^{\mcal{I}^\prime}_2, b^{\mcal{I}^\prime}_{2})$ determined by
\begin{equation}\label{equ_x1yi}
x_{i,2} = \frac{z_{1,i+1}}{z_{1,i}} \cdot \frac{z_{2,i}}{z_{2,i+1}}, \,\,\,\, y_{i,2} = \frac{z_{2,i+1}}{z_{2,i}}.
\end{equation}
\end{proposition}

\begin{proof}
We focus on the algebraic torus $(\C^\times)^{2n-4}$ in~\eqref{equ_clusterchartempt}.
For any fixed $\mcal{I}$, let $\mcal{X}_\varepsilon$ and $\mcal{Y}_{\mcal{I}, \varepsilon}$ be toric degenerations of $\mathrm{Gr}(2,\C^n)$ and $\mathrm{Gr}(2,\C^4)^{|\mcal{I}|} \times \mathrm{Gr}(1,\C^2)^{2n-4 - 4 \cdot | \mcal{I} |}$ explained in Section~\ref{Sec_toricdegg} respectively. 
At each $\varepsilon \in \C$, the expression~\eqref{equ_cpxmapaa} defines a map $\Upsilon_{\textup{cpx},\varepsilon}$.
We have a symplectomorphism $\Xi_\varepsilon$ (resp. $\Xi_{\mcal{I}, \varepsilon}$) from the algebraic torus in $\mcal{X}_\varepsilon$ (resp. $\mcal{Y}_{\mcal{I}, \varepsilon}$) to that in $\mcal{X}_0$ (resp. $\mcal{Y}_{\mcal{I}, 0}$).
We then have the following commutative diagram of $\mathbb{T}^{2n-4}$-equivariant maps (on dense subsets)$\colon$
\begin{equation}\label{equation_toric_degeneration_diagram}
	\xymatrix{
		  (\C^\times)^{2n-4} \, (\textup{in } \mcal{X}_\varepsilon)  \ar[d]_{\Xi_\varepsilon} \ar[rr]^{ \Upsilon_{\textup{cpx},\varepsilon}}
                              & & (\C^\times)^{2n-4}\, (\textup{in } \mcal{Y}_{\mcal{I},\varepsilon})        \ar[d]^{\Xi_{\mcal{I}, \varepsilon}} \\
 (\C^\times)^{2n-4} \, (\textup{in } \mcal{X}_0)  \ar[rr]^{ \Upsilon_{\textup{cpx},0}} & &(\C^\times)^{2n-4} \, (\textup{in } \mcal{Y}_{\mcal{I},0})  }
\end{equation}

For each $\mcal{I}$, recall that the constructed Lagrangian $\scr{L}^\mcal{I}_2$ is  Hamiltonian isotopic to a $\mathbb{T}^{2n-4}$-orbit in $\mathcal{Y}_{\mcal{I},\epsilon}$.
Since $\Upsilon_{\textup{cpx}, \varepsilon}$ sends orbits to orbits, the $\mathbb{T}^{2n-4}$-orbit in $\mathcal{Y}_{\mcal{I},\epsilon}$ corresponds to a monotone $\mathbb{T}^{2n-4}$-orbit in $\mcal{X}_\varepsilon$. It degenerates into a monotone toric fiber in $\mcal{X}_0$.
Passing them into $\mcal{X}_0$ via $\Xi_\varepsilon$, we may compare $\scr{L}^\mcal{I}_2$ and $\scr{L}^{\mcal{I}^\prime}_2$ via the monotone toric fiber therein. Then~\eqref{equ_relationbetweenyijandxyz} leads to~\eqref{equ_x1yi}. 
\end{proof}

When computing the local mirror of $\scr{L}_2^\mcal{I}$, we encounter one issue here, caused by the fact that we are using the product complex structure on $\mcal{Y}_\mcal{I}$ for a local computation for quasi-isomorphisms. 
So, $\mathrm{Gr}(2,\C^n)$ is equipped with an almost complex structure that comes from the pull-backed complex structure on $\mcal{Y}_\mcal{I}$ via $\Upsilon_\textup{sym}$. 
Meanwhile, in order to compute the disk potential function, a specific complex structure associated with the toric degeneration in~\eqref{equ_pluckker} is necessary. Nevertheless, the potential function is invariant under a generic choice of almost complex structures since $\scr{L}^\mcal{I}_2$ is a monotone Lagrangian submanifold. See \cite{EP97} for the invariance of counting invariants on monotone Lagrangians. 
Consequently, Theorem~\ref{theorem_NNUmain} leads to a local mirror of $\scr{L}^\mcal{I}_2$. 

\begin{lemma} The LG mirror from $\scr{L}^\mcal{I}_2$ consists of  
\begin{equation}\label{equLGmm}
\left( \scr{U}_2^\mcal{I} := (\Lambda_{\mathrm{U}})^{2n-4}, 
W_{\scr{L}^\mcal{I}_2}(\textbf{\textup{z}}) = z_{2,1} +  \frac{T^n}{z_{1,n-2}} + \sum_{j=1}^{n-3}
\left( \frac{z_{1,j+1}}{z_{1,j}} +  \frac{z_{2,j+1}}{z_{2,j}} \right) + \sum_{j=1}^{n-2} \frac{z_{1,j}}{z_{2,j}} 
 \right)
 \end{equation}
 for any $\mcal{I} \in \scr{I}_n$.
\end{lemma}

\begin{remark}
Precisely speaking, in order to apply the quasi-isomorphisms~\eqref{equ_coordinatechange10gr2n},~\eqref{equ_coordinatechange20gr2n},~\eqref{equ_coordinatechange12gr2n}, the adjustment of the valuation is necessary. It is because the monotone torus fiber is taken as a reference fiber for our local calculation.
First, apply the coordinate changes $z_{1j} \mapsto T^{(1+j)} z_{1j}$ and  $z_{2j} \mapsto T^{j} z_{2j}$ to~\eqref{equLGmm}. Then apply the quasi-isomorphisms to get the other $W_{\scr{L}^\mcal{I}_\bullet}$'s. Finally, we return back by applying $z_{1j} \mapsto T^{-(1+j)} z_{1j}$ and  $z_{2j} \mapsto T^{-j} z_{2j}$.
\end{remark}

The gluing relations are computed in the local charts $\mathrm{Gr}(2, \C^n) \backslash \scr{D}^\mcal{I}$ of the anti-canonical divisor complement of $\mathrm{Gr}(2,\C^n)$.  They remain the same in $\mathrm{Gr}(2,\C^n)$ since holomorphic disks emanated from the anti-canonical divisor have higher Maslov index which do not contribute to the relations.

\begin{lemma}
	We still have the quasi-isomorphisms in $\mathrm{Gr}(2,\C^n)$ given by ~\eqref{equ_coordinatechange10gr2n},~\eqref{equ_coordinatechange20gr2n},~\eqref{equ_coordinatechange12gr2n}, and~\eqref{equ_x1yi} without modification. 
\end{lemma}

\begin{proof}
	Recall that all the Lagrangians we use are graded under the volume form ~\eqref{equation_merovolgr2n}.  The disks contributing to the quasi-isomorphisms (which have degree zero) must have Chern-Weil Maslov index zero.  If a pair of Lagrangians $L_1$ and $L_2$ constructed above bounds an additional holomorphic strip after compactifying a local chart to $\mathrm{Gr}(2,\C^n)$, the strip lies in the class $\beta+\alpha$ where $\beta$ is a strip class in the chart, and $\alpha$ is a disk class bounded by $L_1$ or $L_2$ in the ambient $\mathrm{Gr}(2,\C^n)$. They have Maslov index zero.  But since $L_i$ is monotone, any Maslov-zero holomorphic disk class is simply zero.  Thus there is no additional strip, and the quasi-isomorphisms remain the same.
\end{proof}

Then by the compatibility between gluing data and Floer potentials (from the $A_\infty$-relation), one can explicitly derive  $W_{\scr{L}^\mcal{I}_k}$ on the other mirror charts from the above relations.

With the explicit description of local charts at hand, we now identify $(\scr{U}_2^\mcal{I}, W_{\scr{L}^\mcal{I}_2})$ and $(\scr{U}_2^{\mcal{I}^\prime}, W_{\scr{L}^{\mcal{I}^\prime}_2})$ by Proposition~\ref{prop_compii1}.
By using the relations~\eqref{equ_coordinatechange10gr2n},~\eqref{equ_coordinatechange20gr2n}, and~\eqref{equ_coordinatechange12gr2n}, local mirror charts together with superpotentials are glued to obtain a desired partially compactified mirror. 

To summarize, we have the following tree of quasi-isomorphisms.  First we have the monotone torus orbit $\mathbb{T}^{2n-4}$.  Second we have quasi-isomorphisms from the tori $\scr{L}_2^\mcal{I}$ to $\mathbb{T}^{2n-4}$.  Third we have the quasi-isomorphisms from the Lagrangian immersions $\scr{L}_0^\mcal{I}$ to $\scr{L}_2^\mcal{I}$.  Finally we have the quasi-isomorphisms from the tori $\scr{L}_1^\mcal{I}$ (on the other sides of the walls) to $\scr{L}_0^\mcal{I}$.  Note that the quasi-isomorphisms between $\scr{L}_1^\mcal{I}$ and $\scr{L}_1^\mcal{I'}$ for different $I$ and $I'$ are obtained by compositions of the quasi-isomorphisms in this tree.
It completes the proof of Theorem~\ref{theroem_diskpotentialimgr26}

\subsection{Identification with the Rietsch's mirrors} \label{sec:WRie=WL}

Finally, we discuss the relation between the glued mirror and the Rietsch's mirror. 

\begin{theorem}\label{theorem_RietschSYZgr2n}
The Rietsch's mirror $(\check{X}, W_{\textup{Rie}})$ of the Grassmannian $X = \mathrm{Gr}(2, \C^n)$ is isomorphic to the LG model glued from Maurer-Cartan deformation spaces of the Lagrangians $\{ \scr{L}_j^\mcal{I} ~\colon j = 0,1,2, \mcal{I} \in \scr{I}_n^\textup{max} \}$ via quasi-isomorphisms given in Lemma~\ref{lemma_localchartcomq}.
\end{theorem}

By~\eqref{equ_coordinateidentifi}, the LG mirror of the monotone GC torus can be embedded into the Rietsch's mirror by the following identification$\colon$ 
\begin{equation}\label{equ_ccoorrx}
x_{i,1} = x_{i,2} = u_i v_i - 1 = \frac{z_{1,i+1}}{z_{1,i}} \cdot \frac{z_{2,i}}{z_{2,i+1}}\mapsto \frac{p_{n-i-2,n-i-1} \cdot p_{n-i, n}}{p_{n-i-2,n} \cdot p_{n-i-1,n-i}},
\end{equation}
\begin{equation}\label{equ_ccoorry}
y_{i,1} = u_i  \mapsto \frac{p_{n-i-2, n-i}}{p_{n-i-1,n-i}}, \quad
y_{i,2} = \frac{z_{2,i+1}}{z_{2,i}}  = v_i^{-1}  \mapsto \frac{p_{n-i-2,n}}{p_{n-i-1,n}}
\end{equation}
together with~\eqref{equ_coordinateidentifi} define an embedding of $\scr{U}^\mcal{I}_0, \scr{U}^\mcal{I}_1$, and $\scr{U}^\mcal{I}_2$ into the glued mirror. 

\begin{figure}[h]
	\begin{center}
	\includegraphics[scale=0.75]{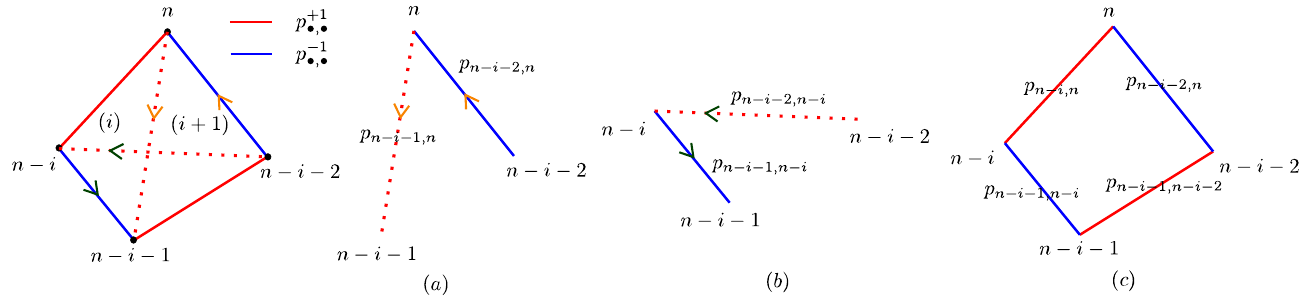}
	\vspace{-0.3cm}
	\caption{\label{fig_coordinatess} Combinatorial process reading the coordinate changes~\eqref{equ_ccoorrx} and~\eqref{equ_ccoorry}.}	
	\end{center}
\end{figure}

The coordinate changes~\eqref{equ_ccoorrx} and~\eqref{equ_ccoorry} can be read off from a triangulation of a planar $n$-gon as follows. 
Let $\mcal{I} \in \scr{I}_n$. 
Consider the subdivision of the $n$-gon associated with the immersed chart $\mcal{U}_\mcal{I}$ in Section~\ref{subimmersedre}.
For each $(i,i+1) \in \mcal{I}$, consider the $4$-gon.
There are two paths of length two starting at $n-i-2$ and ending at $n-i-1$$\colon$ one is via $n-i$ (Figure~\ref{fig_coordinatess} (a)) and the other is via $n$  (Figure~\ref{fig_coordinatess} (b)). 
Those two paths produces the formulae~\eqref{equ_ccoorry}, that is,
$$
v_i \mapsto \frac{p_{n-i-1,n}}{p_{n-i-2,n}}, \,\, u_i \mapsto \frac{p_{n-i-2,n-i}}{p_{n-i-1,n-i}}.
$$
Also, the loop around the $4$-gon (Figure~\ref{fig_coordinatess} (c)) produces~\eqref{equ_ccoorrx}. 

\vspace{-0.2cm}
\vspace{-0.1cm}

\begin{example}
Let recall Example~\ref{example_immersedpote} for $\mathrm{Gr}(2,\C^6)$.
Set $\mcal{I} = \{(2,3)\}$ and $\mcal{I}^\prime = \{(1,2), (3,4)\}$. 
By Theorem~\ref{theroem_diskpotentialimgr26}, 
$$
W_{\scr{L}_0^\mcal{I}} = W_{\mcal{I}}, \quad \mbox{and} \quad W_{\scr{L}_0^{\mcal{I}^\prime}} = W_{\mcal{I}^\prime}.
$$  
From Figure~\ref{fig_gr26tri}, the coordinate changes~\eqref{equ_ccoorrx} and~\eqref{equ_ccoorry} can be read off. 
Applying the coordinate change, we recover~\eqref{align_gri2} and~\eqref{align_gri} respectively.

\vspace{-0.2cm}
\begin{figure}[h]
	\begin{center}
	\includegraphics[scale=0.7]{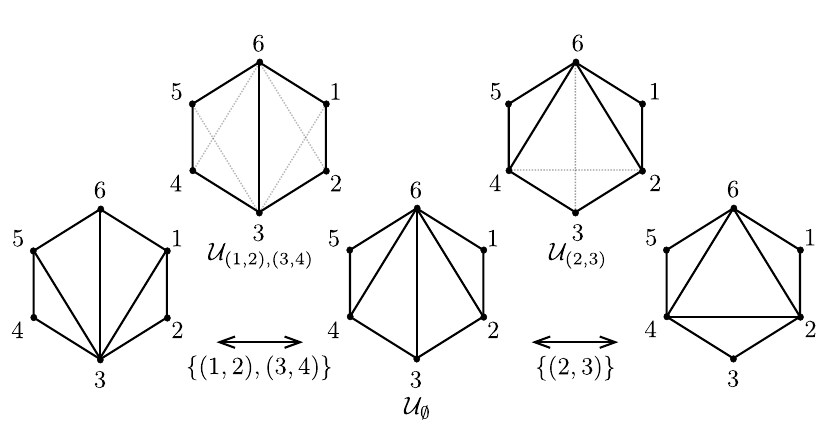}
	\vspace{-0.3cm}
	\caption{\label{fig_gr26tri} Immersed charts of $\check{X}$ for $\mathrm{Gr}(2,\C^6)$.}	
	\end{center}
\end{figure}
\vspace{-0.1cm}
\end{example}

Two LG models $(\scr{U}^\mcal{I}_0, W_{\scr{L}^\mcal{I}_0})$ and $(\mcal{U}, W_\mcal{I})$ in~\eqref{equ_immersedchar} are related by~\eqref{equ_ccoorrx},~\eqref{equ_ccoorry}, and~\eqref{equ_coordinateidentifi}. 
According to Lemma~\ref{lemma_enoughRIe}, the set $\{(\mcal{U}_\mcal{I}, W_\mcal{I}) ~\colon~ \mcal{I} \in \scr{I}^\textup{max}_n\}$ of immersed charts is enough to cover $\check{X}$. 
Thus, the glued mirror of the local mirrors $\{\scr{U}^\mcal{I}_k ~\colon~ \mcal{I} \in \scr{I}^\textup{max}_n, k = 0, 1, 2\}$ covers the Rietsch's mirror. It completes the proof of Theorem~\ref{theorem_RietschSYZgr2n}. 

We summarize as follows.
Take an $n$-gon, and label its vertices clockwisely by $1,\ldots,n$. Consider the triangulation by adding the edges $(i,n)$ for $i=2,\ldots,n-2$, see Figure~\ref{fig_referenceclus}. This corresponds to the monotone Lagrangian torus fiber.  Denote the triangles with vertices $i,i+1,n$ by $\Delta_{n-1-i}$.

We take away some of the edges so that there are only triangles and quadrilaterals in the decomposition.  Each triangle corresponds to $\mathbb{T}^2$, while each quadrilateral corresponds to $\cS_2 \times \mathbb{T}^2$ (where $\cS_2$ denotes the immersed sphere with only one self nodal point).  The polygonal decomposition corresponds to a Lagrangian immersion which is the product of these $\mathbb{T}^2$ and $\cS_2\times \mathbb{T}^2$.  The triangles are denoted as $\Delta_{n-1-i}$ as above, while the quadrilateral which is combined from the triangle $\Delta_{j}$ and $\Delta_{j+1}$ is denoted as $\Box_{j,j+1}$.

For each triangle $\Delta_j$, the corresponding $\mathbb{T}^2$-factor has two holonomy variables (parametrizing its flat connections) denoted by $z_{kj}$ for $k=1,2$.  For each quadrilateral $\Box_{j,j+1}$, the corresponding factor $\cS_2\times \mathbb{T}^2$ has immersed variables $u_j, v_j$ and holonomy variables $z_{1j}, z_{2,j+1}$.

The following extracts open Gromov-Witten invariants of immersed Lagrangians from the Rietsch mirror. Namely, if we write the Rietsch potential in geometric variables associated with immersed generators and holonomies, then each term appearing in the expansion is the count of disks with a topological type indicated by $u,v$ and $z$. The fractions in the expressions below can be understood as two-step paths in the polygons as in~\eqref{equ_ccoorrx} and~\eqref{equ_ccoorry}.

\begin{corollary} \label{thm:W^L}
	Let $W_\textup{Rie}$ be the Rietsch superpotential and $W_{\scr{L}_0}$ the disk potential for the monotone immersed Lagrangian $\scr{L}_0$ corresponding to the above polygonal decomposition.  $W_{\scr{L}_0}=T\cdot W_\textup{Rie}|_{q\mapsto 1}$ under the change of coordinates (where $p_{kl}$ is identified as $p_{lk}$ for the Pl\"ucker coordinates):
	\begin{align*}
	z_{1j} =& \frac{p_{n-j-1,n-j}}{p_{n-j,n}}; z_{2j} = \frac{p_{n-j-1,n}}{p_{n-1,n}}\\
	u_j =& \frac{p_{n-j-2,n-j}}{p_{n-j-1,n-j}}; v_j = \frac{p_{n-j-1,n}}{p_{n-j-2,n}}.
	\end{align*}
\end{corollary}

\section{Complete SYZ mirror of $\mathrm{OG}(1,\C^5)$.}
In this section, we construct a mirror of the orthogonal Grassmannian $\mathrm{OG}(1, \C^5)$ (see below for its definition) applying our gluing method developed in Section \ref{sec_wallcrossingimmersedLag}. While the method is not significantly different from the one in the previous sections, the choice of Lagrangian immersions is a bit tricky since this is a flag variety of type B.  
It also involves more work in the reduction procedure.
We will see that the resulting mirror agrees with the Lie-theoretical mirror in \cite{Rie, PRW}.

\subsection{Review of $\mathrm{OG}(1,\C^4)$}\label{sec_reviewog14}

The orthogonal Grassmannian $\mathrm{OG}(1, \C^4)$ is a partial flag manifold parametrizing the isotropic complex subspaces of one dimension in $\C^4$ equipped with a non-degenerate symmetric bilinear form.
As an algebraic variety, by choosing a suitable non-degenerate symmetric bilinear form, $\mathrm{OG}(1, \C^4)$ is a quadric hypersurface
$$
\mcal{Q}_{2} = \left\{ [Z_0: Z_1: Z_2: Z_3] \in \C\mathbb{P}^{3} ~\colon~ Z_0 Z_1 + Z_2^2 = Z_3^2 \right\}.
$$
It has the toric degeneration $\pi \colon \mcal{X} \to \C$ where $\mcal{X} = \bigcup_{\varepsilon \in \C}  \mcal{X}_\varepsilon$ and
$$
\mcal{X}_\varepsilon := \left\{ [Z_0: Z_1: Z_2: Z_3]  \in \C\mathbb{P}^{3} ~\colon~ Z_0 Z_1 + Z^2_2 = \varepsilon^2 Z_3^2 \right\}.
$$
A generic fiber $\mcal{X}_\varepsilon$ ($\varepsilon \neq 0$) is isomorphic to $\CP^1 \times \CP^1$. 

Passing $\mcal{X}_\varepsilon$ to $Z_3 \neq 0$, we generically have the smoothing of $A_1$-singularity 
$$
\mcal{X}^\textup{aff}_\varepsilon := \{ (X_0, X_1, X_2) \in \C^3 ~:~ X_0 X_1 + X_2^2 = \varepsilon^2 \},
$$
which degenerate into the $A_1$-singularity when $\epsilon=0$. 
The central fiber $\mcal{X}^\textup{aff}_0$ is a (singular) toric variety, whose moment polytope is defined by 
\begin{equation}\label{equ_momentpolytopea1singularity}
u_{1,2} - u_{1,1} \geq 0, u_{1,2} + u_{1,1} \geq 0.
\end{equation}

The smoothing comes with the projection $\Pi_2 \colon \mcal{Y}^\textup{aff}_\varepsilon \to \C$ to the $X_2$-component, which defines a conic fibration having two singular fibers located at $X_2 = \pm \varepsilon$. Moreover, the fibration admits a fiberwise $\mathbb{S}^1$-action given by
\begin{equation}\label{equ_s1orbit}
\mathbb{S}^1 \times \mcal{Y}^\textup{aff}_\varepsilon \to \mcal{Y}^\textup{aff}_\varepsilon, \quad (\theta, (X_0, X_1, X_2)) \mapsto (e^{-i \theta} X_0, e^{i \theta} X_1, X_2).
\end{equation}
By collecting the orbits satisfying $|X_0| = |X_1|$ over a simple closed curve $\gamma$ in the base $\C$, we obtain a Lagrangian torus (or a pinched torus) $\mathbb{T}_\gamma$.
Also, there is a Lagrangian sphere $\mathbb{S}^2$ over the line segment $[- \varepsilon, \varepsilon]$ obtained as a matching cycle. 

Take three simple closed curves$\colon$
\begin{itemize}
\item $\gamma_{0}$ passes through $\pm \varepsilon$,
\item $\gamma_1$ contains $\pm \varepsilon$ in its exterior, 
\item $\gamma_2$ contains $\pm \varepsilon$ in its interior,
\end{itemize}
as in Figure~\ref{Fig_base2} (a). They intersect cleanly with each other along disjoint circles. 

Observe that $\mathbb{T}_{\gamma_1}$ bounds a unique holomorphic disk of Maslov index two, a holomorphic section over the region bounded by $\gamma_1$. Thus the potential function of $\mathbb{T}_{\gamma_1}$ is a Laurent \emph{monomial}. Since a Lagrangian torus goes through two parallel walls under a Lagrangian isotopy from $\mathbb{T}_{\gamma_1}$ to  $\mathbb{T}_{\gamma_2}$, one can compute the potential function of $\mathbb{T}_{\gamma_2}$ by the wall-crossing. 

\begin{theorem}[\cite{auroux09, FOOO10, L14} ] \label{theorem_potentiaal2}
The potential function of $\mathbb{T}_{\gamma_2}$ in the smoothing of $A_1$-singularity is given by
$$
W_{\mathbb{T}_{\gamma_2}}(\textbf{\textup{y}}) = \frac{y_{1,2}}{y_{1,1}}(1+y_{1,1})^2.
$$
where $y_{1,j}$ is the exponential variable corresponding to the loop generated by $u_{1,j}$ in~\eqref{equ_momentpolytopea1singularity}.
\end{theorem}

\begin{remark}
Considering ${\mathbb{T}_{\gamma_2}}$ inside $\mathrm{OG}(1,\C^4) \simeq \CP^1 \times \CP^1$, one obtains its LG mirror.  
It only has two critical points, while the sum of the Betti number of $\mathrm{OG}(1,\C^4)$ is four. 
Thus, the LG mirror is incomplete.
In this case, however, we do \emph{not} need to take the immersed Lagrangian into account because the Lagrangian torus ${\mathbb{T}_{\gamma_1}}$ beyond the wall is Hamiltonian isotopic to a toric fiber. 
Furthermore, the LG mirror of ${\mathbb{T}_{\gamma_1}}$ is complete.
But, in $\mathrm{OG}(1,\C^5)$, the role of immersed Lagrangian is crucial. Namely, the glued mirror from Chekanov and Clifford tori is \emph{not} sufficient. 
$\mathrm{OG}(1,\C^5)$ will be discussed in the remaining sections. 
\end{remark}

\subsection{Construction of immersed Lagrangian}\label{sec_constructionofimmersedog15}

The orthogonal Grassmannian $\mathrm{OG}(1, \C^5)$ is a partial flag manifold parametrizing the isotropic complex subspaces of one dimension in $\C^5$ equipped with a non-degenerate symmetric bilinear form.
Again by taking a suitable non-degenerate symmetric bilinear form, $\mathrm{OG}(1, \C^5)$ is expressed as the quadric hypersurface given by
\begin{equation}\label{equ_q3var}
\mcal{Q}_{3} = \left\{ [Z_0: Z_1: Z_2: Z_3: Z_{4}] \in \C\mathbb{P}^{4} ~\colon~ Z_0 Z_1 + Z_2^2 = Z_3 Z_4 \right\}.
\end{equation}
Restricting to $Z_4 \neq 0$, we obtain  
$$
\mcal{Q}_{3}^\textup{aff} := \left\{ \left( X_0, X_1, X_2, X_3 \right) \in \C^4 ~\colon~
X_0 X_1 + X^2_2 = X_3 \right\}.
$$
where $X_i = Z_i / Z_4$ for $i = 0, 1, 2, 3$.  

The projection to the $X_3$-component given by
\begin{equation}\label{equ_X3projection}
\Pi_3 \colon \mcal{Q}_{3}^\textup{aff} \to \C^\times, \quad  \left( X_0, X_1, X_2, X_3 \right) \mapsto X_3
\end{equation}
has the smoothing of the $A_1$-singularity as a generic fiber, and the $A_1$-singularity itself sits at the origin $X_3 = 0$. Note that $\Pi_3$ has a \emph{non-trivial} monodromy around the origin
\begin{equation}\label{equ_monodromy}
X_i \leftrightarrow - X_i\,\,\, \mbox{for $i = 0, 1, 2$}.
\end{equation}

\begin{figure}[h]
	\begin{center}
		\includegraphics[scale=0.3]{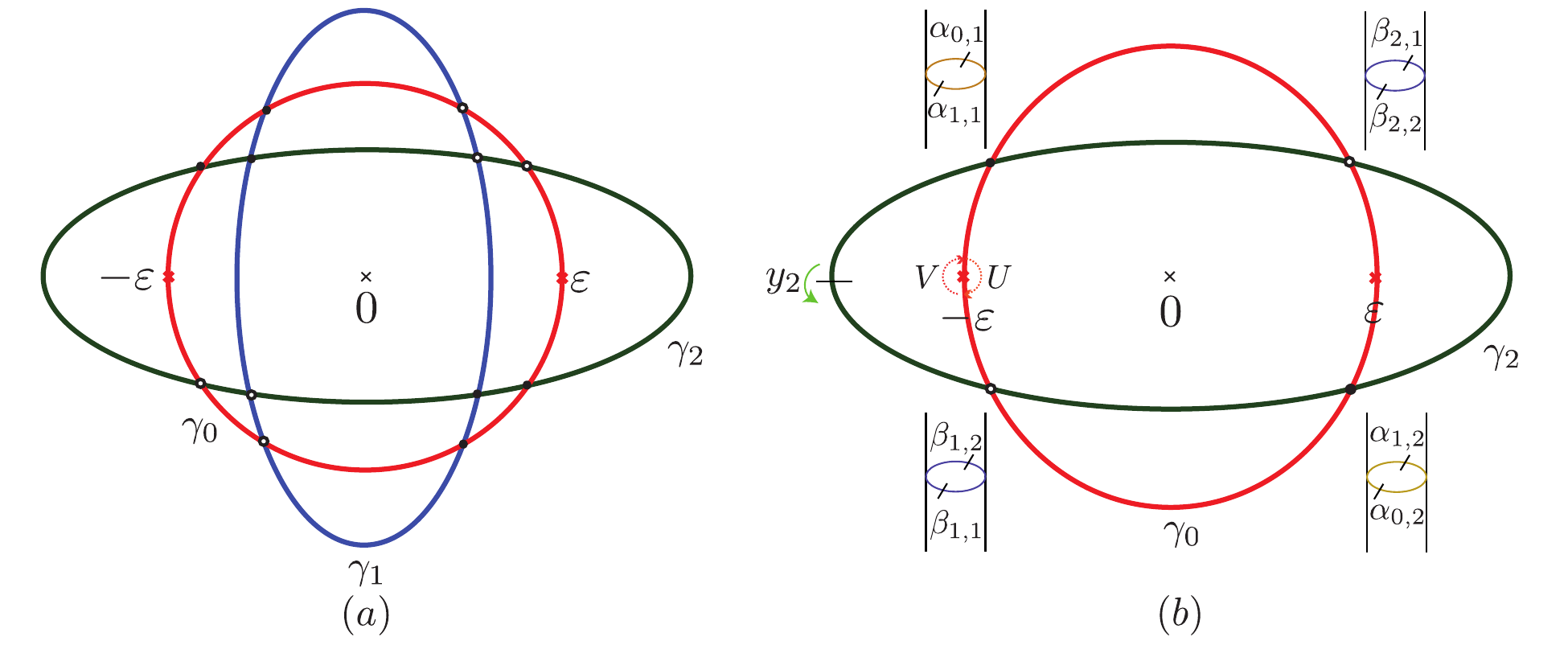} \vspace{-0.3cm}
		\caption{Simple closed curves on the base $\C_{X_2}$}
		\label{Fig_base2}
	\end{center}
\end{figure}

\begin{figure}[h]
	\begin{center}
		\includegraphics[scale=0.45]{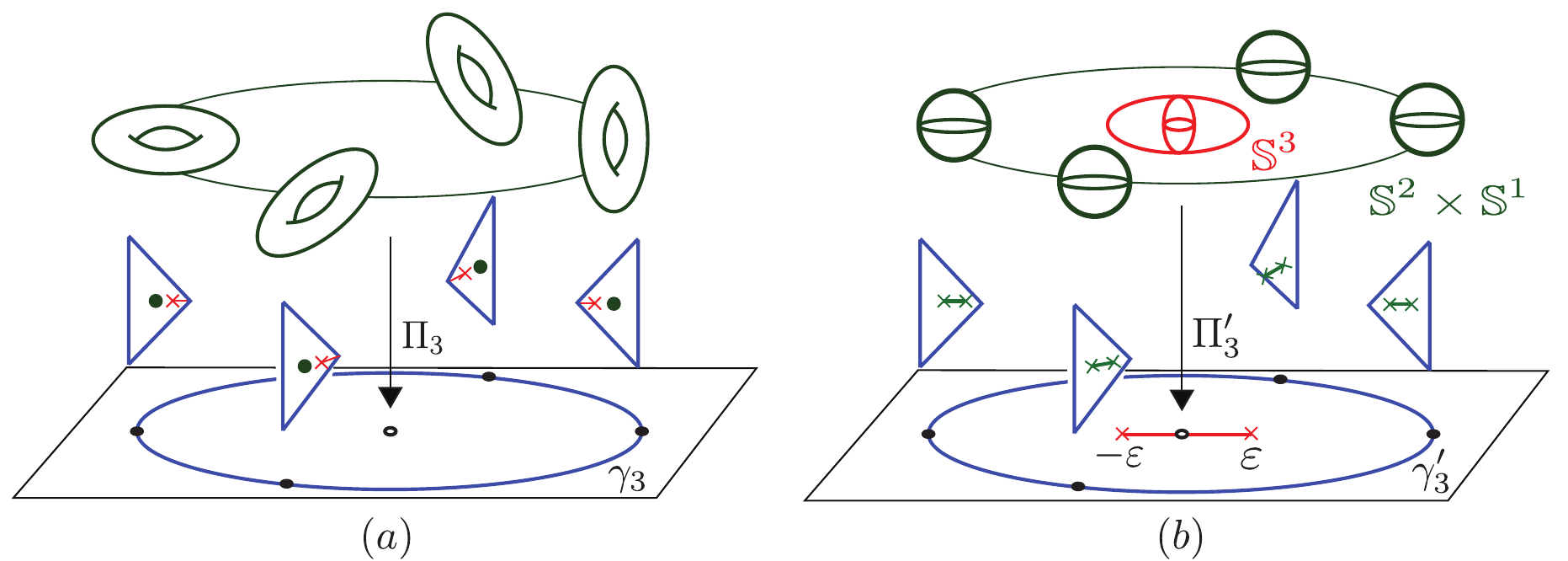}
			\vspace{-0.3cm}
		\caption{Two construction of Lagrangians}
		\label{Fig_Lag_Q3}
	\end{center}
	\vspace{-0.2cm}
\end{figure}

We then take simple closed curves $\gamma_0, \gamma_1,$ and $\gamma_2$ in $\C_{X_2}$ that are symmetric with respect to the origin as in Figure~\ref{Fig_base2} (a) and the $\mathbb{S}^1$-orbits generated by~\eqref{equ_s1orbit}. 
By parallel transporting the Lagrangian torus $\mathbb{T}_{\gamma_0}$, the Lagrangian torus $\mathbb{T}_{\gamma_1}$, and the immersed Lagrangian $\mathbb{T}_{\gamma_2}$ over a circle $\gamma_3$ on the base $\C_{X_3}$ as in Figure~\ref{Fig_Lag_Q3} (a), we obtain the following set of matching Lagrangians:
\begin{itemize}
\item an immersed Lagrangian $\scr{L}_0$,
\item a Lagrangian torus $\scr{L}_1$,
\item a Lagrangian torus $\scr{L}_2$.
\end{itemize}
Furthermore, by requiring the simple closed curves to bound the certain area, we may assume that the constructed Lagrangians are monotone.  

\vspace{-0.1cm}

\vspace{-0.3cm}

The quadric hypersurface has a Hamiltonian $\mathbb{T}^2$-symmetry.  
\begin{equation}\label{equ_t2syog15}
(\theta_1, \theta_2) \mapsto \left(e^{- \sqrt{-1} \theta_1} Z_0, e^{\sqrt{-1} \theta_1} Z_1, Z_2, e^{\sqrt{-1} \theta_2} Z_3, e^{-\sqrt{-1} \theta_2} Z_4  \right).
\end{equation}
The constructed Lagrangian $\scr{L}_j$ is invariant under the $\mathbb{T}^2$-symmetry if taking $\gamma_j$ and $\gamma_3$ as circles. 
Thus, $\scr{L}_j$ is Lagrangian isotopic to a $\mathbb{T}^2$-invariant Lagrangian. 

\begin{figure}[h]
	\begin{center}
			\vspace{-0.3cm}
		\includegraphics[scale=0.2]{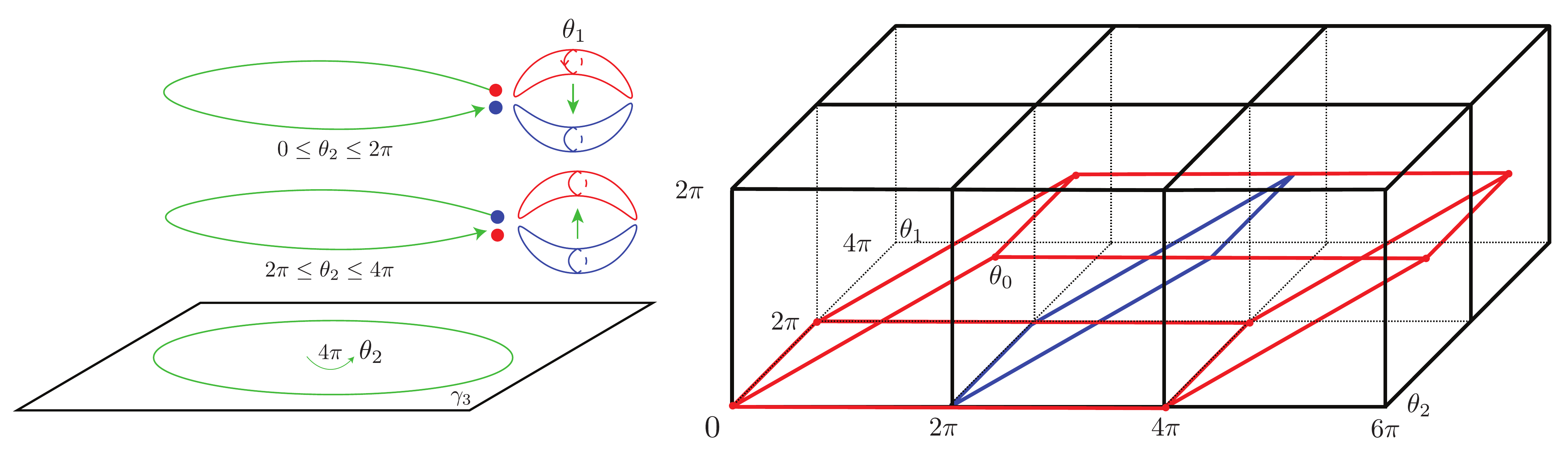}
			\vspace{-0.3cm}
		\caption{Fundamental domain of $\scr{L}_2$.}
		\label{Fig_fundomain}
	\end{center}
	\vspace{-0.3cm}
\end{figure}

\begin{remark}\label{remark_trivializationofshere}
At first glance, $\scr{L}_0$ produces a non-commutative mirror because $\mathbb{L}_0$ over $\gamma_0$ in Figure~\ref{Fig_base2} does, see \cite{CHLnc}. 
Yet the monodromy on the fibers of $\Pi_3$ along $\gamma_3$ swaps the immersed loci at $\pm \varepsilon$.
If trivializing the constructed Lagrangians along the orbit $\theta_2$,
$\scr{L}_0$ is indeed the product of immersed two sphere with a \emph{single} nodal self-intersection and circle as in Figure~\ref{Fig_fundomain}.
\end{remark}

\subsection{Local mirrors of the Lagrangians $\scr{L}_0, \scr{L}_1$, and $\scr{L}_2$}

To relate the Floer theoretical relation, we begin by choosing the meromorphic volume form 
$$
\Omega = \frac{dX_1 \wedge dX_2 \wedge dX_3}{X_1X_2X_3}, 
$$
which gives the $\Z$-grading on the Lagrangians $\scr{L}_\bullet$ on the complement of $\scr{D} = \{ Z_2 \cdot Z_3 \cdot Z_4 = 0\}$.

Equip the Lagrangians $\scr{L}_0, \scr{L}_1$, and $\scr{L}_2$ with the flat $\Lambda_\mathrm{U}$-connection $\nabla^{x_j, z_j}$ where $x_j$ and $z_j \in \Lambda_\mathrm{U}$ are the holonomy along the $\theta_1$ and $\theta_2$-orbits for $j = 0, 1, 2$ respectively. 
Let $\theta_0$ be a loop whose arguments of $X_i$ for $i = 0,1,2,3$ change from $0$ to $2\pi$.  
Also, adorn the torus $\scr{L}_2$ with the flat $\Lambda_\mathrm{U}$-connection $\nabla^{y_2}$ where $y_2$ is the holonomy along the $\theta_0$-orbit.

Regarding $\scr{L}_0$ as the product of immersed two sphere with a \emph{single} nodal self-intersection and circle, set $U$ and $V$ to be the degree one generators in $\scr{L}_0$, the products of immersed generators and the unit class of a $\theta_2$-orbit. 
As in the Section~\ref{subsec:supple}, the flat $\Lambda_\mathrm{U}$-connection $\nabla^{y_1}$ is chosen so that the standard wall-crossing formula in Theorem~\ref{thm:q-iso} is derived.
Geometrically, the flat connection can be realized as a choice of the gauge hypertori. The flat connection associated with the holonomy along the $\theta_j$-orbit is realized as a co-dimension one torus generated by the other $\theta_\bullet$'s with the co-orientation given by $\theta_j$.

We then have the formally deformed Lagrangians$\colon$
\begin{enumerate}
\item $(\scr{L}_i, b_{\scr{L}_i} := \nabla^{x_i, y_i, z_i})$ with the deformation space $\scr{U}_i \simeq (\Lambda_{\mathrm{U}})^3$  for $i = 1, 2$
\item $(\scr{L}_0, b_{\scr{L}_0} := (uU + vV, \nabla^{z_0}))$ with the deformation space $\scr{U}_0 \simeq (\Lambda_0 \times \Lambda_+) \cup (\Lambda_+ \times \Lambda_0) \times  \Lambda_{\mathrm{U}}$.
\end{enumerate}
By analyzing the cocycle condition on the local model (with which the Lagrangians become isomorphic), the relation~\eqref{equ_standardwallcross3},~\eqref{equ_standardwallcross2},~\eqref{equ_standardwallcross}, and $z_0 = z_1 = z_2$ can be derived. 

An explanation emphasizing on issues that do not appear in the Grassmannians of $A$-type in the previous sections is in order. 
Consider the projection from $\mcal{Q}_3 \backslash \scr{D}$ to $X_3$. 
By the maximum principle, there are two types of holomorphic strips bounded by $\scr{L}_i$ and $\scr{L}_j$ in the complement $\mcal{Q}_3 \backslash \scr{D}\colon$ the first type lies in the fiber of the projection and the second type does not. 
Because the origin is missing in the $X_3$-plane, the second type does \emph{not} appear. 
Recall that the curve rotating along $\gamma_3$ in the base $X_3$-plane twice counterclockwise lifts to a simple closed curve in $\scr{L}_i$ because of the monodromy~\eqref{equ_monodromy}.
It represents the cycle generated by a $\theta_2$-orbit, see Figure~\ref{Fig_fundomain}. 

To analyze the first type holomorphic strips, we choose Morse functions reflecting the monodromy~\eqref{equ_monodromy} on cleanly intersection loci as follows. 
For any pair of $\scr{L}_i$ and $\scr{L}_j$ ($i \neq j$), observe that they intersect cleanly at \emph{two} $\mathbb{T}^2$-orbits generated by $\theta_1$ and $\theta_2$. 
We shall take a product type Morse function on each $\mathbb{T}^2$-orbit satisfying the following. We choose a usual height function on the $\theta_1$-orbit and a $\Z/2$-invariant function on the $\theta_2$-orbit where the $\Z/2$-action is given by the monodromy.
Each $\mathbb{T}^2$-orbit is required to have two minimum points, four saddle points, and two maximum points. 
Also, assume that a minimum point of one torus and a maximum point of the other torus project to the same point $P$ in $X_3$. 

The critical points of the $\mathbb{T}^2$-orbits over the point $P$ are denoted as follows$\colon$ 
\begin{itemize}
\item two minimum points $\alpha^{\scr{L}_i, \scr{L}_j}_{0,1}, \alpha^{\scr{L}_i, \scr{L}_j}_{0,2}$, two saddle points $\alpha^{\scr{L}_i, \scr{L}_j}_{1,1}, \alpha^{\scr{L}_i, \scr{L}_j}_{1,2}$ in the first torus, 
\item two saddle points $\beta^{\scr{L}_i, \scr{L}_j}_{1,1}, \beta^{\scr{L}_i, \scr{L}_j}_{1,2}$, two maximum points $\beta^{\scr{L}_i, \scr{L}_j}_{2,1}, \beta^{\scr{L}_i, \scr{L}_j}_{2,2}$ in the second torus.
\end{itemize}
See Figure~\ref{Fig_base2} (b).

We now investigate the coordinate change on Maurer-Cartan deformation spaces so that $\alpha_0^{\scr{L}_i, \scr{L}_j} := \alpha_{0,1}^{\scr{L}_i, \scr{L}_j}  + \alpha_{0,2}^{\scr{L}_i, \scr{L}_j}$ becomes a quasi-isomorphism. 
Then the cocycle condition ${m}_1 (\alpha_0^{\scr{L}_i, \scr{L}_j}) = 0$ leads to the desired relation. For instance, the four strips shown in Figure~\ref{Fig_base2}. (b) contribute to $\langle {m}_1 (\alpha_0^{\scr{L}_i, \scr{L}_j}), \beta_1^{\scr{L}_i, \scr{L}_j} \rangle = 0$ where $\beta_1^{\scr{L}_i, \scr{L}_j} := \beta^{\scr{L}_i, \scr{L}_j}_{1,1} + \beta^{\scr{L}_i, \scr{L}_j}_{1,2}$.

Finally, by the local-to-global argument as in Theorem~\ref{theorem:q-iso-hat}, the derived relation gives us isomorphisms between the constructed Lagrangians in $\mathrm{OG}(1,\C^5)$.

\subsection{Computation of the potential functions}

To obtain the glued mirror, it remains to compute the potential function $W_{\scr{L}_2}$ of $\scr{L}_2$.
We will exploit the toric degeneration of~\eqref{equ_q3var} to compute $W_{\scr{L}_2}$. 
But, due to the monodromy~\eqref{equ_monodromy}, it is hard to relate the Lagrangian $\scr{L}_2$ with a toric fiber directly. 
Instead, regarding $\mathrm{OG}(1,\C^5)$ as a suspension of $\mathrm{OG}(1,\C^4)$, we construct a (monotone) Lagrangian torus ${\scr{L}}^\prime_2$, which degenerates into a toric fiber at $\mcal{X}_0$ and whose disk potential coincides with $W_{\scr{L}_2}$. 

Consider 
$$
\mcal{Q}_{3}^\prime := \left\{ [Z_0: Z_1: Z_2: Z_3: Z_{4}] \in \C\mathbb{P}^{4} ~\colon~ Z_0 Z_1 + Z_2^2 = Z_3^2 - Z_4^2 \right\}.
$$
We will employ the toric degeneration $\pi \colon \mcal{X} \to \C$ where $\mcal{X} = \bigcup_{\varepsilon \in \C}  \mcal{X}_\varepsilon$ and
\begin{equation}\label{equ_toricdeg2}
{\mcal{X}}^\prime_\varepsilon := \pi^{-1}(\varepsilon) = \left\{ [Z_0: Z_1: Z_2: Z_3: Z_{4}]  \in \C\mathbb{P}^{4} ~\colon~ Z_0 Z_1 + Z^2_2 = \frac{\varepsilon^2}{2} \left( Z_3^2 - \varepsilon^2 Z_4^2 \right) \ \right\}.
\end{equation}
Passing it to the affine chart $Z_4 \neq 0$, we obtain
$$
{\mcal{X}}^{\prime, \textup{aff}}_\varepsilon := \left\{ \left( X_0, X_1, X_2, X_3 \right) \in \C^4  ~\colon~
X_0 X_1 + X^2_2 = \frac{\varepsilon^2}{2} (X_3^2 - \varepsilon^2) \right\}.
$$

The projection to the $X_3$-component
\begin{equation}\label{equ_proprime3}
{\Pi}^\prime_3 \colon {\mcal{X}}^{\prime, \textup{aff}}_\varepsilon \to \C^\times, \quad  \left( X_0, X_1, X_2, X_3 \right) \mapsto X_3 
\end{equation}
defines a fibration which has a smoothing of $A_1$-singularity as a generic fiber. There are two singular fibers (both have $A_1$-singularity) at two points $X_3 = \pm \varepsilon$. Note that $\Pi_3$ has a trivial monodromy over the circle ${\gamma}^\prime_3$ enclosing $\pm \varepsilon$ at the interior. 

Using the fibration, we have the following three kinds of Lagrangians, see Figure~\ref{Fig_Lag_Q3} (b).
\begin{itemize}
\item Lagrangian $\mathbb{T}^3$ by parallel-transporting $\mathbb{L}_2$ in Section~\ref{sec_reviewog14} over ${\gamma}^\prime_3$
\item Lagrangian $\mathbb{S}^3$ by suspending the vanishing sphere over the line segment $[-\varepsilon, \varepsilon]$,
\item Lagrangian $\mathbb{S}^2 \times \mathbb{S}^1$ by parallel-transporting the vanishing sphere  over ${\gamma}^\prime_3$.
\end{itemize}

\begin{figure}[h]
	\begin{center}
			\vspace{-0.3cm}
		\includegraphics[scale=0.7]{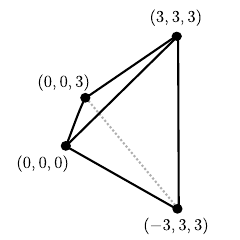}
			\vspace{-0.2cm}
		\caption{Gelfand-Cetlin polytope $\Delta^\textup{B}$}
		\label{Fig_GCpolytope}
	\end{center}
	\vspace{-0.2cm}
\end{figure}

\begin{remark}
The constructed Lagrangians can viewed as analogues of Gelfand-Cetlin Lagrangian fibers in $\mathrm{OG}(1,\C^5)$. 
The image of the Gelfand-Cetlin systems is the polytope $\Delta^\textup{B}$ in in Figure~\ref{Fig_GCpolytope}.
The Lagrangian $\mathbb{S}^3$-fiber is sitting over the origin and the Lagrangian $\mathbb{S}^2 \times \mathbb{S}^1$-fibers are sitting over the relative interior of the edge containing $(0,0,0)$ and $(0,0,3)$.
\end{remark}

The toric moment map of $\mcal{X}^\prime_0$ can be taken as
\begin{equation}\label{equ_momentmaps}
(u_{1,1}, u_{1,2}, u_{1,3}) := \left( \frac{- |Z_0|^2 + |Z_1|^2}{ \|Z\|^2}, \frac{|Z_0|^2 + |Z_1|^2 + |Z_2|^2}{\|Z\|^2}, \frac{|Z_0|^2 + |Z_1|^2 + |Z_2|^2 + |Z_3|^2}{\|Z\|^2} \right)
\end{equation}
where $\|Z\|^2 = |Z_0|^2 + |Z_1|^2 + |Z_2|^2 + |Z_3|^2 + |Z_4|^2$ at $\varepsilon = 0$, see \cite{NNU12} for instance.
The moment polytope is then given by 
\begin{equation}\label{equ_polytope}
|u_{1,1}| \leq u_{1,2} \leq u_{1,3} \leq 3,
\end{equation}
which will be denoted by $\Delta^\textup{B}$, see Figure~\ref{Fig_GCpolytope}.

\begin{theorem}\label{thm_potentialscrl2}
The potential function of ${\scr{L}}_2^\prime$ in $\mcal{Q}_3^\prime$ is
\begin{equation}\label{equ_potentialq3l2}
W_{{\scr{L}}^\prime_2}(\textbf{\textup{y}}) = \frac{1}{y_{1,3}} + \frac{y_{1,3}}{y_{1,2}} + \frac{y_{1,2}}{y_{1,1}} (1 + y_{1,1})^2.
\end{equation}
where $y_{1,1}, y_{1,2}$, and $y_{1,3}$ are the exponential variables corresponding to $u_{1,1}, u_{1,2}$, and $u_{1,3}$, respectively.
\end{theorem}

\begin{proof}
By applying the maximum principle to the holomorphic map ${\Pi}^\prime_3 \colon {\mcal{X}}^{\prime, \textup{aff}}_\varepsilon \to \C^\times$ in~\eqref{equ_proprime3}, we see that there are two kinds of holomorphic disks with Maslov index two: 
\begin{itemize}
\item[(i)]
the first one is fully contained in the fiber $\Pi_3^{\prime -1}(P)$ over a point $P$ of the base $\C^\times$, 
\item[(ii)]
 the second one is such that the image of the disk boundary under $\Pi_3$ is ${\gamma}^\prime_3$.
 \end{itemize}
For the first kind (i), the following three classes
$$
\beta_{1,2} - \beta_{1,1}, \, \beta_{1,2}, \, \beta_{1,2} + \beta_{1,1}
$$
can be represented by such holomorphic disks of Maslov index two by Theorem~\ref{theorem_potentiaal2}.
Each of the moduli spaces for the above classes is indeed an $\mathbb{S}^1$-family of holomorphic disks contained in the fiber over each point of ${\gamma_3^\prime}$.  
Moreover, the corresponding open Gromov-Witten invariants are respectively $n_{\beta_{1,2} - \beta_{1,1}} = 1, n_{\beta_{1,2}} =2$, and $n_{\beta_{1,2} + \beta_{1,1}} = 1$.

As $\varepsilon \to 0$, the spheres collapse and 
${\scr{L}}^\prime_2$ degenerates into a Lagrangian torus isotopic to a toric fiber as $\varepsilon \to 0$. One obtains two additional classes of the second kind (ii),
\begin{equation}\label{equ_classhittinginfinity}
H - \beta_{1,3}, \, \beta_{1,3} - \beta_{1,2}
\end{equation}
represented by Maslov two disks. They correspond to either the facet $3 - u_{1,3} = 0$ or the facet $u_{1,3} - u_{1,2} = 0$ in~\eqref{equ_polytope} respectively. The Fredholm regularity is persistent under a \emph{small} perturbation. Their open Gromov-Witten invariants are $n_{ H - \beta_{1,3}} = 1$ and $n_{\beta_{1,3} - \beta_{1,2}} = 1$. 

The argument so far confirms that $W_{\scr{L}^\prime_2}(\textbf{\textup{y}})$ contains at least all five terms in~\eqref{equ_potentialq3l2}. 
The following lemma completes the proof of Theorem~\ref{thm_potentialscrl2} by asserting that the five terms are indeed all in $W_{\scr{L}^\prime_2}(\textbf{\textup{y}})$.
\end{proof}

\begin{lemma}
The only possible classes in $\pi_2(\mcal{Q}^\prime_3, \scr{L}^\prime_2)$ that can be realized as a holomorphic disk of Maslov index two are $H - \beta_{1,3}, \beta_{1,3} - \beta_{1,2} , \beta_{1,2} - \beta_{1,1}, \beta_{1,2}, \beta_{1,2} + \beta_{1,1}$.
\end{lemma}

\begin{proof}
Consider $\widehat{\mcal{X}^\prime_0}$ be a crepant resolution of the toric variety ${\mcal{X}^\prime_0}$.
By abuse of notation, let $\scr{L} := \scr{L}^\prime_2$ in ${\mcal{X}^\prime_\varepsilon}$. 
Let $\widehat{\scr{L}}$ be a toric fiber of  $\widehat{\mcal{X}^\prime_0}$. 
Note that $H_2 (\widehat{\mcal{X}^\prime_0}, \widehat{\scr{L}}; \mathbb{Q}) \simeq \Q^5$, which is generated by the classes of holomorphic disks bounded by $\widehat{\scr{L}}$ and intersecting the toric anti-canonical divisor exactly once. For simplicity we abuse  
\begin{equation}\label{equ_fourclassesblowup}
H - \beta_{1,3},\, \beta_{1,3} - \beta_{1,2},\, \beta_{1,2} - \beta_{1,1},\,  \beta_{1,2} + \beta_{1,1}.
\end{equation}
to denote the classes corresponding to the facets in~\eqref{equ_polytope} in the polytope associated to $\widehat{\mcal{X}^\prime_0}$. Apart from these, we have one more class intersecting the exceptional divisor, which is denoted by $\widehat{\beta}_{1,2}$. Setting
$
\overline{H}_2(\widehat{\mcal{X}^\prime_0}, \widehat{\scr{L}}) := {H_2(\widehat{\mcal{X}^\prime_0}, \widehat{\scr{L}})}/{(\widehat{\beta}_{1,2} \sim \beta_{1,2})} \simeq \Z^4, 
$
the Maslov homomorphism factors through $\overline{H}_2(\widehat{\mcal{X}^\prime_0}, \widehat{\scr{L}})$ as follows$\colon$
\begin{equation}
	\xymatrix{
		  {H}_2(\widehat{\mcal{X}^\prime_0}, \widehat{\scr{L}}) \ar[dr] \ar[rr]^{ \mu}
                              & & \Z \\
  & \overline{H}_2(\widehat{\mcal{X}^\prime_0}, \widehat{\scr{L}})  \ar[ur] &}.
\end{equation}

For each $X_3 \in \C^\times_{X_3}$, recall that the fiber $\Pi^{^\prime -1}_3(X_3)$ is the smoothing of $A_1$-singularity generically. 
By collapsing the vanishing cycle fiberwise, we obtain a contraction map $\psi \colon ({\mcal{X}^\prime_t} , {\scr{L}}) \to ({\mcal{X}^\prime_0}, {\scr{L}}).$
We consider the resolution map (blow-up) $\phi \colon (\widehat{\mcal{X}^\prime_0} , \widehat{\scr{L}}) \to ({\mcal{X}^\prime_0}, {\scr{L}}).$
Then these two maps induce homomorphisms
$$
\begin{cases}
\psi_* \colon H_2({\mcal{X}^\prime_t} , {\scr{L}}; \Q) \to  H_2({\mcal{X}^\prime_0} , {\scr{L}}; \Q)\\
\phi_* \colon H_2(\widehat{\mcal{X}^\prime_0} , \widehat{\scr{L}}; \Q) \to  H_2({\mcal{X}^\prime_0} , {\scr{L}}; \Q).
\end{cases} 
$$
Observe that $\psi_*$ is an isomorphism since $H_2({\mcal{X}^\prime_t} , {\scr{L}}; \Q) \simeq \Q^4$ is generated by 
\begin{equation}\label{equ_fourclassesss}
H - \beta_{1,3},\, \beta_{1,3} - \beta_{1,2},\, \beta_{1,2} - \beta_{1,1},\,  \beta_{1,2} + \beta_{1,1},
\end{equation}
and the images of those classes under $\psi_*$ still generate $H_2({\mcal{X}^\prime_0} , {\scr{L}}; \Q)$.
The map $\phi_*$ is not an isomorphism, yet it induces an isomorphism 
$$
\overline{\phi}_* \colon \overline{H}_2(\widehat{\mcal{X}^\prime_0}, \widehat{\scr{L}}; \Q) \to  H_2({\mcal{X}^\prime_0} , \scr{L}; \Q).
$$
Furthermore, the composition $(\overline{\phi}_*)^{-1} \circ \psi^{\vphantom{-1}}_*$ preserves the Maslov indices because~\eqref{equ_fourclassesss} maps into~\eqref{equ_fourclassesblowup} respectively. 

Suppose  that there exists a holomorphic disk bounded by ${\scr{L}}$ with Maslov index two, representing a class other than $\beta_{1,1} + \beta_{1,2}, \beta_{1,2},$ and $\beta_{1,1} - \beta_{1,2}$. Let $a_0 \, H + a_1 \, \beta_{1,1} + a_2 \, \beta_{1,2} + a_3 \, \beta_{1,3}$
be the class represented by the disk. 
Note that the disk must be of the second kind (ii), which implies that $a_3$ must be non-zero. 

Observe that there exists a holomorphic disk (without any sphere bubbles) of Maslov index two in $((\overline{\phi}_*)^{-1} \circ \psi^{\vphantom{-1}}_*) (a_0 \, H + a_1 \, \beta_{1,1} + a_2 \, \beta_{1,2} + a_3 \, \beta_{1,3})$. Let $\varphi_t$ be a holomorphic disk in $a_0 \, H + a_1 \, \beta_{1,1} + a_2 \, \beta_{1,2} + a_3 \, \beta_{1,3}$. By considering sequence of holomorphic disks as $t \to 0$, we obtain a limit of holomorphic curve in $\mcal{X}_0^\prime$. Taking its strict transformation, we have a holomorphic curve in $((\overline{\phi}_*)^{-1} \circ \psi^{\vphantom{-1}}_*) (a_0 \, H + a_1 \, \beta_{1,1} + a_2 \, \beta_{1,2} + a_3 \, \beta_{1,3})$. By the classification of disks bounded by a toric fiber, the only possible Maslov index two holomorphic disks (with $a_3 \neq 0$) are those in~\eqref{equ_classhittinginfinity}. 
\end{proof}

\begin{figure}[h]
	\begin{center}
		\includegraphics[scale=0.5]{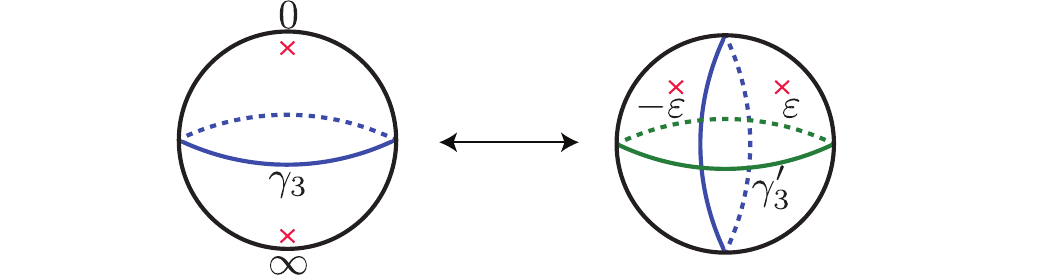}
			\vspace{-0.2cm}
		\caption{Comparison between two simply closed curves $\gamma_3$ and $\gamma_3^\prime$}
		\label{Fig_LagrangianonQ3}
	\end{center}
\end{figure}

With the potential function \eqref{equ_potentialq3l2} of $\scr{L}^\prime_2$ in hand, we next compute the local mirrors arising from $\scr{L}_0, \scr{L}_1,$ and $\scr{L}_2$. 
Two homogeneous coordinates for $\mcal{Q}_3$ and $\mcal{Q}^\prime_3$ are related by the following transformation
\begin{equation}\label{equ_mobius}
Z_0 \mapsto {Z}_0, \, Z_1 \mapsto {Z}_1,  \, Z_2 \mapsto {Z}_2, \,\,
\begin{bmatrix}
    {Z}_3       \\
    {Z}_4       
\end{bmatrix}
\mapsto
\frac{1}{\sqrt{2}}
\begin{bmatrix}
   1 & - \varepsilon      \\
   1 &  \varepsilon      \\
\end{bmatrix}
\cdot
\begin{bmatrix}
    {Z}_3       \\
    {Z}_4       
\end{bmatrix}.
\end{equation}
Under the M{\"o}bius transformation~\eqref{equ_mobius}, 
the base circle $\gamma_3$ in the sphere $[Z_3 \colon Z_4]$ in $\mcal{Q}_3$ maps into a circle in the sphere $[Z_3 \colon Z_4]$ in $\mcal{Q}^\prime_3$ enclosing only one $\pm \varepsilon$, see Figure~\ref{Fig_LagrangianonQ3}.

Note that $\scr{L}_2$ is related to $\scr{L}^\prime_2$ by a Lagrangian isotopy, which lies over an isotopy of simple closed curves on the base $\Pi^\prime_3 (\mcal{Q}_3^\prime)$, see the second sphere in Figure~\ref{Fig_LagrangianonQ3}. 
Any Lagrangian torus through the isotopy does \emph{not} bound any holomorphic disks of Maslov index zero because the simple closed curve $\gamma_2$ in $\C_{X_2}$ is fixed through the isotopy. Thus, the potential function of $\scr{L}_2$ is \emph{same} as that of $\scr{L}^\prime_2$. 

To translate the potential function $W_{{\scr{L}}^\prime_2}$ in~\eqref{equ_potentialq3l2} in terms of holonomy variables for $\scr{L}_2$, we analyze the arguments of variables of holonomy cycles. 
It leads to the following coordinate change. 
\begin{equation}\label{equ_coordinatechangel2l2tilde}
x_2 = y_{1,1}, \,\, y_2 = {y_{1,3}}, \,\, z_2 = \frac{y^2_{1,3}}{y_{1,2}}.
\end{equation}
Therefore
\begin{equation}\label{equ_potentialwscl2og15}
W_{\scr{L}_2} (x_2, y_2, z_2) = \frac{1}{y_2} + \frac{z_2}{y_2} + \frac{y_2^2}{x_2z_2} (x_2 + 1)^2.
\end{equation}
We then obtain the following theorem.

\begin{theorem}\label{theorem:q-iso-hatog}
In $\mathrm{OG}(1, \C^5)$, the potential functions of $\scr{L}_0, \scr{L}_1$, and $\scr{L}_2$ are 
\begin{itemize}
\item $\left( \scr{U}_0 = \left(\Lambda_0 \times \Lambda_+ \right) \cup \left( \Lambda_+ \times \Lambda_0\right) \times \Lambda_{\rm U}, W_{\scr{L}_0} \right)$ where
\begin{equation}\label{equn_potentil0og}
W_{\scr{L}_0} (u,v, z_0) = v + {vz_0} + \frac{u^2}{z_0(uv - 1)}.
\end{equation}
\item $\left(\scr{U}_1 = (\Lambda_{\rm U})^3, W_{\scr{L}_1} \right)$ where
\begin{equation}\label{equn_potentil1og}
W_{\scr{L}_1} (x_1, y_1, z_1) = \frac{1}{y_1} + \frac{x_1}{y_1} + \frac{z_1}{y_1} + \frac{x_1z_1}{y_1} + \frac{y_1^2}{x_1z_1}.
\end{equation}
\item $\left(\scr{U}_2 = (\Lambda_{\rm U})^3, W_{\scr{L}_2} \right)$ where
\begin{equation}\label{equn_potentil2og}
W_{\scr{L}_2} (x_2, y_2, z_2) =  \frac{1}{y_2} + \frac{z_2}{y_2} + \frac{y_2^2}{x_2z_2} (x_2 + 1)^2.
\end{equation}
\end{itemize}
In particular, by gluing them, a partially compactified mirror of $\mathrm{OG}(1, \C^5)$ is obtained.
\end{theorem}

\begin{corollary}
	The disk potential~\eqref{equn_potentil0og} of the immersed Lagrangian $\scr{L}_0$ in $\OG(1,5)$ has a critical point at $u=v=0, z_0=-1$.  Hence $\scr{L}_0$ is non-displaceable.
\end{corollary}

\subsection{Identification with the Rietsch's mirror of $\mathrm{OG}(1, \C^5)$}\label{sec_RietschmirrorQ3}

In \cite{Rie, PRW}, Pech, Rietsch and Williams constructed a LG mirror of the quadric hypersurface $\mcal{Q}_n$. 
The equivalent LG mirror for $\mcal{Q}_{2m-1}$ also appeared in Gorbounov-Smirnov~\cite{GoS15}.
It is a partial compactification of the LG mirror consisting of an algebraic torus together with a Laurent polynomial in Przyjalkowski \cite{Prz}. In this section, we provide an enumerative meaning of their LG models. 

The Pech-Rietsch-Williams' mirror is a LG model $(\check{X}, W_{\textup{Rie}})$ described as follows. Setting
$$
\scr{D}_0 := \{ p_0 = 0\},\, \scr{D}_1 := \{ p_1p_2 - p_0p_3 = 0 \},\, \scr{D}_2 := \{p_{3} = 0\},
$$
the mirror space $\check{X}$ is defined by the complement of the anti-canonical divisor $\scr{D}:= \scr{D}_0 + \scr{D}_1 + \scr{D}_2$ in $\CP^{3}$. 
The superpotential is given by
\begin{equation}\label{equ_supoerPRW}
W_{\textup{Rie}} = \frac{p_1}{p_0} + \frac{p_2^2}{p_1p_2 - p_0 p_3} + q \frac{p_1}{p_{3}}\colon \check{X} \to \C.
\end{equation}
Here is the list of its critical points.

\begin{lemma}\label{lemma_critcal}
The superpotential~\eqref{equ_supoerPRW} has four critical points$\colon$ 
\begin{enumerate}
\item Let $\xi := e^{i \frac{2\pi}{3}}$.
For $j = 0, 1, 2$, 
$$
p_{0} = 1,\, p_{1} = \sqrt[3]{4} \, \xi^j  q^{\frac{1}{3}},\, p_{2} = \sqrt[3]{2} \, \xi^{2j} q^{\frac{2}{3}},\, p_{3} = q
$$ 
The critical values are respectively $3 \sqrt[3]{4} \, \xi^j q^{\frac{1}{3}}$ for $j = 0, 1, 2$.

\item $(p_0,p_1,p_2,p_3) = (1,0,0,-q)$ whose
associated critical value is $0$. 
\end{enumerate}
\end{lemma}

We are now ready to state the main theorem of this section.

\begin{theorem}\label{theorem_PRWrecover}
The Rietsch's mirror $(\check{X}, W_{\textup{Rie}})$ of the quadric hypersurface $\mcal{Q}_3 \simeq \mathrm{OG}(1, \C^5)$ can be recovered by gluing the LG models arising from the immersed Lagrangian $\scr{L}_0$, the monotone Lagrangian tori $\scr{L}_1$ and $\scr{L}_2$ as deformation spaces of their Lagrangian Floer theory. 
\end{theorem}

\begin{proof}
After suitably adjusting valuations of variables for~\eqref{equn_potentil0og},~\eqref{equn_potentil2og}, and~\eqref{equn_potentil1og}, 
we can identify~\eqref{equ_supoerPRW} and~\eqref{equn_potentil0og} by the following relation 
$$
x_1 = x_2 = uv - 1 \mapsto \frac{p_1p_2 - p_0p_3}{p_0p_3}, \,\,  y_1 = u \mapsto \frac{p_2}{p_3}, \,\,  y_2 = v^{-1} \mapsto \frac{p_0}{p_1}, \,\, z_2 = z_1 = z_0 \mapsto \frac{p_0}{p_3}.
$$
\end{proof}

\section{Toward Homological Mirror Symmetry for $\mathrm{Gr}(2,\C^4)$ and $\mathrm{OG}(1,\C^5)$}\label{sec:HMS}

We have constructed the mirror LG models  for $\mathrm{Gr}(2,\C^4)$ and $\mathrm{OG}(1,\C^5)$ which are described as unions of local patches coming from Maurer-Cartan deformations of several reference Lagrangians. As the mirror potentials are Morse in these cases, their matrix factorization categories are generated by objects supported at the critical points. Since each critical point lies in one chart, it is enough to consider our local functors to study homological mirror symmetry in these cases.
 
More precisely, we will proceed as follows.
Recall that we have a LG model $(Y,W) $ which is obtained by gluing local patches $(U_i:=\mathcal{MC} (\scr{L}_i), W_i)$ where $Y= \cup U_i$ and each $W_i$ is Morse (and of course, $W_i |_{U_i \cap U_j} = W_j|_{U_i \cap U_j}$). For each critical point $y$ of $W$, there exists $U_{i(y)}$ that contains $y$. We fix it once and for all. We have a functor
$$\mathcal{F}_{i(y)} : \Fuk (X) \to \MF(U_{i(y)},W_{i(y)} ) \,\, (\hookrightarrow \MF(Y,W))$$
induced by Maurer-Cartan deformation of $\scr{L}_i$ (where $X= \mathrm{Gr}(2,\C^4)$ or $\mathrm{OG}(1,\C^5)$). Let $P_y$ be the image of $(\scr{L}_{i(y)}, b(y))$ where $b(y)$ is the weak bounding cochain of $\scr{L}_{i(y)}$ corresponding to $y \in U_i = \mathcal{MC} (\scr{L}_{i(y)})$. Note that we also have an $A_\infty$-algebra homomorphisms
$$ CF((\scr{L}_{i(y)}, b(y)), (\scr{L}_{i(y)}, b(y))) \to \hom (P_y, P_y).$$

Now consider the sub-Fukaya category $\mathcal{A}$ generated by the objects $(\scr{L}_{i(y)}, b(y))$ for all critical point $y$  of $W$. Since two objects $(\scr{L}_{i(y_1)}, b(y_1))$ and $(\scr{L}_{i(y_2)}, b(y_2))$ corresponding to different critical points $y_1$ and $y_2$ do not have a nontrivial morphism space, we obtain a well-defined $A_\infty$ functor
$$\mathcal{A} \to \MF(Y,W)$$
which sends $(\scr{L}_{i(y)}, b(y))$ to $P_y$ on the object level. 
In this section, we show that this functor establishes an equivalence both for $\mathrm{Gr}(2,\C^4)$ and $\mathrm{OG}(1,\C^5)$ (after being derived), which proves Theorem \ref{prop:HMSpart}.

\subsection{Mirror matrix factorizations for $\mathrm{OG}(1,\C^5)$} \label{subsec:maog15}
We first compute the mirror matrix factorizations of reference Lagrangians in $\mathrm{OG}(1,\C^5)$ together with (weak) bounding cochains that correspond to the critical points of $W$. By the definition of our local functor, they are simply Floer complexes with suitable boundary deformation. We only spell out the argument of the immersed Lagrangians $\scr{L}_0$, and for smooth tori $\scr{L}_i$ ($i=1,2$), we refer readers to \cite[Section 9]{CHLtoric}. (Indeed, the same argument for $\scr{L}_0$ works for $\scr{L}_1$ and $\scr{L}_2$ without much modification.)

We first set up the following notations for the standard generators of $CF((\scr{L}_0,b),(\scr{L}_0,\underline{b}))$ where $b$ varies over $\mathcal{MC} (\scr{L}_0)$ and $\underline{b}$ corresponds to the critical point $(u,v,z_0) = (0,0,-1)$. 
$$ \one:= (\one_{\BL_0}, \one_{\mathbb{S}^1}), \Theta_1:= (U, \one_{\mathbb{S}^1}), \Theta_2:=(V, \one_{\mathbb{S}^1}), \Theta_3:= (\one_{\BL_0}, \pt_{\mathbb{S}^1}), $$
$$ \Theta_1 \wedge \Theta_2:= (\pt_{\BL_0}, \one_{\mathbb{S}^1}), \Theta_2 \wedge \Theta_3 := (V, \pt_{\mathbb{S}^1}), \Theta_1 \wedge \Theta_3:= (U, \pt_{\mathbb{S}^1}),$$
$$ \Theta_1 \wedge \Theta_2 \wedge \Theta_3 := (\pt_{\BL_0}, \pt_{\mathbb{S}^1}).$$
(Recall that $\scr{L}_0$ is a product of an immersed 2-sphere $\BL_0$ with $\mathbb{S}^1$.) Here $\one_{\mathbb{S}^1}$ and $\pt_{\mathbb{S}^1}$ actually mean the critical points of the Morse function on $\mathbb{S}^1$-factor with the corresponding degrees, and the same applies to $\one_{\BL_0}$ and $\pt_{\BL_0}$.
In addition, we define wedge products of $\Theta_i$'s in other orders by imposing usual skew-commuting relations among $\Theta_i$'s, for e.g., $\Theta_2 \wedge \Theta_1:= - \Theta_1 \wedge \Theta_2$. Then $CF((\scr{L}_0,b),(\scr{L}_0,\underline{b}))$ can be identified with the exterior algebra generated by $\Theta_1, \Theta_2,\Theta_3$ as a vector space. Note that the differential $\delta:=m_1^{b,\underline{b}}$ on this complex can be decomposed as
$$ \delta = \delta_{+1} + \delta_{-1} + \delta_{-3}$$
where the sub-indices are the degrees of operators with respect to the natural degree from the exterior algebra. $\delta_{i}$ with $i\leq -5$ vanishes by degree reason.

We can identify $\delta_{+1} = u \Theta_1 \wedge (-). + v \Theta_2 \wedge (-) + (z_0 - (-1)) \Theta_3 \wedge (-)$  since the contributing holomorphic disks are those which contribute to the weak Maurer-Cartan equation (and hence have Maslov index 0). Here the signs are determined by the Koszul convention with respect to the product, which is similar to the sign rule for Morse flows appearing in \cite[Appendix A]{CHLtoric}. For instance, if a flow runs along the second factor keeping the first factor constant, there comes an additional sign coming from (the parity of) the degree of the first component of an input.
Therefore we have
$$ \delta (\one) = u \Theta_1 + v \Theta_2 + (z_0+1) \Theta_3$$
$$ \delta (\Theta_1) = -v \Theta_1 \wedge \Theta_2 + (z_0 +1) \Theta_3 \wedge \Theta_1 + f _1 \cdot \one $$
$$ \delta (\Theta_2) =  u \Theta_1 \wedge \Theta_2 - (z_0 +1) \Theta_2 \wedge \Theta_3 + f _2 \cdot \one $$
$$ \delta (\Theta_3) = v \Theta_2 \wedge \Theta_3 - u \Theta_3 \wedge \Theta_1 + f_3 \cdot \one $$
where $f_1,f_2,f_3$ are unknown, but should satisfy
\begin{equation}\label{eqn:wcone}
\delta^2 ( \one) = u f_1 + v f_2 + (z_0 +1) f_3 = W - W(\underline{b}).
\end{equation}
In fact, $f_i$ is contributed by the same set of the holomorphic polygons which also contribute the potential. 

\begin{figure}[h]
	\begin{center}
		\includegraphics[scale=0.6]{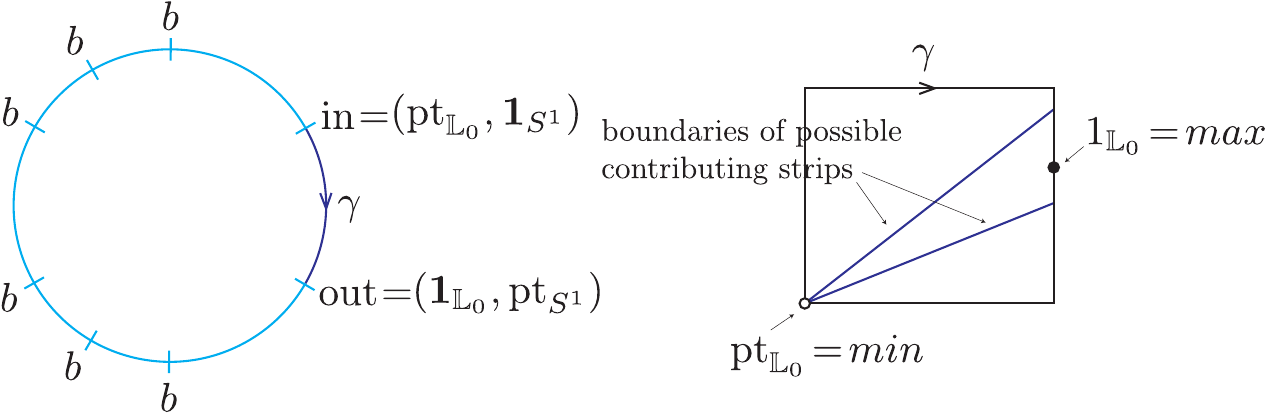}
			\vspace{-0.2cm}
		\caption{boundary shapes of the strips from $\Theta_1 \wedge \Theta_2$ to $\Theta_3$}
		\label{fig:mfnonex}
	\end{center}
\end{figure}

Let us proceed to the next degree. We have
$$ \delta (\Theta_1 \wedge \Theta_2 )= (z_0+1) \Theta_1 \wedge \Theta_2 \wedge \Theta_3 - g_2 \cdot \Theta_1 + g_1 \cdot \Theta_2 + \tilde{g} \cdot \Theta_3$$
for some $g_1, g_2, \tilde{g} $, where the first term on the right hand side is simply $\delta_{+1}$ applied  to $\Theta_1 \wedge \Theta_2$. 
We claim that $\tilde{g}$ is zero. Suppose to the contrary that there exists a strip from $(\pt_{\BL_0}, \one_{\mathbb{S}^1})$ to $(\one_{\BL_0}, \pt_{\mathbb{S}^1})$. Such a strip should come from a Maslov 2 disk with the configuration as in Figure \ref{fig:mfnonex}. Take a boundary path $\gamma$ of the strip running from the input to the output in clockwise direction. Notice that $\gamma$ does not have any corner since $u=v=0$ for $\underline{b}$. Therefore, one can trivialize $\BL_0$ along $\gamma$ to have a circle fibration over $\gamma$ as in the right diagram in Figure \ref{fig:mfnonex}. (These fiber circles are the ones appearing in conic fibers in our local picture.)

The moduli of Maslov two disks passing through $(\pt_{\BL_0},\one_{\mathbb{S}^1})$ is one dimensional, and this one dimensional deformation indeed comes from the (global) ${S}^1$-action given by $(\theta, [Z_i]) \mapsto [Z_0: Z_1: Z_2: e^{i \theta} Z_3, e^{-i \theta} Z_4]$, which rotates the ${S}^1$-factor in $\scr{L}_0 = \BL_0 \times {S}^1$. Moreover, the action preserves the input/output condition on the second factor $\mathbb{S}^1$ of $\scr{L}_0$ since there is a flexibility from flow lines from $\one_{\mathbb{S}^1}$ to $\pt_{\mathbb{S}^1}$. Therefore, all possible contributions to the coefficient of $(\one_{\BL_0}, \pt_{{S}^1})$ in $\delta(\pt_{\BL_0}, \one_{{S}^1})$ give discrete set of curves in this fibration after taking the boundary segment from the input to the output clockwisely. Note that these curves should exactly hit $\pt_{\BL_0}$ (minimum) and $\one_{\BL_0}$ in order to actually contribute. Obviously a generic Morse function avoid such a contribution, which proves the claim. (This is essentially the same argument as in the proof of \cite[Lemma 8.6]{CHLtoric}.)

Also, there is no strip from $\Theta_2 \wedge \Theta_3=(V,\pt_{\mathbb{S}^1})$ to $\Theta_1=(U,\one_{\mathbb{S}^1})$. A similar argument as above still works for this case, and in fact, for the case of $\Gr(2,\C^4)$ as well, showing that there is no (Maslov 2) polygon from $\Theta_{I}$ to $\Theta_{J}$ unless $J \subset I$, which we will not repeat afterwards. Here $I,J \subset \{1,2,3,4\}$ and $|J| = |I| -1$.

Since $(U,\one_{\mathbb{S}^1})$ is an output, a contributing polygon has two $V$-corners at two ends of the boundary  path $\gamma'$ joining the input and the output in clockwise direction. ($\gamma'$ is analogous to $\gamma$ in the previous case.) This time, we consider the $\mathbb{S}^1$-fibration over $\gamma'$ with fibers coming from the $\mathbb{S}^1$-factor in $\scr{L}_0 = \mathbb{L}_0 \times \mathbb{S}^1$. As before, boundary paths of all possible contributing strips give a finite set of curves on this trivialization, and hence $min(=\pt_{\mathbb{S}^1})$ and $max(=\one_{\mathbb{S}^1})$ of a generic Morse function can avoid all such curves. (Again, boundary paths should exactly hit these two points to actually contribute.) Note that this argument can not apply to a strip from $(V,\pt_{\mathbb{S}^1})$ to $(V,1)$ since two corners of a strip associated to input and output may appear in a constant disk component (possibly attached to a non-constant disk by a flow line).

Therefore $\delta(\Theta_i \wedge \Theta_j)$ does not involve $\Theta_k$ if $k$ is not either $i$ or $j$, and hence we have
$$\delta (\Theta_2 \wedge \Theta_3 )= u \Theta_1 \wedge \Theta_2 \wedge \Theta_3 -  h_3 \cdot \Theta_2 +  h_2 \cdot \Theta_3$$
$$\delta (\Theta_3 \wedge \Theta_1 )= v \Theta_1 \wedge \Theta_2 \wedge \Theta_3 -  k_1 \cdot \Theta_3 +  k_3 \cdot \Theta_1$$
for some functions $h_i$ and $k_i$. Observe that
\begin{equation*}
\begin{array}{lcl}
\delta^2 (\Theta_1) &=& f_1 ( u \Theta_1 + v \Theta_2 + (z_0+1) \Theta_3) -v \delta_{-1} (\Theta_1 \wedge \Theta_2) + (z_0+1)  \delta_{-1} (\Theta_3 \wedge \Theta_1) \\
&=& (u f_1 + v g_2 + (z_0 +1) k_3 ) \Theta_1 + v( f_1 - g_1) \Theta_2 + (z_0+1) (f_1 - k_1) \Theta_3
\end{array}
\end{equation*}
Since $\delta^2 (\Theta_1) = (W- W(\underline{b})) \cdot \Theta_1$ from the $A_\infty$-relation, we see that $g_1 = k_1 = f_1$. Likewise, we have gives $g_2 = h_2 = f_2$ and $h_3= k_3 = f_3$.

We next compute $\delta(\Theta_1 \wedge \Theta_2 \wedge \Theta_3)$. We first show $\delta_{-3}$ vanishes. Let $\delta_{-3} (\Theta_1 \wedge \Theta_2 \wedge \Theta_3) = \xi \cdot \one$ for some function $\xi$.
By observing the coefficient of $\one$ in $\delta^2 (\Theta_1 \wedge \Theta_2)$ (which should be zero), one obtain $-g_2 f_1 + g_1 f_2 +(z_0 +1) \xi =0$, which implies $\xi$ is zero since $g_i =f_i$ for $i=1,2$.
Let us write 
$$\delta (\Theta_1 \wedge \Theta_2 \wedge \Theta_3 ) = l_1 \Theta_2 \wedge \Theta_3 + l_2 \Theta_3 \wedge \Theta_1 + l_3 \Theta_1 \wedge \Theta_2.$$
for some $l_1,l_2,l_3$. From $\delta^2(\Theta_1 \wedge \Theta_2) = (W- W(\underline{b}) ) \Theta_1 \wedge \Theta_2$, we have 
$$ u f_1 + v f_2 + (z_0 +1) l_3 = W - W(\underline{b}) = u f_1 + v f_2 + (z_0 +1) f_3,$$
and hence $l_3 = f_3$. Likewise $l_1 = f_1$ and $l_2 = f_2$. 

In summary, one can identify $\delta$ in terms of $\Theta_i$'s as
$$\delta =u \Theta_1 \wedge (-). + v \Theta_2 \wedge (-) + (z_0 - (-1)) \Theta_3 \wedge (-) + \sum_{i=1}^3 f_i \frac{\partial}{\partial \Theta_i}$$
which is precisely a wedge-contraction type matrix factorization known to generate the component $\MF_{W(\underline{b})} (W)$ of $\MF(W)$ (consisting of matrix factorizations of $W-W(\underline{b})$) by the work of Dyckerhoff \cite{Dy}. Also the morphism level functor
$$ \mathcal{F}^{\scr{L}_0} : CF( (\scr{L}_0, \underline{b}), (\scr{L}_0, \underline{b}) ) \to \hom (P_{\underline{b}}, P_{\underline{b}})  $$
is injective by \cite[Theorem 4.9]{CHLnc}, and hence an isomorphism since the cohomology groups of both sides have dimension $8$. The same argument appears in \cite[Section 7]{CHLabc} to prove homological mirror symmetry for $\mathbb{P}^1_{a,b,c}$.

Finally, for each critical value $\lambda$ of $W$, we know that the corresponding component in $QH(\mathrm{OG}(1,\C^5))$ is 1-dimensional (i.e. the $\lambda$-eigenspace in $QH(\mathrm{OG}(1,\C^5))$ with respect to the operator $c_1 (\mathrm{OG}(1,\C^5)) \wedge$). Having closed string mirror symmetry, one can see this from the mirror LG model in which there is only one critical point over each critical value of $W$ (see Lemma \ref{lemma_critcal}). By \cite[Corollary 2.19]{Sheridan-Fano} tells us that the Fukaya category is generated by our reference Lagrangians together with weak bounding cochains that correspond to critical points of $W$. This completes the proof of the statement in Theorem \ref{prop:HMSpart} concerning $\mathrm{OG}(1,\C^5)$.

\subsection{Mirror symmetry for $\mathrm{Gr}(2,\C^4)$}
We next compute the mirror matrix factorization for $\mathrm{Gr}(2,\C^4)$.
As in the previous case, we will only exhibit the argument for immersed reference Lagrangians $\scr{L}_0$, which is the product of the immersed 2-sphere $\BL_0$ with the 2-torus $\mathbb{T}^2$. Let $\underline{b}$ be a critical point of $W_{\scr{L}_0}$, which is either of the following two
$$
\underline{b}_j= \{ u = 0 , v = 0, z_0 = - 1, w_0 = \sqrt{-1} \, \xi^{2j} \} \quad \mbox{for } j = 1, 2.
$$
($(\scr{L}_0, \underline{b}_1)$ and $(\scr{L}_0,\underline{b}_2)$ have trivial Floer cohomology between them. Indeed, one can easily check that the unit class, which is $\one$ below, is in the image of the Floer differential.) For simplicity, let us write $\underline{b} = (0,0,\underline{z_0}, \underline{w_0})$ where $\underline{z_0} = -1$ and $\underline{w}_0 = \sqrt{-1} \xi^{2j}$ (for $j=1$ or $2$).

As before, we set the generator of $CF( (\scr{L}_0,b),(\scr{L}_0,\underline{b}) )$ as follows.
$$ \one:= (\one_{\BL_0}, \one_{\mathbb{T}^2}), \,\,  \Theta_1:= (U, \one_{\mathbb{T}^2}), \,\,\Theta_2:=(V, \one_{\mathbb{T}^2}), \,\, \Theta_3:= (\one_{\BL_0}, \alpha),\,\, \Theta_4:= (\one_{\BL_0}, \beta)$$
$$ \Theta_1 \wedge \Theta_2:= (\pt_{\BL_0}, \one_{\mathbb{T}^2}), \,\, \Theta_1 \wedge \Theta_3 := (U, \alpha), \,\, \Theta_1 \wedge \Theta_4:= (U, \beta), $$
$$\Theta_2 \wedge \Theta_3 := (V , \alpha ),\,\, \Theta_2 \wedge \Theta_4 := (V ,\beta ), \,\,\Theta_3 \wedge \Theta_4 := ( \one_{\BL_0} , \pt_{\mathbb{T}^2} )$$
$$ \Theta_1 \wedge \Theta_2 \wedge \Theta_3 := (\pt_{\BL_0}, \alpha), \,\, \Theta_1 \wedge \Theta_2 \wedge \Theta_4:=(\pt_{\BL_0}, \beta),$$
 $$ \Theta_1 \wedge \Theta_3 \wedge \Theta_4:= (U, \pt_{\mathbb{T}^2}), \,\, \Theta_2 \wedge \Theta_3 \wedge \Theta_4 := (V, \pt_{\mathbb{T}^2})$$
$$\Theta_1 \wedge \Theta_2 \wedge \Theta_3 \wedge \Theta_4:= (\pt_{\BL_0}, \pt_{\mathbb{T}^2} )$$
where $\alpha$ and $\beta$ are (dual to) the cycles in $\mathbb{T}^2$ around which holonomies are $z_0$ and $w_0$ for $(\scr{L}_0,b)$.
Again, we can decompose $\delta:=m_1^{b,\underline{b}}$ into 
$$\delta = \delta_{+1} + \delta_{-1} + \delta_{-3}.$$
We will see later that $\delta_{-3}$ vanishes as in the case of $\OG(1,\C^5)$.

By the same argument as in \ref{subsec:maog15}, 
$$\delta_{+1} = u \Theta_1 \wedge (-) + v \Theta_2 \wedge (-) + (z_0-\underline{z_0}) \Theta_3 \wedge (-) + (w_0 - \underline{w_0}) \Theta_4 \wedge (-).$$
Also, the previous argument is still valid to show that $\delta$ restricted to $\deg \leq 2$-component is given by
$$\delta|_{\deg \leq 2} = \delta_{+1} + \sum_{i=1}^4 f_i \frac{\partial}{\partial \Theta_i}$$
for some $f_i$ satisfying 
$$ u f_1 + v f_2 + (z_0 - \underline{z_0}) f_3 + (w_0 - \underline{w_0}) f_4 = W - W(\underline{b}),$$
and that $\delta_{-3}|_{\deg 3} =0$.

We next consider the action of $\delta$ on the $\deg=3$ component. We set  
$$\delta_{-1} (\Theta_1 \wedge \Theta_2 \wedge \Theta_3) = g_3 \Theta_1 \wedge \Theta_2 - g_2 \Theta_1 \wedge \Theta_3 + g_1 \Theta_2 \wedge \Theta_3.$$ 
It does not contain other $\Theta_i \wedge \Theta_j$'s due to the same reason as in the case of $\OG(1,\C^5)$. 
Since $\left( \delta_{-1} \right)^2 (\Theta_1 \wedge \Theta_2 \wedge \Theta_3) = 0$. We have
$$ -g_3 f_2  + g_2 f_3 =0, \,\, - g_3 f_1 + g_1 f_3 = 0, \,\, g_2 f_3 - g_3 f_2 = 0.$$
Thus $g_1 = c f_1, g_2= c f_2, g_3 = c f_3$ for some constant $c$. Likewise we have
$$\delta_{-1} (\Theta_1 \wedge \Theta_2 \wedge \Theta_4) = df_1 \Theta_2 \wedge \Theta_4 -df_2 \Theta_1 \wedge \Theta_4 + df_4 \Theta_1 \wedge \Theta_2 $$
$$\delta_{-1} (\Theta_1 \wedge \Theta_3 \wedge \Theta_4) =  k f_1 \Theta_3 \wedge \Theta_4 - k f_3 \Theta_1  \wedge \Theta_4 +k f_4 \Theta_1 \wedge \Theta_3 $$
$$\delta_{-1} (\Theta_2 \wedge \Theta_3 \wedge \Theta_4) = l f_2 \Theta_3 \wedge \Theta_4 - l f_3 \Theta_2 \wedge \Theta_4  + l f_4 \Theta_2 \wedge \Theta_3$$
Now we use the equation $\delta^2 (\Theta_i \wedge \Theta_j ) = (W - W(\underline{b})) \cdot \Theta_i \wedge \Theta_j$, and get the following system of linear equations.
\begin{equation*}
\begin{array}{l}
(c-1) (z_0 - \underline{z_0} f_3 + (d-1) (w_0 - \underline{w_0} ) f_4 = 0, \\
(c-1) v f_2 + (k-1) (w_0 - \underline{w_0}) f_4 =0, \\
(d-1) v f_2 + (k-1) (z_0 - \underline{z_0}) f_3 = 0, \\
(c-1) u f_1 + (l-1) (w_0 - \underline{w_0}) f_4 =0,  \\
(d-1) u f_1 + (l-1) (z_0 - \underline{z_0}) f_3 = 0, \\
(k-1) u f_1 + (l-1) v f_2 =0,
\end{array}
\end{equation*}
and it is tedious, but elementary to check that the only possible solution is $c=d=k=l=1$, which proves 
$$\delta|_{\deg \leq 3} = \delta_{+1} + \sum_{i=1}^4 f_i \frac{\partial}{\partial \Theta_i}.$$

It only remains to analyze $\delta$ on the $\deg=4$ component. We set
$$\delta (\Theta_1 \wedge \Theta_2 \wedge \Theta_3 \wedge \Theta_4) = \sum_{i=1}^4 \xi_i \Theta_i + \eta_1  \Theta_2 \wedge \Theta_3 \wedge \Theta_4- \eta_2 \Theta_1  \wedge \Theta_3 \wedge \Theta_4 + \eta_3 \Theta_1 \wedge \Theta_2 \wedge \Theta_4 - \eta_4 \Theta_1 \wedge \Theta_2 \wedge \Theta_3,$$
and compute $\delta (\Theta_1 \wedge \Theta_2 \wedge \Theta_3)$ to get
$$\delta(\Theta_1 \wedge \Theta_2 \wedge \Theta_3) = -(w_0 - \underline{w_0}) \Theta_1 \wedge \Theta_2 \wedge \Theta_3 \wedge \Theta_4 + f_3 \Theta_1 \wedge \Theta_2 + f_1 \Theta_2 \wedge \Theta_3 - f_2 \Theta_1 \wedge \Theta_3.$$
Since $\delta^2 = (W - W(\underline{b})) \cdot id)$, any $\Theta_i$-component of $\delta^2 (\Theta_1 \wedge \Theta_2 \wedge \Theta_3)$ should vanish. Since such a component from $\delta_{-1} (f_3 \Theta_1 \wedge \Theta_2 + f_1 \Theta_2 \wedge \Theta_3 - f_2 \Theta_1 \wedge \Theta_3) =0$,  this implies that $\xi_i=0$.
By looking at the coefficient of $\Theta_1 \wedge \Theta_2 \wedge \Theta_3$ in $\delta^2 (\Theta_1 \wedge \Theta_2 \wedge \Theta_3)$, we have $\eta_4 = f_4$. Similar argument show that $\eta_i = f_i$ for $i=1,2,3$.
We conclude that the matrix factorization $P_{\underline{b}}$ mirror to $(\scr{L}_0, \underline{b})$ is of wedge-contraction type. The rest of argument is the same as in \ref{subsec:maog15}.

Similarly as in the case of $\mathrm{OG}(1,\C^5)$, eigenspaces of ${\rm QH} (\mathrm{Gr}(2,\C^4))$ are 1-dimensional except the zero eigenspace (see Lemma \ref{lemma_critcal}). Therefore the argument appearing at the end of \ref{subsec:maog15} now proves the equivalence between the full subcategory of the Fukaya category of $\mathrm{Gr}(2,\C^4)$ consisting of Lagrangians with potential value $\lambda$ and $\oplus_\lambda \MF(W-\lambda)$, where $\lambda$ runs over all nonzero critical values of $W$.

\bibliographystyle{amsalpha}
\bibliography{geometry}

\end{document}